\documentclass{amsart}
\usepackage{latexsym,amssymb}
\usepackage{mathptmx}
\usepackage[dvips]{graphicx}
\usepackage{subfigure}
\usepackage{pstricks,pst-plot}

 \DeclareMathAlphabet{\mathcal}{OMS}{cmsy}{m}{n}

\makeatletter
 \@addtoreset{equation}{section}

 \@addtoreset{figure}{section}
 
\def\fps@figure{h,t}

\makeatother

\newtheorem{theorem}{Theorem}[section]
\newtheorem{proposition}[theorem]{Proposition}

\newtheorem{corollary}[theorem]{Corollary}
\theoremstyle{definition}

\newtheorem{remark}[theorem]{Remark}

\newcommand{\R}{\mathbb{R}}

\newcommand{\CC}{\mathbf{C}}
\newcommand{\DD}{\mathbf{D}}
\newcommand{\N}{\mathbb{N}}
\newcommand{\NN}{\mathcal{N}}
\newcommand{\<}{\left<}
\newcommand{\m}{\right>}
\newcommand{\p}{\mathcal{P}}
\newcommand{\J}{\mathbf{J}}
\newcommand{\Fix}{\mathop{\rm Fix}}
\newcommand{\Ker}{\mathop{\rm Ker}}
\newcommand{\Coad}{\mathop{\rm Coad}\nolimits}
\newcommand{\CoAd}{\mathop{\rm Coad}\nolimits}

\renewcommand{\Re}{\mathop{\rm Re}\nolimits}
\newcommand{\re}{\textsc{re}}
\newcommand{\diag}{\mathop{\rm diag}\nolimits}
\newcommand{\half}{\textstyle\frac12}
\newcommand{\sfrac}[2]{{\textstyle\frac{#1}{#2}}}
\newcommand{\uu}{\mathbf{u}}
\newcommand{\Z}{\mathbb{Z}}
\newcommand{\dd}{\mathsf{d}}
\newcommand{\ii}{{\mathsf i}}
\newcommand{\ee}{{\mathsf e}}

\newcommand{\so}{\mathfrak{so}(3)}
\newcommand{\SO}{\mathsf{SO}}
\newcommand{\OO}{\mathsf{O}}

\def\restr#1{\vrule height1.2ex width.4pt
               depth1.4ex\lower1.0ex\hbox{\scriptsize $\,#1$}}

\def\paragraph#1{\par\medskip\noindent{\it #1\hskip1em}}
\definecolor{gold}{rgb}{0.8,0.5,0.2}
\definecolor{LyapounovColour}{rgb}{0.9,0.3,0.3}
\definecolor{EllipticColour}{rgb}{0.4,0.6,0.9}
\definecolor{UnstableColour}{rgb}{0.6,0.5,0}

\newgray{midgray}{0.75}
\newgray{darkgray}{0.25}

\title[Stability of point vortices on the sphere]{Point Vortices on the Sphere:  \\ Stability of Symmetric Relative Equilibria} 
\author[F.~Laurent-Polz, J.~Montaldi \& M.~Roberts]{}


 \email{laurentpolzfj@voila.fr}
 \email{j.montaldi@manchester.ac.uk}
 \email{m.roberts@surrey.ac.uk}

\begin{document}

\maketitle

%AUTHORS
\centerline{\scshape Frederic Laurent-Polz}
\smallskip
{\footnotesize
 \centerline{Institut Non Lin\'eaire de Nice} 
\centerline{1361 route des Lucioles} 
\centerline{06560 Valbonne, France}
}

\bigskip

\centerline{\scshape James Montaldi}
\smallskip
{\footnotesize
 \centerline{School of Mathematics}
   \centerline{University of Manchester}
   \centerline{Manchester, M13 9PL, UK}
}

\bigskip

\centerline{\scshape Mark Roberts}
\smallskip
{\footnotesize
\centerline{Department of Mathematics} 
\centerline{University of Surrey}
\centerline{Guildford GU2 7XH, UK}
}

\bigskip

\centerline{\emph{Dedicated to Tudor Ratiu on the occasion of his 60th birthday} }

\bigskip

%%%%%%%%%%%%%%%%%%%%%%%%%%%%%%%%%%%%%%%%%%%%%%%%%%%%%%%%%%%%%%%%%%%%%%%%%%%%%

\begin{abstract}
We describe the linear and nonlinear stability and instability of
certain symmetric configurations of point vortices on the sphere forming
relative equilibria.  These configurations consist of one or two
rings, and a ring with one or two polar vortices.  Such configurations have
dihedral symmetry, and the symmetry is used to block diagonalize the relevant matrices, to distinguish the subspaces on which their eigenvalues need to be calculated, and also to 
describe the bifurcations that occur as eigenvalues pass through zero.
\end{abstract}

%%%%%%%%%%%%%%%%%%%%%%%%%%%%%%%%%%%%%%%%%%%%%%%%%%%%%%%%%%%%%%%%%%%%%%%%%%%%%

\tableofcontents

\section{Introduction}   \label{sec:Introduction} %

Since the work of Helmholtz \cite{H}, and later Kirchhoff \cite{K}, 
systems of point vortices on
the plane have been studied as finite-dimensional approximations to
vorticity dynamics in ideal fluids.  For a general survey of
patterns of point vortices see Aref \emph{et al} \cite{ANSTV03}.
Point vortex systems on the sphere, introduced by Bogomolov
\cite{B77}, provide simple models for the dynamics of concentrated
regions of vorticity, such as cyclones and hurricanes, in planetary
atmospheres.  In this paper we consider a non-rotating sphere, since
the rotation of the sphere induces a non-uniform background
vorticity which makes the whole system infinite-dimensional.

As in the planar case, the equations governing the motion of $N$
point vortices on a sphere are Hamiltonian \cite{B77} and this
property has been used to study them from a number of different
viewpoints.  Phase space reduction shows that the three vortex
problem is completely integrable on both the plane and the sphere:
the motion of three vortices of arbitrary vorticity on a sphere is
studied in \cite{KN98}.  The stability of some of the relative
equilibria described in \cite{KN98} are computed in \cite{PM98} and
numerical simulations are presented in \cite{MPS99}.  The existence
of relative equilibria of $N$ vortices is treated in \cite{LMR01},
and the nonlinear stability of a latitudinal ring of $N$ identical
vorticities is computed in \cite{BC03}, and independently in the
present paper. In fact the linear stability results of such a ring
obtained in \cite{PD93} coincide with the Lyapounov stability
results. The stability of a ring of vortices on the sphere together
with a central polar vortex is studied in \cite{CMS03}, and again
independently in the present paper, though with different methods
(and stronger results). The existence and nonlinear stability of
relative equilibria of $N$ vortices of vorticity $+1$ together with
$N$ vortices of vorticity $-1$ are studied in \cite{LP02}. It has
also been proved in \cite{LP} that relative equilibria formed of
latitudinal rings of identical vortices for the non-rotating sphere
persist to relative equilibria when the sphere rotates. However, the
question of stability becomes much more delicate: for motions that
are not relative equilibria, the vorticity of a point vortex is no
longer preserved as it interacts with the background vorticity, and
the problem becomes fundamentally infinite-dimensional.  In
\cite{Ku04} Kurakin studies the stability of equilibrium
configurations of identical vortices placed at the vertices of
regular polyhedra; he finds that the tetrahedron, octahedron and
icosahedron are stable, while the other two are unstable.  Finally,
studies of periodic orbits of point vortices on the sphere can be
found in \cite{ST02,To,LPf, LP04}.

Our study of the stability of relative equilibria is based on the
symmetries of the system, and especially the isotropy subgroups of
the relative equilibria.  The Hamiltonian is invariant under
rotations and reflexions of the sphere and permutations of
identical vortices.  However, some of these symmetries (eg
reflexions) are not symmetries of the equations of motion: they
are \emph{time-reversing} symmetries.  From Noether's theorem, the
rotational symmetry provides three conserved quantities, the
components of the momentum map $\Phi: \p \to \R^3$ where $\p$ is
the phase space.

Relative equilibria (\re) are dynamical trajectories that are generated by
the action of a 1-parameter subgroup of the symmetry group. More
intuitively, they correspond here to motions of the point vortices
which are stationary in a steadily rotating frame.  In other words,
the motion of a relative equilibrium corresponds to a rigid rotation
of $N$ point vortices about some axis (which we always take to be
the $z$-axis). In the same way as equilibria are critical points of
the Hamiltonian $H$, relative equilibria are critical points of the
restrictions of $H$ to the level sets of $\Phi$. Section \ref{intro}
is devoted to a description of the system of point vortices on the
sphere, and to an outline of stability theory for relative
equilibria.  The appropriate concept of stability for relative
equilibria of Hamiltonian system is Lyapounov stability \emph{modulo
a subgroup}.  The stability study is realized using on one hand the
\emph{energy-momentum method} \cite{Pa92,Or98} which consists of
computing the eigenvalues of a certain Hessian, and leads to
\emph{nonlinear} stability results, and on the other hand a linear
study computing the eigenvalues of the linearization of the
equations of motion. To both these ends, we block diagonalize these
matrices using a suitable basis, the \emph{symmetry adapted basis}
(Section \ref{sec:basis}), which makes use of the specific dihedral
symmetry of the relative equilibrium. This is equivalent to noting
that the matrices (or certain submatrices) are \emph{circulant}, as
used in \cite{CMS03}.  The symmetry adapted basis splits the
tangent space to the phase space into a number of modes, from
$\ell=0$ to $\ell=[n/2]$ (where $n$ is the number of vortices in the
ring), and calculations can proceed separately for each mode.
Moreover, the symmetry is also used to apply the energy-momentum
method as it helps distinguish on which subspaces computations are
needed. In Section \ref{sec:symmetry+bifurcations} we outline the 
different bifurcations which occur in this system: symmetric pitchforks, 
Hamiltonian Hopf bifurcations and those bifurcations which occur due 
to the `geometry of reduction'---those occurring in a neighbourhood of points 
with zero momentum. 

The remaining five sections each treat one of five different types of relative equilibrium, consisting of rings of identical vortices together with possible vortices at the poles, whose existence were proved in \cite{LMR01}. The notation for the different configurations is taken from the same source and is described at the end of the introduction.  We now outline the main stability results.

We begin in Section \ref{sec:Cnv(R)} by computing the stability of the relative equilibria consisting of a single ring of identical vortices, a configuration denoted $\CC_{nv}(R)$ (Figure \ref{fig:ring,pole}(a)).  We show in Theorem~\ref{polvani} that for $n\geq 7$, they are unstable for all co-latitudes\footnote{the co-latitude $\theta$ of a point on the sphere is the angle subtended with the North pole} of the ring, while for $n<7$ there exist ranges of Lyapounov stability when the ring is near a pole. These results are not new \cite{PD93} (for linear stability results) and \cite{BC03} (for nonlinear stability), but serve to demonstrate the method used in later sections.

In Section \ref{sec:Rp-stability}, we study the stability of the relative equilibria $\CC_{nv}(R,p)$ (Figure~\ref{fig:ring,pole}(b)) which are configurations formed of a ring of $n$ identical vortices together with a single polar vortex.  For $n\geq 7$ they are all unstable if the vorticity $\kappa$ of the polar vortex has opposite sign to that of the ring.  However if the vorticities have the same sign then for each co-latitude of the ring there exists a range of $\kappa$ for which the relative equilibrium is Lyapounov stable. Adding polar vortices can therefore \emph{stabilize} the unstable pure ring relative equilibria. The detailed results are contained in Theorem~\ref{ring+pole/stab}, its corollary and the following discussion. Our results are consistent with those of \cite{CMS03} (aside from an error in their Figure~7 where the wrong curves are plotted), though the present methods are stronger as they give more regions of stability than obtained in \cite{CMS03}---see Remark \ref{rmk:CMS}.

%%%%%%%%%%%%%%%%%%%%%%%%%%%%%%%
\begin{figure}[tbp]
\begin{center}\psset{unit=0.85}
\begin{pspicture}(-3,-3)(3.5,3)
 \pscircle(0,0){2.5}
 \psline[linestyle=dashed](0,-2.7)(0,2.2)
 \psline{->}(0,2.236067978)(0,3.2)
 \psellipse[linewidth=0.1pt](0,2.8)(0.36,0.12)
 \psline[arrowsize=.2]{->}(0,2.68)(0.1,2.68)
 \psellipse[linestyle=dotted](0,1.118033989)(2.165063510,0.9682458370)
 \psdot[dotstyle=x](0,0)
 \psdots[dotsize=.2](1.786904023, 1.664746714)(-1.222486816, 1.917161762)(-1.786904023, .5713212643)(1.222486817, .3189062162)
 \rput(1.7,1.3){\small +1}
 \rput(-1,1.7){\small +1}
 \rput(-1.6,.3){\small +1}
 \rput(1.5,.1){\small +1}
 \psdots[dotsize=.1](0,-2.236067978)(0,2.236067978) % poles
\end{pspicture}
\begin{pspicture}(-3.5,-3)(3,3)
 \pscircle(0,0){2.5}
 \psline[linestyle=dashed](0,-2.7)(0,2.5)
 \psline{->}(0,2.236067978)(0,3.2)
 \psellipse[linewidth=0.1pt](0,2.8)(0.36,0.12)
 \psline[arrowsize=.2]{->}(0,2.68)(0.1,2.68)
 \psellipse[linestyle=dotted](0,1.118033989)(2.165063510,0.9682458370)
 \psdot[dotstyle=x](0,0)
 \psdots[dotsize=.2](1.786904023, 1.664746714)(-1.222486816, 1.917161762)(-1.786904023, .5713212643)(1.222486817, .3189062162)
 \psdot[dotsize=.3](0,2.236067978)
 \psdots[dotsize=.1](0,-2.236067978)(0,2.236067978) % poles
 \rput(1.7,1.3){\small +1}
 \rput(-1,1.7){\small +1}
 \rput(-1.6,.3){\small +1}
 \rput(1.5,.1){\small +1}
 \rput(.3,2.1){\small $\kappa$}
\end{pspicture}
\caption{The $\CC_{nv}(R)$ and $\CC_{nv}(R,p)$ relative equilibria, here shown for $n=4$. The $\CC_{nv}(R,2p)$ configurations have a vortex at the South pole as well.}
\label{fig:ring,pole}
\end{center}
\end{figure}
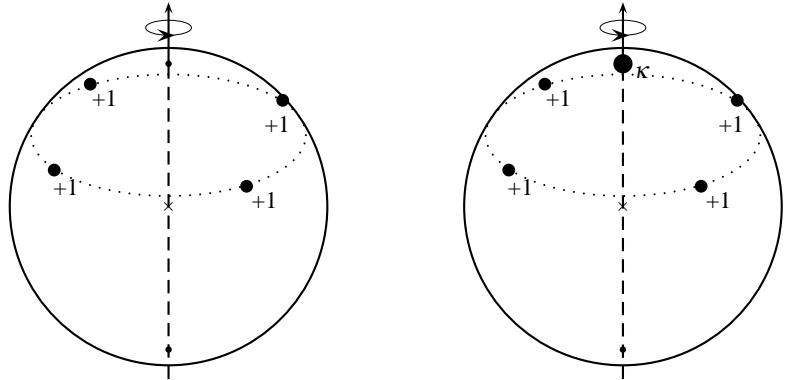
%%%%%%%%%%%%%%%%%%%%%%%%%%%%%%%%%%%

\begin{figure}[tbp]
\begin{center}\psset{unit=0.85}
\begin{pspicture}(-3,-3)(3.5,3)
 \pscircle(0,0){2.5}
 \psline[linestyle=dashed](0,-2.7)(0,2.2)
 \psline{->}(0,2.236067978)(0,3.2)
 \psellipse[linestyle=dotted](0,1.118033989)(2.165063510,0.9682458370)
 \psdot[dotstyle=x](0,0)
 \psdots[dotsize=.2](1.786904023, 1.664746714)(-1.222486816, 1.917161762)%
   (-1.786904023, .5713212643)(1.222486817, .3189062162)
 \psline[linewidth=0.1pt](1.786904023, 1.664746714)(-1.786904023, .5713212643)
 \psline[linewidth=0.1pt](-1.222486816, 1.917161762)(1.222486817, .3189062162)
 \rput(1.7,1.3){\small +1}
 \rput(-1,1.7){\small +1}
 \rput(-1.6,.3){\small +1}
 \rput(1.5,.1){\small +1}
 \psellipse[linestyle=dotted](0,-1.581138830)(1.767766953,0.7905694152)
 \psdots[dotsize=.25](1.459001025, -1.134749760)(-.9981563053, -.9286537357)%
 (-1.459001025, -2.027527900)(.9981563060,-2.233623924)
 \psline[linewidth=0.1pt](1.459001025, -1.134749760)(-1.459001025, -2.027527900)
 \psline[linewidth=0.1pt](-.9981563053, -.9286537357)(.9981563060,-2.233623924)
 \rput(1.3,-1.3){\small $\kappa$}
 \rput(-0.8,-1.2){\small $\kappa$}
 \rput(-1.2,-1.8){\small $\kappa$}
 \rput(1,-2){\small $\kappa$}
 \psdots[dotsize=.1](0,-2.236067978)(0,2.236067978) % poles
\end{pspicture}
\begin{pspicture}(-3.5,-3)(3,3)
 \pscircle(0,0){2.5}
 \psline[linestyle=dashed](0,-2.7)(0,2.5)
 \psline{->}(0,2.236067978)(0,3.2)
 \psellipse[linestyle=dotted](0,1.118033989)(2.165063510,0.9682458370)
 \psdot[dotstyle=x](0,0)
 \psdots[dotsize=.2](1.786904023, 1.664746714)(-1.222486816, 1.917161762)(-1.786904023, .5713212643)(1.222486817, .3189062162)
 \psline[linewidth=0.1pt](1.786904023, 1.664746714)(-1.786904023, .5713212643)
 \psline[linewidth=0.1pt](-1.222486816, 1.917161762)(1.222486817, .3189062162)
 \rput(1.7,1.3){\small +1}
 \rput(-1,1.7){\small +1}
 \rput(-1.6,.3){\small +1}
 \rput(1.5,.1){\small +1}
 \psellipse[linestyle=dotted](0,-1.581138830)(1.767766953,0.7905694152)
 \psdots[dotsize=.25](.3258664255, -.8041174565)(-1.737472611,-1.435406934)%
 (-.3258664248, -2.358160204)(1.737472611, -1.726870725)
 \psline[linewidth=0.1pt](.3258664255, -.8041174565)(-.3258664248, -2.358160204)
 \psline[linewidth=0.1pt](-1.737472611,-1.435406934)(1.737472611, -1.726870725)
 \rput(2.0,-1.8){\small $\kappa$}
 \rput(0.6,-1){\small $\kappa$}
 \rput(-1.4,-1.6){\small $\kappa$}
 \rput(-0.4,-2.7){\small $\kappa$}
 \psdots[dotsize=.1](0,-2.236067978)(0,2.236067978) % poles
\end{pspicture}
\caption{Configurations of types $\CC_{nv}(2R)$ (2 aligned rings)
and $\CC_{nv}(R,R')$ (2 staggered rings).  As in the previous figure,  this is depicted for $n=4$.} \label{fig:2rings}
\end{center}
\end{figure}
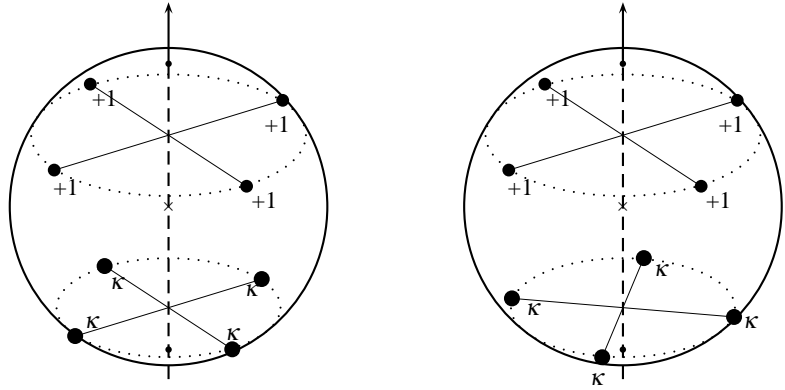

In Sections \ref{sec:2R} and \ref{sec:RR'} we investigate
configurations formed of two rings of arbitrary vorticities (each
ring, as always, consisting of identical vortices).  In \cite{LMR01}
it was shown that two rings of $n$ vortices can be relative
equilibria if and only if they are either aligned or staggered.
These two arrangements are denoted $\CC_{nv}(2R)$ and
$\CC_{nv}(R,R')$ respectively (see Figure~\ref{fig:2rings}). Here we
show that for almost all pairs of ring latitudes there is a unique
ratio of the ring vorticities for which these configurations are
relative equilibria. Numerical computations of their stability
suggest that these relative equilibria can only be stable if $n\leq
6$, and in the aligned case the two rings must be close to opposite
poles, and hence have opposite vorticities. In some cases, staggered
rings may also be stable when in the same hemisphere.

Finally, in Section \ref{sec:R2p} we study the relative equilibria $\CC_{nv}(R,2p)$ which are configurations formed of a ring of $n$ identical vortices together with \emph{two} polar vortices. The two polar vorticities make this a 2-parameter family of systems.  We obtain analytic (in)stability criteria, but only with respect to certain modes ($\ell\geq2$). As in the case of a single polar vortex, the two polar vortices play the role of control parameters for the stability.  The details are contained in Theorem \ref{thm:R2p-stability}.  We then continue with numerical investigations for the remaining mode ($\ell=1$) in order to provide full stability criteria; the conclusions are shown in a series of stability diagrams. Further diagrams are available on the second author's website \cite{Mo-web}.

The numerical computations all consist simply of computing the eigenvalues of the Hessian of the Hamiltonian and the linear vector field on the symplectic slice, using the block diagonalization of the matrices to simplify the calculations. These are precisely the same calculations as in the earlier sections where they can be performed analytically. The numerics are all done using \textsc{Maple}, and a selection of the code can also be downloaded from \cite{Mo-web}.

%%%%%%%%%%%%%%%%%%%%%%%%%%%%%%%
%%%%%%%%%%%%%%%%%%%%%%%%%%%%%%%%%%%

In principle the method applies to larger numbers of rings but the
algebraic problem of diagonalizing the matrices in general becomes
intractable; however numerical studies for particular (numerical)
values of the vorticities in the rings would be feasible.

\paragraph{Symmetry group notation}
All possible symmetry types of configurations of point vortices on
the sphere were classified in \cite{LMR01}.  The symmetry group of
the system is of the form $\OO(3)\times S$, where $S$ is a group of
permutations, and a particular configuration with symmetry, or
isotropy, subgroup $\Sigma<\OO(3)\times S$ is denoted $\Gamma(A)$,
where $\Gamma$ is the projection of $\Sigma$ to $\OO(3)$ and $A$
represents the way $\Sigma$ permutes the point vortices, so describing 
their geometry (for example $A=R$ means they are in a ring, $A=p$ 
means a polar vortex).  The classical \emph{Sch\"{o}nflies-Eyring} 
notation for subgroups of $\OO(3)$ is used.

In this paper we single out configurations consisting of concentric
rings of identical vortices, with the same number of vortices in
each ring, and with possible polar vortices. These configurations
have cyclic symmetry (in the `horizontal plane'), and the
Sch\"{o}nflies-Eyring notation for this subgroup of $\OO(3)$ is
$\CC_n<\SO(3)$.  In fact we only consider the cases where the rings
are either aligned (the vortices lie on the same longitudes) or
staggered (they lie on intermediate longitudes, out of phase by
$\pi/n$). In these two cases the symmetry group is the larger
\emph{dihedral} group $\CC_{nv}$ ($n$ being the number of vortices
in each ring, and $v$ denoting the fact that there are vertical
planes of reflexion). For such configurations, we write
$\CC_{nv}(k_1R,k_2R',k_pp)$ to mean that there are $k_r=k_1+k_2$
rings and $k_p$ polar vortices. The difference between $R$ and $R'$
is that the $k_1$ rings R are all aligned and the $k_2$ rings $R'$
are staggered with respect to the first (and so aligned with each
other). Of course $k_p=0,1$ or 2.

In particular, the following table shows the five types of configuration we consider in detail and the relevant section in this paper:

\medskip

\begin{center}
\begin{tabular}{rll}
\hline
\rule{0pt}{12pt}$\CC_{nv}(R)$ & a single ring of $n$ identical vortices & \S\,5\\
$\CC_{nv}(R,p)$ & a single ring as above, together with a polar vortex  & \S\,6\\
$\CC_{nv}(2R)$ &  a pair of aligned rings of $n$ identical vortices& \S\,7\\
%, on different latitudes\\
$\CC_{nv}(R,R')$ &  a pair of staggered rings of $n$ identical vortices& \S\,8 \\
$\CC_{nv}(R,2p)$ & a single ring as above, together with two polar vortices& \S\,9 \\ 
\hline
\end{tabular}
\end{center}

\medskip

We also refer on occasion to other symmetry types. In particular $\DD_{nh}$ is the subgroup of $\OO(3)$ generated by the dihedral group $\CC_{nv}$ together with a reflexion in the horizontal plane: it is the symmetry group of a pair of aligned rings of vortices lying symmetrically about the equator. Similarly, $\DD_{nd}$ is the subgroup generated by $\CC_{nv}$ together with a combined reflexion-rotation: reflexion about the horizontal plane combined with a rotation by $\pi/n$ (the rotations in $\CC_{nv}$ are through multiples of $2\pi/n$). This is the symmetry group of a pair of staggered rings of vortices lying at opposite latitudes.

For the symmetry arguments to hold, the vortices in a single ring
must be identical, but we do not assume that the vortices in
different rings or at the poles are identical.

\paragraph{Historical note} An early version of this paper was written in 2002, with most of the stability calculations performed by the first author (FLP). It was submitted soon after, but  the editor asked for more extensive numerics.  This lay in the to-do pile of the other two authors for the intervening years (FLP having left academic research in the meantime), and the occasion of Tudor Ratiu's birthday celebration was the impetus required to produce a final version, with the more extensive numerics---especially those in the final section. 

%%%%%%%%%%%%%%%%%%%%%%%%%%%%%%%%%%%%%%%%%%%%%%%%%%%%%%%%%%%%%%%%%%%%%%%%%%

\section{Point vortices on the sphere and stability theory}
\label{intro}

In this section we briefly recall that the system of point vortices
on a sphere is an $n$-body Hamiltonian system with symmetry and we
review the stability theory for relative equilibria.

%%%%%%%%%%%%%%%%%%%%%%%%%%%%%%%%%%%%%

\subsection{Point vortices}

Consider $n$ point vortices $x_1,\dots,x_n\in S^2$ with vorticities
$\kappa_1,\dots,\kappa_n\in\R$.  Let $\theta_i,\phi_i$ be
respectively the co-latitude and the longitude of the vortex $x_i$
(so $\theta=0$ corresponds to the North pole). The dynamical system
is Hamiltonian with Hamiltonian given by
$$
H \ =\  -\sum_{i<j}\kappa_i\kappa_j\ln (
1-\cos\theta_i\cos\theta_j-\sin\theta_i\sin\theta_j\cos(\phi_i-\phi_j)
)
$$
and conjugate variables given by
$q_i=\sqrt{\vert\kappa_i\vert}\cos\theta_i$ and
$p_i=\mathrm{sign}(\kappa_i)\sqrt{\vert\kappa_i\vert}\phi_i$.
Note that many authors have a constant factor (involving $\pi$) in the Hamiltonian; this will not effect the existence or stability of the relative equilibria, but is equivalent to a rescaling of time by that factor.

The phase space is $\p = \lbrace (x_1,\dots,x_{n})\in
S^2\times\cdots\times S^2\mid x_i\neq x_j\ \mbox{if}\ i \neq j
\rbrace$ endowed with the symplectic form $\omega=\sum_{i}
\kappa_i \sin\theta_i\ d\theta_i\land d\phi_i$.  The Hamiltonian
vector field $X_H$ satisfies $\omega(\ \cdot\ ,X_H(x))=dH_x$.  If
we consider $S^2$ as a subset of $\R^3$, so each vortex $x_j\in
\R^3$, then we obtain
\begin{equation}\label{eq:ode}
\begin{array}{c}
\displaystyle \dot {x_i} \ =\  X_H(x)_i\ =\ \sum_{j,j\neq i} \kappa_j\,
 \frac{x_{j}\times x_{i}}{1-x_{i}\cdot x_{j}}\ ,\ i=1,\dots,N,\\[12pt]
\displaystyle H \ =\  - \sum_{i<j}\kappa_i\kappa_j\ln(\|x_i-x_j\|^2/2).
\end{array}
\end{equation}
It follows that $H$ is invariant under $\OO(3)$ acting by rotations and
reflexions of the sphere.

The rotation subgroup $\SO(3)$ leaves the symplectic form invariant,
and so the Hamiltonian vector field $X_H$ is $\SO(3)$-equivariant.
On the other hand, the reflexions in $\OO(3)$ reverse the sign of
the symplectic form and so give rise to time-reversing symmetries of
$X_H$. Moreover $H$, $\omega$ and $X_H$ are all invariant (or
equivariant) with respect to permutations of vortices with equal
vorticity.

The rotational symmetry implies the existence of a momentum map
$\Phi:\p\to\so^*\simeq\R^3$:
\begin{equation}\label{eq:momentum map}
 \Phi(x)\ =\ \sum_{j=1}^{N}\kappa_j x_j\qquad (x_j\in S^2\subset\R^3)
\end{equation}
which is conserved under the dynamics. In other words each of the
three components of $\Phi(x)$ is a conserved quantity. We will
always orient the solution we are studying so that $\Phi(x)$ lies on
the $z$-axis, and write $\Phi(x)=(0,0,\mu)$. There is no intrinsic difference between $\mu>0$ and $\mu<0$---they are rotationally equivalent.  However when studying bifurcations in a neighbourhood of a point with $\mu=0$, then the two sides should be distinguished.

%%%%%%%%%%%%%%%%%%%%%%%%%%%%%%%%%%

\subsection{Relative equilibria}

A point $x_e\in\p$ is a relative equilibrium if and only if there
exists $\xi\in\so\simeq\R^3$ (the angular velocity) such that
$x_e$ is a critical point of the function $H_{\xi}(x)=H(x)-\left<
\Phi(x),\xi\right>$, where the pairing $\left<\ ,\ \right>$
between $\R^3$ and its dual is identified with the canonical
scalar product on $\R^3$.  Equivalently, relative equilibria are
critical points of the restriction of $H$ to $\Phi^{-1}(\mu)$,
since the level set $\Phi^{-1}(\mu)$ are always non-singular for
point vortex systems of more than two vortices.  The function
$H_\xi$ is called the \emph{augmented Hamiltonian}.

Since the momentum is conserved, we can choose a frame for $\R^3$ such
that $\Phi$ is parallel to the $z$-axis (provided the momentum is
non-zero).  It follows from the symmetry that the angular velocity
$\xi\in\R^3$ is also parallel to the $z$-axis.  We can therefore
identify $\xi$ and $\Phi$ with their $z$-components and the
augmented Hamiltonian becomes simply $H_{\xi}(x)=H(x)-\xi\Phi(x)$.

Let $f:\p\to\R$ be a $K$-invariant function with $K$ a compact group.
Recall that $\Fix(K)=\lbrace x\in\mathcal{P} \mid g\cdot x = x,
\ \forall g\in K \rbrace$.  The \emph{Principle of Symmetric
  Criticality} \cite{P79} states that a critical point of the
restriction of a $K$-invariant function $f$ to $\Fix(K)$ is a
critical point of $f$.  As a corollary, if the Hamiltonian is
invariant under $K$ and $x_e$ is an isolated point in
$\Fix(K)\cap\Phi^{-1}(\mu)$, then $x_e$ is a relative equilibrium.
It follows in particular that all configurations of type
$\CC_{nv}(R)$, $\CC_{nv}(R,p)$ (Figure~\ref{fig:ring,pole}), and
$\CC_{nv}(R,2p)$ are relative equilibria: take $K$ such that
$\pi(K)=\CC_{nv}$ where $\pi:\OO(3)\times S_n\to \OO(3)$ is the
Cartesian projection.

Finally, one can show that if $x_e$ is a relative equilibrium with
angular velocity $\xi$, then $H_\xi$ is a $G_{x_e}$-invariant
function, where $G_{x_e}<\OO(3)\times S$ is the symmetry group
(isotropy group) of the configuration $x_e$.

%%%%%%%%%%%%%%%%%%%%%%%%%%%%%%

\subsection{Stability theory and isotypic decomposition}
\label{sec:isotypic}

Stability is determined by the energy-momentum method together with
an isotypic decomposition of the symplectic slice. We recall the
main points of the method.

Let $x_e\in\p$ be a relative equilibrium, $\mu=\Phi(x_e)$, and $\xi$ be
its angular velocity.  The energy-momentum method consists of
determining the symplectic slice
$$
\NN\ =\ (\so_{\mu}\cdot x_e)^\bot \cap \Ker D\Phi(x_e)
$$
transverse to $\so_{\mu}\cdot x_e$, where
$$
\SO(3)_{\mu}\ =\ \lbrace g\in \SO(3)\mid\CoAd_g\cdot \mu=\mu \rbrace
$$
and then examining the definiteness of the restriction
$\dd^2H_{\xi}|_{\NN}(x_e)$ of the Hessian $\dd^2H_\xi(x_e)$ to $\NN$.
(In practice we will represent $\mu$ as a vector, in which case
$\Coad_g\mu = g\mu$ is just matrix multiplication.) If $K$ is a
group acting on the phase space a relative equilibrium $x_e$ is said
to be \emph{Lyapounov stable modulo K} if for all $K$-invariant open
neighbourhoods $V$ of $K\cdot x_e$ there is an open neighbourhood
$U\subseteq V$ of $x_e$ which is invariant under the Hamiltonian
dynamics.  The \emph{energy-momentum theorem} of Patrick \cite{Pa92}
holds since $\SO(3)$ is compact, and so we have:
\begin{center} \textit{If $\dd^2H_{\xi}|_{\NN}(x_e)$ is definite, then
  $x_e$ is Lyapounov stable modulo $\SO(3)_{\mu}$.}
\end{center}
For $\mu\neq 0$, $\SO(3)_{\mu}$ is the set of rotations with axis
$\<\mu\m$, and so isomorphic to $\SO(2)$, while for $\mu = 0$,
$\SO(3)_{\mu} = \SO(3)$.  If $\mu\neq 0$ Lyapounov stability modulo
$\SO(3)_{\mu}$ of a relative equilibrium with non-zero angular
velocity coincides with ordinary (orbital) stability of the
corresponding periodic orbit.

The second ingredient consists of performing an isotypic
decomposition of the symplectic slice $\NN$ in order to block
diagonalize $\dd^2H_{\xi}|_{\NN}(x_e)$.  Let $V$ be a finite
dimensional representation of a compact Lie group $K$.  Recall that
a $K$ invariant subspace $W\subset V$ of $K$ is said to be
\emph{irreducible} if $W$ has no proper $K$ invariant subspaces.
Since $K$ is compact, $V$ can be expressed as a direct sum of
irreducible representations: $V=W_1\oplus\cdots\oplus W_n$. In
general this is not unique.  There are a finite number of
isomorphism classes of irreducible representations of $K$ in $V$,
say $U_1,\dots,U_\ell$.  Let $V_k$ ($k=1,\dots,\ell$) be the sum of
all those irreducible representations $W_j$ in the above sum for which 
$W_j$ is isomorphic to $U_k$.  Then $V=V_1\oplus\cdots\oplus
V_\ell$. This decomposition of $V$ is unique and is called the
\emph{$K$-isotypic decomposition} of $V$ \cite{Se78}.  By Schur's
Lemma, the matrix of a $K$-equivariant linear map $f:V\to V$ block
diagonalizes with respect to a basis
$\mathcal{B}=\{\mathcal{B}_1,\dots,\mathcal{B}_l \}$ where
$\mathcal{B}_k$ is a basis of $V_k$, each block corresponding to a
subspace $V_k$.  The basis $\mathcal{B}$ is called a \emph{symmetry
adapted basis}. In this paper, the isotropy group $K$ is a dihedral
group, and the corresponding symmetry adapted basis consists of
`Fourier modes', defined below.

Let $G$ denote the group of all symmetries of the Hamiltonian $H$ and vector field $X_H$ and $G^0$ the subgroup consisting of time-preserving symmetries, i.e.\ those acting symplectically. The elements of $G\setminus G^0$ reverse the symplectic form, and hence reverse the vector field $X_H$, so effectively reversing time.  In the case of $N$ identical vortices we have $G=\OO(3)\times S_N$ and $G^0=\SO(3)\times S_N$.  Since $H_\xi$ is a $G_{x_e}$-invariant function, the restricted Hessian $\dd^2H_{\xi}|_{\NN}(x_e)$ is $G_{x_e}$-equivariant as a matrix.  Moreover the symplectic slice $\NN$ is a $G_{x_e}$-invariant subspace and so we can implement a $G_{x_e}$-isotypic decomposition of $\NN$ to block diagonalize $\dd^2H_{\xi}|_{\NN}(x_e)$. This block diagonalization of $\dd^2H_{\xi}|_{\NN}(x_e)$ simplifies the computation of its eigenvalues, and hence of its definiteness. If it is definite then the relative equilibrium is Lyapounov stability modulo $\SO(3)_{\mu}$.

If $\dd^2H_{\xi}|_{\NN}(x_e)$ is {\it not} definite then we study the
spectral stability of $x_e$. In particular we examine the
eigenvalues of $L_\NN$, the matrix of the linearized system on the
symplectic slice, that is $L_\NN=\J_{\NN}\;\dd^2H_{\xi}|_{\NN}(x_e)$,
where $\J_{\NN}^{-1}$ is the matrix of $\omega|_{\NN}$.  The matrix
$L_\NN$ is $G_{x_e}^0$-equivariant and so we perform a
$G_{x_e}^0$-isotypic decomposition of $\NN$ to obtain a block
diagonalization of $L_\NN$, and so to determine the spectral
stability of $x_e$. In particular, if $L_\NN$ has eigenvalues with
non-zero real part, then $x_e$ is linearly unstable.  Note that the
block diagonalization of $\dd^2H_{\xi}|_\NN(x_e)$ refines that of
$L_\NN$ since $G_{x_e}^0\subset G_{x_e}$.

Throughout this paper, we will align the vortices so that $\phi(x_e)$ lies along the $z$-axis, and we write $\phi(x_e)=(0,0,\mu)$. If $\mu\neq0$ then \emph{Lyapounov stable} will mean \emph{Lyapounov stable modulo $\SO(2)$}, while if $\mu=0$ it means \emph{Lyapounov stable modulo $\SO(3)$}.

%%%%%%%%%%%%%%%%%%%%%%%%%%%%%%%%%%%%%%%%%%%%%%%%%%%%%%%%%%%%%%%%%%%%%%%%%%%

\section{Symmetry adapted bases for rings and poles}
\label{sec:basis}

In this section we give the ingredients needed to determine the symmetry adapted bases for the symplectic slice at the configurations described above, that is those of type $\CC_{nv}(k_1R,$ $k_2R', k_pp)$.  In the first subsection we give a general symmetry adapted basis for the tangent space $T_{x_e}\p$ to the phase space at such a configuration, and express the derivative of the momentum map and tangent space to the group orbit in this basis. In the following two subsections we describe the isotypic decomposition of $T_{x_e}\p$, first for a single ring and then in general. Recall that the isotropy subgroup $G_{x_e}$ is always a dihedral group $\CC_{nv}$ (which has order $2n$) and that the irreducible representations of this group are of dimension $1$ or $2$. The cyclic subgroup $\CC_n$ acts symplectically and is denoted
$G_{x_e}^0$.

The actual symmetry adapted bases of the symplectic slices will be given case-by-case in the following sections. We do not give the proof of the results in the first subsection, since they can be easily deduced from the proofs of Propositions 4.1--4.4 in \cite{LP02}.

%%%%%%%%%%%%%%%%%%%%%%%%%%%%%%%%%

\subsection{Description of the symplectic slice}

Let $x_e$ be a $\CC_{nv}(k_1R,k_2R', k_pp)$ configuration. Let $k_r=k_1+k_2$ be the total number of rings. The total number of vortices is then $N = nk_r + k_p$, where $k_p=0,1$ or $2$ is the number of polar vortices.  We suppose the vorticities in ring $j$ is $\kappa_j$ for $j = 1,\ldots,k_r$ while the vorticities of the possible polar vortices are $\kappa_N$ for the North pole and $\kappa_S$ for the South pole. In this paper we only consider in detail $k_r=1$ or $2$, but here we describe the more general case.

For each ring $j=1,\dots,k_r$ let $s=1,\dots,n$ label the vortices in the ring in cyclic order (in an easterly direction), and define \emph{real} tangent vectors $\alpha^{(\ell)}_{j,\theta}$ etc.\ in $T_{x_e}\p$ by
\begin{equation} \label{eqn:alphas,betas}
\begin{array}{rcl}
  \alpha^{(\ell)}_{j,\theta} + \ii\,\beta^{(\ell)}_{j,\theta}
  &=&
  \sum_{s=1}^n\ee^{\ii\ell\phi_{j,s}}\,\delta\theta_{j,s}
  \\[4pt]
  \alpha^{(\ell)}_{j,\phi} + \ii\,\beta^{(\ell)}_{j,\phi}
  &=&
  \sum_{s=1}^n\ee^{\ii\ell\phi_{j,s}} \,\delta\phi_{j,s}
\end{array}
\end{equation}
where $\ell=0,\dots,[n/2]$, $\ii=\sqrt{-1}$ and $\phi_{j,s} = 2\pi s/n-\phi_j^0$, which is the longitude of the $s^\mathrm{th}$ vortex in ring $j$, where $\phi_j^0=0$ or $\pi/n$ depending on whether the $j^{\rm th}$ ring of vortices is of type $R$ or $R'$.  Note that $\beta^{(\ell)}_{j,\theta},\beta^{(\ell)}_{j,\phi}$ vanish for $\ell=0$ and $n/2$ (for $n$ even).  The vectors $\alpha^{(\ell)}_{j,\theta} $ etc.\ defined in (\ref{eqn:alphas,betas}) are called the \emph{Fourier modes}, and are depicted (for $n=4$ and $5$ and rings of type $R$) in Figure~\ref{fig:Fourier} at the end of the paper.

For each pole $i=1,\dots,k_p$ (where $k_p=0,1,2$) we also have tangent vectors $\delta x_i$ and $\delta y_i$.

The tangent vectors defined in the last two paragraphs are almost canonical, in the sense that
\begin{equation} \label{eq:omega on alphas,betas}
\begin{array}{rcl} 
\omega\left(\alpha^{(\ell)}_{j,\theta},\, \alpha^{(\ell)}_{j,\phi}\right) &=&
  \begin{cases}n\sin\theta_j\,\kappa_j & \text{if $\ell=0,\, n/2$} \\ 
        \frac12{n\,\sin\theta_j}\,\kappa_j & \text{otherwise} \end{cases} \\[8pt]
\omega\left(\beta^{(\ell)}_{j,\theta},\, \beta^{(\ell)}_{j,\phi}\right) &=&
   \frac12{n\sin\theta_j}\,\kappa_j \\[6pt]
\omega\left(\delta x_i,\,\delta y_i \right)&=&{\rm
    sign}(z_i)\,\kappa_i,
\end{array}
\end{equation}
while all other pairings vanish.

The cyclic subgroup $\CC_n$ of the isotropy subgroup $\CC_{nv}< \OO(3)\times S$ (see the end of the Introduction) acts symplectically, and in terms of the Fourier modes, the action of the generator $\sigma\in \CC_n$ is
\begin{equation} \label{eq:action on alphas,betas}
\begin{array}{rcl}
  \sigma\cdot(\alpha^{(\ell)}_{j,\theta} + \ii\beta^{(\ell)}_{j,\theta}) %
  &=&
  \exp\left(\frac{2\pi\ii\ell}{n}\right)(\alpha^{(\ell)}_{j,\theta} + \ii\beta^{(\ell)}_{j,\theta}) %
  \\[4pt]
  \sigma\cdot(\alpha^{(\ell)}_{j,\phi} + \ii\beta^{(\ell)}_{j,\phi})
  &=&
  \exp\left(\frac{2\pi\ii\ell}{n}\right)(\alpha^{(\ell)}_{j,\phi} + \ii\beta^{(\ell)}_{j,\phi})\\[4pt]
  \sigma\cdot(\delta x_j+\ii\delta y_j)
  &=&
  \exp\left(\frac{2\pi\ii}{n}\right)(\delta x_j+\ii\delta y_j).
\end{array}
\end{equation}

The reflexions in $\CC_{nv}\simeq \DD_n$ (dihedral group) act antisymplectically. One of the reflexions $\kappa\in\DD_n\setminus\CC_n$ acts on the sphere by $(\theta,\,\phi)\mapsto(\theta,\,-\phi)$ and by permuting the particles within each ring by,
$$\kappa\cdot(1,\dots,n) = \begin{cases}(n-1,\dots,1,n) & \text{for rings of type $R$},\cr
  (n,\dots,2,1) & \text{for rings of type $R'$}.\end{cases}
$$
The difference arises because for rings of type $R$, vortex $n$ lies in the half-plane $\phi=0$, while for rings of type $R'$ it lies in the half-plane $\phi=-\pi/n$. 
In terms of the Fourier modes, this reflexion acts by, 
$$\begin{array}{lll}
\alpha^{(\ell)}_{j,\theta} \mapsto \alpha^{(\ell)}_{j,\theta}\,, && 
    \beta^{(\ell)}_{j,\theta} \mapsto -\beta^{(\ell)}_{j,\theta}\,,\\[6pt]
\beta^{(\ell)}_{j, \phi}\mapsto \beta^{(\ell)}_{j, \phi}\, , &&  
    \alpha^{(\ell)}_{j,\phi} \mapsto  -\alpha^{(\ell)}_{j,\phi}, \\[6pt]
\delta x \mapsto \delta x\,, && \delta y \mapsto -\delta y\,. 
\end{array}
$$

Each of the subspaces $\left<\alpha^{(\ell)}_{j,\theta} ,\,\beta^{(\ell)}_{j,\theta} \right>$ and similarly with $\phi$ are $\DD_n$-invariant isotropic subspaces. 

\medskip

In order to compute a specific basis for the symplectic slice it is necessary to have expressions for the derivative of the momentum map and the tangent space to the group orbit at a $\CC_{nv}(k_1R,k_2R',k_pp)$ configuration.  These expressions are given in the next two propositions.  Since $\CC_{nv}$ refers to fixed vertical reflexion planes, the values of $\phi_j^0$ in (\ref{eqn:alphas,betas}) above can be taken to be:
$$
 \phi_j^0 = \begin{cases}0 & \text{if $j=1\dots k_1$}\cr \pi/n & \text{if
 $j=(k_1+1)\dots k_r$}\end{cases}
$$
where $k_r=k_1+k_2$ is the total number of rings.

\begin{proposition} \label{prop:kernel} 
At a $\CC_{nv}(k_1R,k_2R',k_pp)$ configuration, the
differential of the momentum map is as follows: let %
$$\uu=\sum_{j,\ell}\left(a_{j,\theta}^{(\ell)}\alpha_{j,\theta}^{(\ell)} +   a_{j,\phi}^{(\ell)}\alpha_{j,\phi}^{(\ell)} +
b_{j,\theta}^{(\ell)}\beta_{j,\theta}^{(\ell)} +
b_{j,\phi}^{(\ell)}\beta_{j,\phi}^{(\ell)}\right) +
\sum_{j\;\textrm{polar}}(c_j\delta x_j + d_j\delta y_j),$$%
be a tangent vector (the $a_{j,\theta}^{(\ell)}$ etc.\ are its
components with respect to the Fourier basis, where $j$ numbers the rings and $\ell$ the modes in each ring), then
\begin{eqnarray*}
\dd\Phi(\uu) &=& \sum_{j=1}^{k_r}
\kappa_j\cos\theta_j(a_{j,\theta}^{(1)} + \ii b_{j,\theta}^{(1)}) + \ii\sum_{j=1}^{k_r}
\kappa_j\sin\theta_j(a_{j,\phi}^{(1)} + \ii b_{j,\phi}^{(1)}) + 
\sum_{j\;\textrm{polar}}\kappa_j(c_j + \ii d_j)
\\
&&\quad \bigoplus -\sum_{j=1}^{k_r} \kappa_j\sin\theta_j a^{(0)}_{j,\theta}
\end{eqnarray*}
where the direct sum corresponds to the $\CC_{nv}$-invariant decomposition of $\so^*$ as a direct sum of a plane (the `$x$-$y$-plane') and the line $\Fix(\CC_{nv},\,\so^*)$ (the `$z$-axis').
\end{proposition}

This is proved from the definition of the momentum map (\ref{eq:momentum map}), and some Fourier type calculations. The kernel of the momentum map is readily read off from this. The following proposition results from similar calculations involving the infinitesimal rotation matrices in $\R^3$:

\begin{proposition} \label{prop:tangentspace}
Let $x_e$ be a $\CC_{nv}(k_1R,k_2R',k_pp)$ configuration, and $\mu=\Phi(x_e)$.
If $\mu\neq 0$, then the tangent space to the orbit $\so_{\mu}\cdot
x_e$ is generated by the vector
$$\sum_{j=1}^{k_r} \alpha^{(0)}_{j,\phi}.$$ %
If $\mu = 0$, then $\so_{\mu}\cdot x_e = \so\cdot x_e$ and is generated by the vector above, and the following two other vectors,
\begin{eqnarray*}
\sum_{j=1}^{k_r}\left(\beta^{(1)}_{j,\theta}  + \cot\theta_j\,\alpha^{(1)}_{j,\phi}\right)\; + 
  \sum_{j\;\textrm{polar}} {\rm sign}(z_j)\,\delta y_j,\\[6pt]
\sum_{j=1}^{k_r} \left(\alpha^{(1)}_{j,\theta} - \cot\theta_j\,\beta^{(1)}_{j,\phi}\right)\;+ 
  \sum_{j\;\textrm{polar}} {\rm sign}(z_j)\,\delta x_j.
\end{eqnarray*}
\end{proposition}

Now, $\so_\mu\cdot x_e\subset \ker\dd\Phi_{x_e}$ and the symplectic slice $\NN$ is the orthogonal complement to the former inside $\ker\dd\Phi$ (or indeed any complementary subspace, not necessarily orthogonal, though to take advantage of the symmetry group the subspace should chosen to be invariant under $G_{x_e}$). 

%%%%%%%%%%%%%%%%%%%%%%%%%%%%%%%%%%%%%%%%%%%%%%%%%

\subsection{Single ring}
\label{singlering}

Let $x_e$ be a $\CC_{nv}(R,k_pp)$ configuration, a single ring of $n$ vortices together with $k_p$ polar vortices where, of course, $k_p=0,1$ or $2$. Since there is only one ring, we write
$\alpha^{(\ell)}_{\theta},\beta^{(\ell)}_{\theta}$ etc.\ instead of $\alpha^{(\ell)}_{1,\theta},\beta^{(\ell)}_{1,\theta}$ etc.

\begin{proposition} \label{prop:isotypic}
For these single ring configurations, the symplectic slice can be
decomposed by Fourier modes (or isotypic representations) as
$$\NN =  \bigoplus_{\ell=1}^{[n/2]}V_\ell,$$
where  for $\ell\geq 2$,
  $$V_\ell = \<\alpha^{(\ell)}_{\theta},\alpha^{(\ell)}_{\phi},\beta^{(\ell)}_{\theta},\beta^{(\ell)}_{\phi}\right>$$
(with the understanding that $\beta_\theta^{(n/2)} = \beta_\phi^{(n/2)}=0$
when $n$ is even). So $\dim V_\ell=4$ for $1<\ell<n/2$, and if $n$
is even $\dim V_{n/2}=2$. Furthermore,
$$\dim V_1 = \begin{cases}
 2k_p+2 & \text{if $\mu\neq0$ and $n>2$}\cr
 2k_p & \text{if $\mu=0$ and $n>2$}\cr
 2k_p   & \text{if $\mu\neq0$ and $n=2$}\cr
 \max\{0,2k_p-2\} & \text{if $\mu=0$ and $n=2$}.
    \end{cases}$$
A basis for $V_1$ follows from Propositions \ref{prop:kernel} and \ref{prop:tangentspace}, and is given in the relevant section.
\end{proposition}

\paragraph{Representations}
As representations of $G_{x_e}\simeq\CC_{nv}\simeq\DD_n$ and
$G_{x_e}^0\simeq\CC_n$, the $V_\ell$ are pairwise distinct, so by
Schur's Lemma both the Hessian $\dd^2H_{\xi}|_{\NN}$ and the
linearized vector field $L_\NN$ block diagonalize with respect to
this decomposition. We now give a brief description of these
representations. These statements follow directly from the
definition of the Fourier modes in (\ref{eqn:alphas,betas}) above.

For $2<\ell<n/2$, the Fourier subspace $V_\ell$ is 4-dimensional,
and consists of a direct sum of two  2-dimensional
symplectic irreducibles, one isomorphic to the usual representation
of $\CC_{n}$ acting on the plane with `speed' $\ell$ (that is, a
rotation by $\rho$ in $\CC_{n}$ acts on the plane by a rotation
through $\ell\rho$), and the other the dual representation. Of course, the reflexions in $\CC_{nv}$ act antisymplectically.

For $\ell=n/2$ (when $n$ is even), $V_{n/2}$ is 2-dimensional, with
$G_{x_e}^0$ acting by $\{\pm  I\}$, but $G_{x_e}$ acting by 
reflexion in a line. More specifically,
$V_{n/2}=\left<\alpha_\theta^{(n/2)},\,\alpha_\phi^{(n/2)}\right>$
and one generator of $G_{x_e}$ acts trivially on
$\alpha^{(n/2)}_{\theta}$ and by $-1$ on $\alpha^{(n/2)}_{\phi}$, while another acts the other way around.

The subspace $V_1$ is more involved. The representations 
$\<\alpha^{(1)}_{\theta},\beta^{(1)}_{\theta}\m,\; \<
  \alpha^{(1)}_{\phi},\beta^{(1)}_{\phi} \m$ and $\< \delta x_j , \delta
  y_j \m$ ($j=1,\dots, k_p$) 
do not lie in the symplectic slice, and it is necessary to use
Propositions \ref{prop:kernel} and \ref{prop:tangentspace} to write down the
bases of $V_1$ which will be central in later calculations. However,
the underlying representations are all isomorphic (i.e., ignoring the symplectic structure); they are copies of the basic 2-dimensional representations of the dihedral group $\CC_{nv}$, and it follows that $V_1$ is a direct sum of such representations. But as \emph{symplectic} representations this is more subtle: the subspaces $\left<\delta x_j,\delta y_j\right>$ are symplectic, while the others mentioned are isotropic.

%%%%%%%%%%%%%%
\subsection{General case}
\label{tworings}

In the general case where $x_e$ is a $\CC_{nv}(k_1R,k_2R',k_pp)$
configuration ($k_p=0,1$ or 2), we have the following decomposition.

\begin{proposition} For these multiring configurations, the
symplectic slice decomposes as a direct sum of Fourier modes
$$\NN=\bigoplus_{\ell=0}^{[n/2]} V_\ell$$
where  for $2\leq \ell \leq [n/2]$,
$$V_\ell = \<\alpha^{(\ell)}_{j,\theta},\alpha^{(\ell)}_{j,\phi},\beta^{(\ell)}_{j,\theta},\beta^{(\ell)}_{j,\phi}\mid j=1\dots k_r\right>,
$$
which is of dimension $4k_r$ if $\ell<n/2$ and $2k_r$ if $\ell=n/2$, where $k_r$ is the total number of rings. The subspaces $V_0$ and $V_1$ are subspaces of the corresponding
spaces determined by Propositions \ref{prop:kernel} and
\ref{prop:tangentspace} and are given in the relevant sections.
\end{proposition}

The spaces $V_\ell$ as representations of $G_{x_e}$ and $G_{x_e}^0$
and $2\leq\ell\leq n/2$ will just be the sum of several copies of
the corresponding $V_\ell$ for a single ring.

%%%%%%%%%%%%%%%%%%%%%%%%%%%%%%%%%%%%%%%%%%%%%%%%%%%%%%%%%%%%%%%%%%%%%%%%%
\section{Bifurcations and symmetry}
\label{sec:symmetry+bifurcations}

The configurations we consider in this paper depend on a number of parameters, both external and internal (the internal parameter is the conserved momentum, the external ones the vorticity). As a result, one finds many bifurcations and since  the system has symmetry, the bifurcations will reflect this symmetry, and in particular will depend on how the symmetry group of a given relative equilibrium (\re) acts on the eigenspace in which the bifurcation occurs. Of course, the transitions involving a change in stability are such bifurcations. Here we present a brief zoology of the bifurcations that are encountered.

The simplest steady-state bifurcation in a Hamiltonian system is the saddle-node, where two equilibria come together and vanish.  However, because we are considering symmetric relative equilibria (\re), there is always a `central' symmetric \re\ which persists throughout the bifurcation and so the saddle-node will not arise.  Consequently, the simplest bifurcation we encounter is the pitchfork. 

Some details about how the action of the symmetry group influence which steady-state bifurcations to expect generically were elaborated in \cite{GS87}, see also \cite{BLM05} for a summary.  However, only symplectic actions were considered there (without time-reversing symmetries), and the extra time-reversing symmetries restrict what happens generically so their results are not directly applicable.

\begin{remark}\label{rmk:rpos}
When the relative equilibrium (or reduced equilibrium) has a mode with imaginary eigenvalues, one can apply the symmetry methods of \cite{MRS88,MRS90} to find families of relative periodic orbits in a neighbourhood of the relative equilibrium with prescribed spatio-temporal symmetries.  We do not pursue this here.
\end{remark}

%%%%%%%%%%%
\subsection{Symmetric pitchfork bifurcations}
\label{sec:pitchfork}

%%%%%%%%%%%%%%%%%%%%%%%%%%%%%%%%%%%%%%
\begin{figure}[t]
\begin{center}
\subfigure[Transcritical $\DD_4$-pitchfork]%
{\includegraphics[scale=0.15]{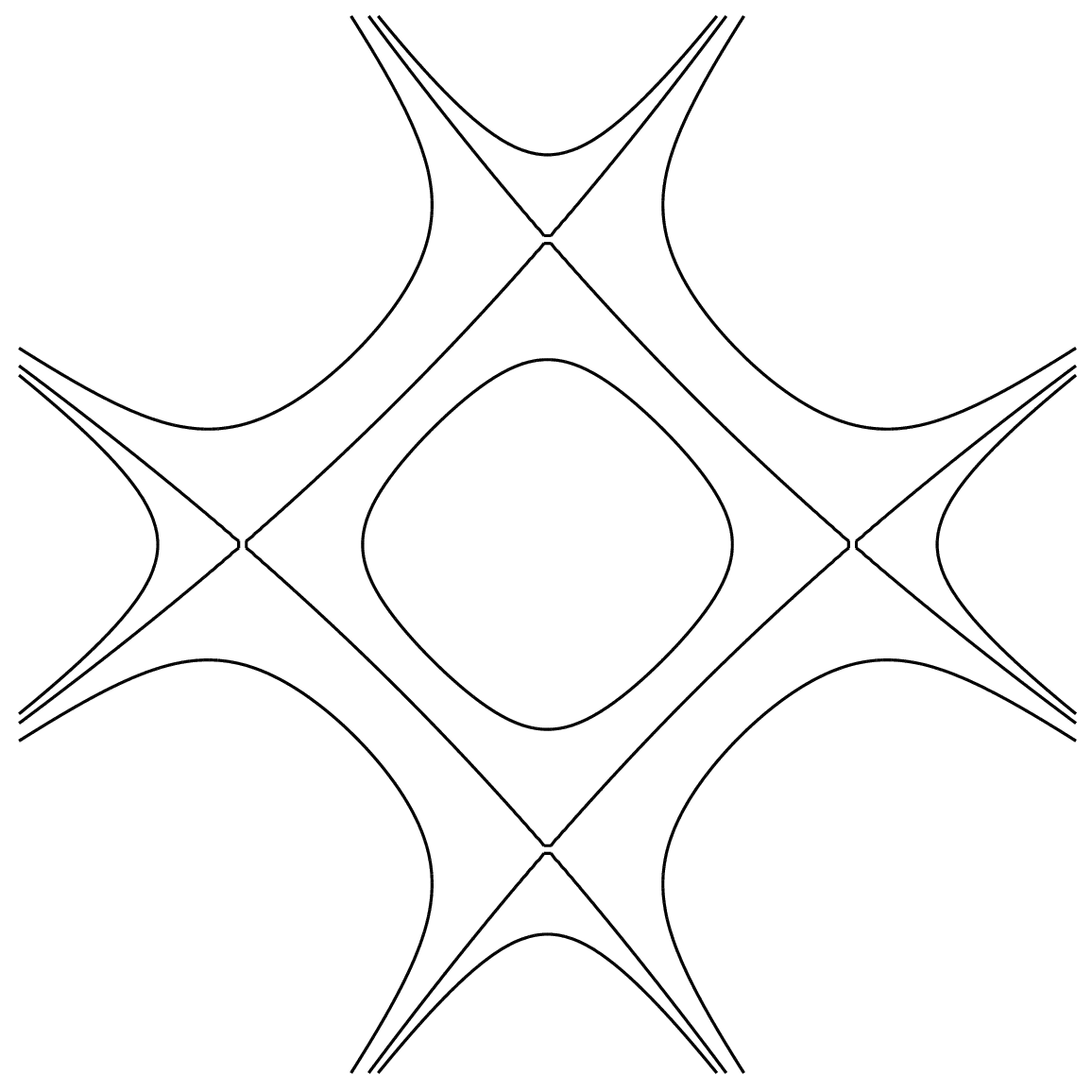}
 \qquad\qquad
\includegraphics[scale=0.15]{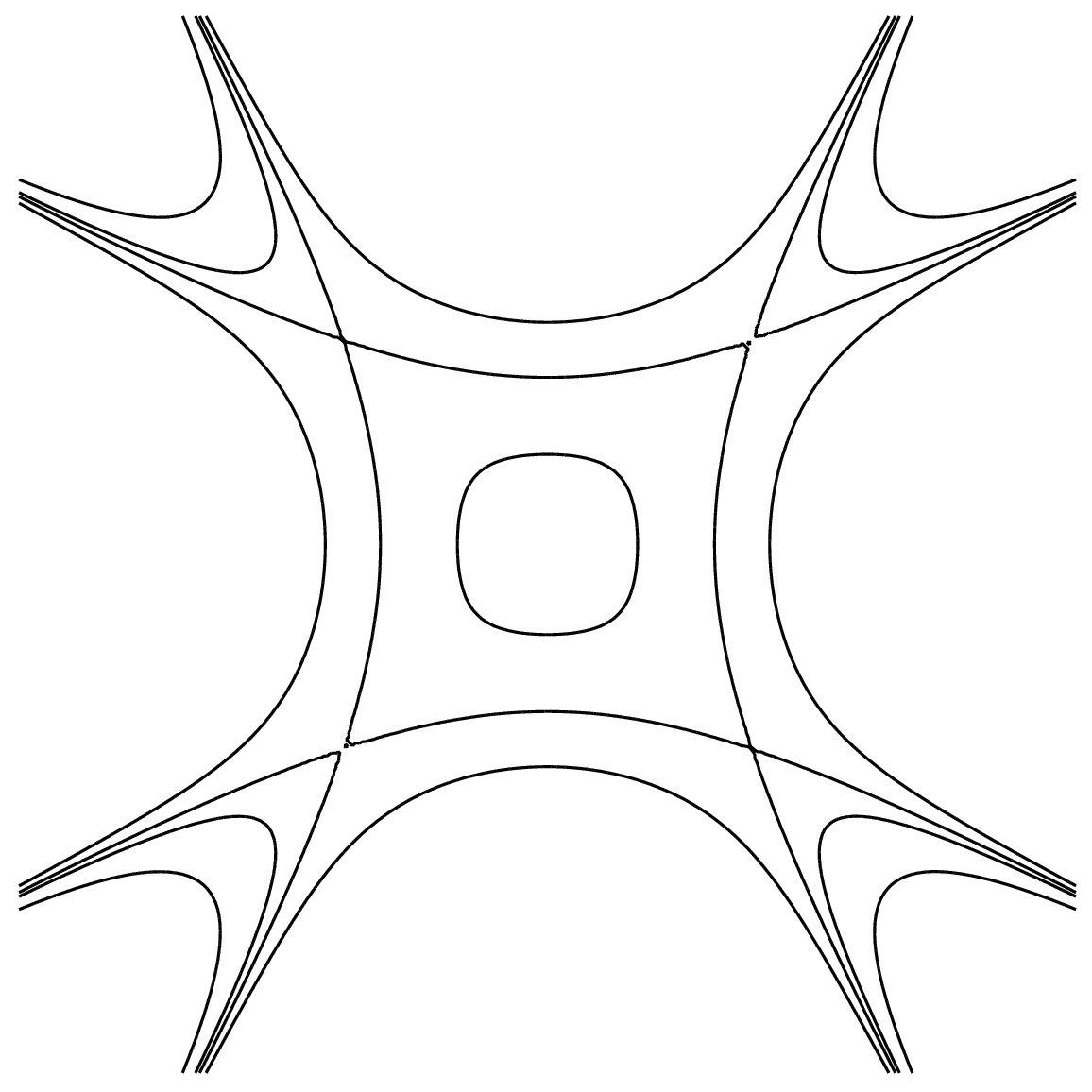}
}

\bigskip

\subfigure[`Standard' $\DD_4$-pitchfork]%
{\includegraphics[scale=0.15]{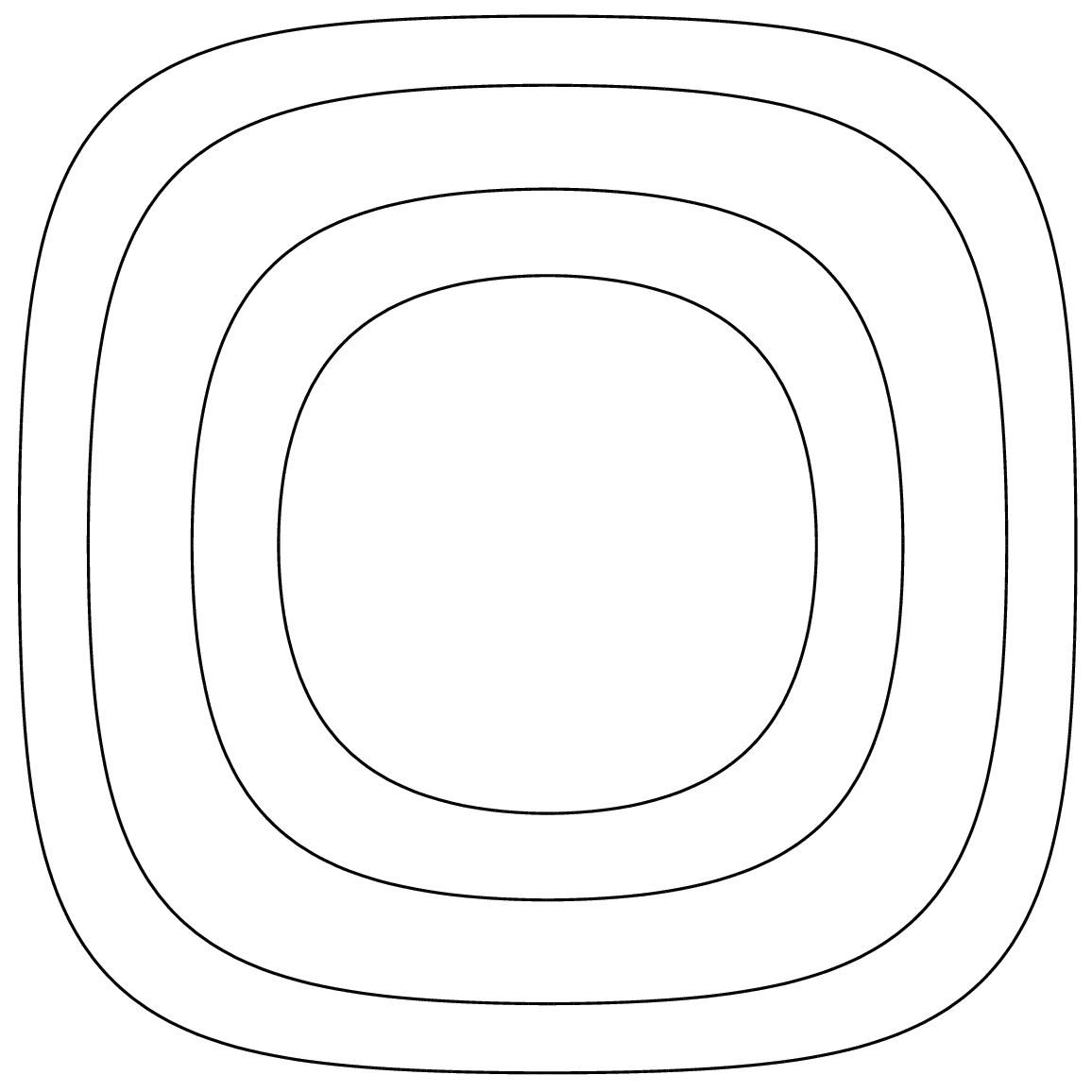}
 \qquad\qquad
\includegraphics[scale=0.15]{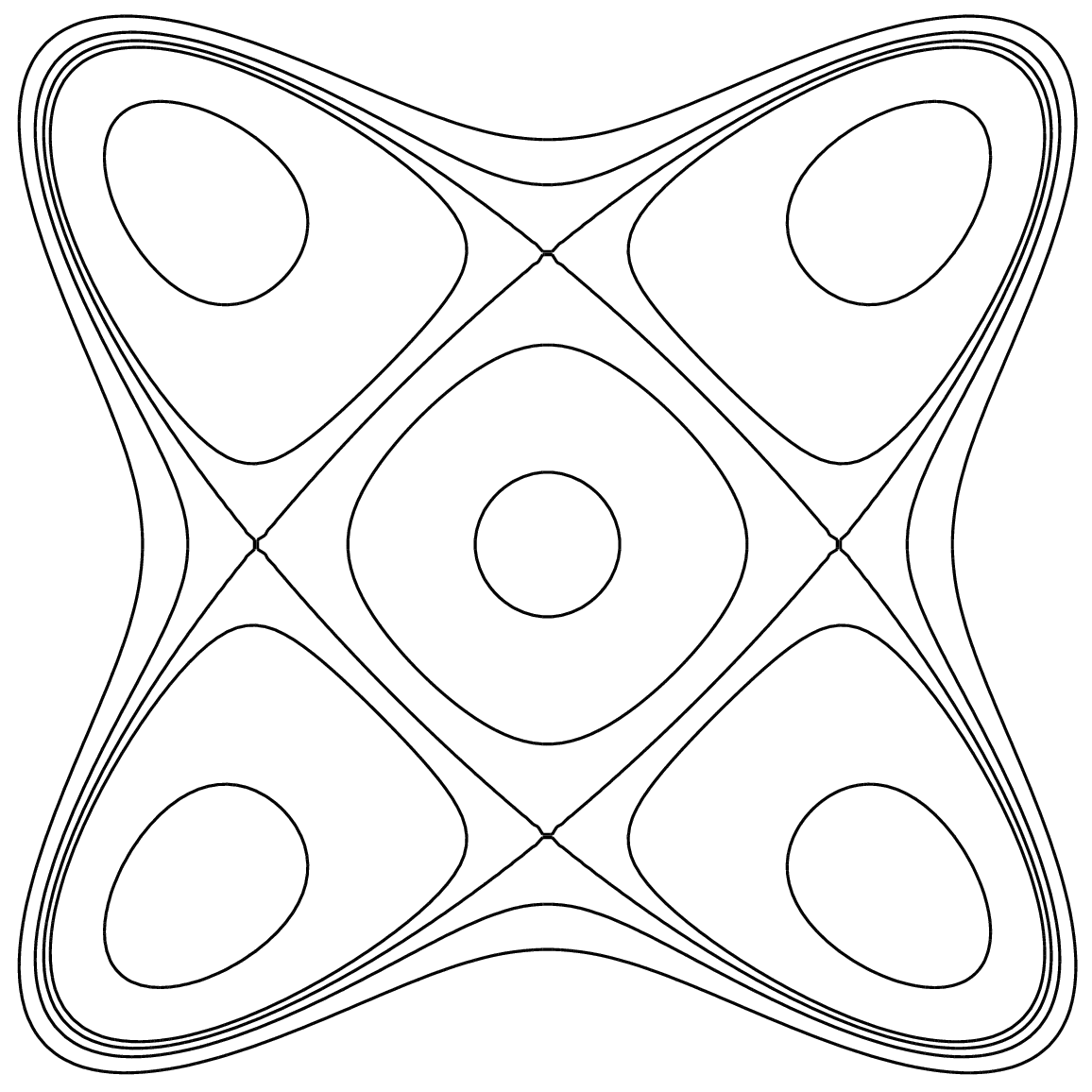}
}

\caption{Contours of the generic 1-parameter family of $\DD_{4}$-invariant functions in the plane 
$H_\lambda(z) = \lambda|z|^2 + a|z|^4+b\Re(z^4) + O(z^5)$. 
The top pair (a) represent a transcritical bifurcation and occur if $|a|<|b|$, while the bottom pair (c) and (d) represent a standard pitchfork bifurcation, occurring if $|a|>|b|$. In both cases the left-hand figure represents $\lambda<0$ and the right-hand one $\lambda>0$.
(Distinguishing between sub- and super-critical cases depends on stability information in other `modes'.) See text for how these relate to other dihedral groups.
}
\label{fig:dihedral bifurcations}
\end{center}
\end{figure}
%%%%%%%%%%%%%%%%%%%%%%%%%%%%%%%%%%%%%%

In this bifurcation, there is a central `core' \re\ which persists throughout the bifurcation, and whose eigenvalues (of the linearization) pass through zero as parameters are varied.  As the transition occurs,  other \re\ branch off from the central one. In our setting this bifurcation occurs with the symmetry of a dihedral group $\DD_k$ (where $k$ depends on both $\ell$ and $n$).  

One can distinguish between transcritical pitchfork bifurcations on the one hand and standard pitchforks on the other; in the former there are bifurcating solutions on both sides of the bifurcation point, while in the latter the bifurcating solutions all coexist for the same parameter values: see Figure~\ref{fig:dihedral bifurcations}.  Which occurs depends principally on the value of $k$.  A typical pitchfork bifurcation with $\DD_3$ symmetry is transcritical, analogous to (a) in the figure but with 3-fold symmetry, while a typical pitchfork with $\DD_k$-symmetry for $k>4$ is a standard pitchfork analogous to (b) but again with $k$-fold symmetry. If $k=4$ then either case can occur, as described in the figure.

Note that if a standard pitchfork bifurcation involves a change in stability then it is either sub- or super-critical, depending on whether the bifurcating solutions are stable or unstable respectively. 

In the bifurcations occurring in this paper, if the symmetry is $\CC_{nv}\simeq\DD_n$, and the bifurcation is in mode $\ell$, then the effective symmetry for the bifurcation is $\DD_k$ where $k=n/(n,\ell)$, where $(n,\ell)$ is the highest common factor of $n$ and $\ell$.  This is because all points in this mode have $\Z_{(n,\ell)}$-symmetry, see (~\ref{eq:action on alphas,betas}), leaving an effective action of $\DD_n/\Z_{(n.\ell)}\simeq \DD_k$. 

For example, if $n$ is even and $\ell=n/2$ then $k=2$ and there is an effective action of $\DD_2\simeq\Z_2\times\Z_2$, and all points of this mode have $\Z_{n/2}$-symmetry. It follows that all if the bifurcation involves this mode, then all bifurcating solutions have at least this symmetry. In fact an analysis of the generic bifurcation shows that they have a further $\Z_2$ giving $\DD_{n/2}$ symmetry.  The generic pitchfork bifurcation in this mode is the well-known one where there are two branches of bifurcating solutions and a $\Z_2$-symmetry is lost; since the central equilibrium has $\DD_n$ symmetry the bifurcating solutions will still have a $\DD_{n/2}$-symmetry.  

At the other extreme, if $n$ and $\ell$ are coprime, then $k=n$ and the bifurcating solutions have only $\Z_2$-symmetry, generated by one of the reflexions in $\DD_n$.  

Figures \ref{fig:bifurcating re for ell = n/2} and \ref{fig:ell=2 bifurcations} show symmetries of bifurcating solutions.

\medskip

\begin{figure}[t]
\begin{center}
 \unitlength=0.4mm \thinlines
\begin{picture}(40,80)(-20,-60)
 \put(-20,0){\line(1,0){40}}
 \put(0,-20){\line(0,1){40}}
 \put(0,-10){\circle*{3}}
 \put(0,10){\circle*{3}}
 \put(5,12){\vector(0,-1){8}}
 \put(5,-12){\vector(0,1){8}}
 \put(-8,-30){$a<0$}
\end{picture}\quad
\begin{picture}(40,80)(-20,-60)
 \put(-20,0){\line(1,0){40}}
 \put(0,-20){\line(0,1){40}}
 \put(0,0){\circle*{5}}
 \put(-8,-30){$a=0$}
 \put(-15,-48){\textsc{splitting}}
\end{picture}\quad
\begin{picture}(40,80)(-20,-60)
 \put(-20,0){\line(1,0){40}}
 \put(0,-20){\line(0,1){40}}
 \put(-10,0){\circle*{3}}
 \put(10,0){\circle*{3}}
 \put(8,-3){\vector(1,0){8}}
 \put(-8,-3){\vector(-1,0){8}}
 \put(-8,-30){$a>0$}
\end{picture}
\qquad\qquad
\begin{picture}(40,80)(-20,-60)
 \put(-20,0){\line(1,0){40}}
 \put(0,-20){\line(0,1){40}}
 \put(0,-10){\circle*{3}}
 \put(0,10){\circle*{3}}
 \put(5,12){\vector(0,-1){8}}
 \put(5,-12){\vector(0,1){8}}
 \put(-8,-30){$a<0$}
\end{picture}\quad
\begin{picture}(40,80)(-20,-60)
 \put(-20,0){\line(1,0){40}}
 \put(0,-20){\line(0,1){40}}
 \put(0,0){\circle*{5}}
 \put(-8,-30){$a=0$}
 \put(-11,-48){\textsc{passing}}
\end{picture}\quad
\begin{picture}(40,80)(-20,-60)
 \put(-20,0){\line(1,0){40}}
 \put(0,-20){\line(0,1){40}}
 \put(0,-10){\circle*{3}}
 \put(0,10){\circle*{3}}
 \put(-8,-30){$a>0$}
 \put(5,8){\vector(0,1){8}}
 \put(5,-8){\vector(0,-1){8}}
\end{picture}
\caption{The two scenarios for the generic movement of eigenvalues in symmetric steady-state pitchfork bifurcation.}
\label{fig:passing/splitting}
\end{center}
\end{figure}
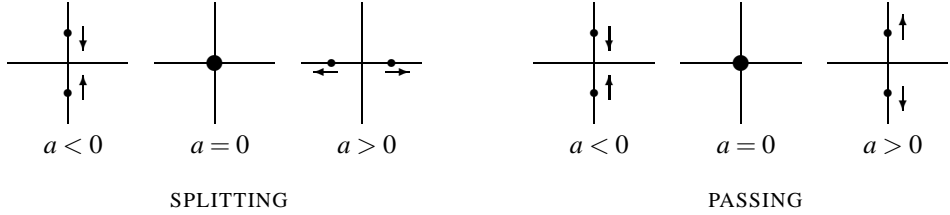

As the eigenvalues pass through zero in either type of pitchfork bifurcation there are two possibilities, the so-called \emph{splitting} and \emph{passing} behaviours (see Figure~\ref{fig:passing/splitting}):

\begin{description}
\item[Pitchfork of splitting type] This is the familiar scenario where a complex conjugate pair of imaginary eigenvalues of the linear system collide at 0 and then become real and opposite. The central (relative) equilibrium is then possibly elliptic or Lyapounov stable on one side of the bifurcation, and linearly unstable on the other.  
\item[Pitchfork of passing type] This is similar, except that the eigenvalues pass through each other and remain on the imaginary axis. See Figure~\ref{fig:passing/splitting}. If the central \re\ is say Lyapounov stable (or elliptic) before the bifurcation, then in a pitchfork of passing type the central \re\ becomes (remains) elliptic, and in the standard pitchfork bifurcation $k$ of the bifurcating solutions are Lyapounov stable (elliptic) while the other $k$ are linearly unstable. In the transcritical case, all the bifurcating \re\ are linearly unstable.
\end{description}

%%%%%%%%%%%%%%%%%%

For example, if a bifurcation takes place due to a degeneracy in the mode space $V_{n/2}$, then generically it will be a standard $\Z_2$-pitchfork bifurcation, and the central relative equilibrium will pass from being (Lyapounov) stable to being linearly unstable. A study of higher order terms will distinguish between standard sub- and super-critical pitchfork bifurcations (something we don't do in this paper). 

If on the other hand a bifurcation occurs in the mode $V_\ell$ with $0<\ell<n/2$ then generically one may see pitchforks of either  passing or splitting type.

%%%%%%%%%%%
\subsection{Hamiltonian Hopf bifurcation}
\label{sec:HamHopf}

This occurs when two imaginary eigenvalues of a (relative) equilibrium collide and move off the imaginary axis.  This can only happen if the (reduced) Hessian matrix at the (relative) equilibrium is not definite.

There are two possible scenarios depending on the higher order terms:  \\
(a) the two Lyapounov modes associated with the two imaginary eigenvalues are globally connected, and the set of all of these collapse to a point as the eigenvalues collide, and\\
(b) the two Lyapounov modes are not part of the same family, and as the eigenvalues move off the imaginary axis, so the two families connect and move away from the (relative) equilibrium. 

For a detailed description of the bifurcations of periodic and relative periodic orbits that can be expected in this case see \cite{MS71, vM85, COR02, BLM05}.

%%%%%%%%%%%%%
\subsection{Bifurcations from zero momentum}
\label{sec:zero momentum}

Consider a relative equilibrium $x_0$ in a system with $\SO(3)$-symmetry with zero momentum, whose reduced Hessian is non-degenerate. 

If the angular velocity of the \re\ is non-zero, then through that \re\ there is  a curve of \re\ parametrized by $\mu$. If the reduced Hessian of $x_0$ is positive definite (so $x_e$ is Lyapounov stable) then generically on one side of the curve the \re\ is also Lyapounov stable, while on the other it is elliptic, \cite{Mo11}. 

The genericity assumption is that the frequency of the rotation of the \re\ (which is $|\xi|$) is distinct from the frequencies of the normal modes 
of the reduced system. This \emph{ro-vibrational resonance} is discussed briefly at the end of \cite{Mo97}, where it is made clear how it is related to the  Hamiltonian Hopf bifurcation; the associated dynamics is the subject of a paper by Patrick \cite{Pa99}, where he calls it group-reduced resonance. This resonance phenomenon can be seen for a ring and 2 polar vortices in Section~\ref{sec:R2p}. The bifurcations arising as parameters are varied have not been studied. 

If on the other hand the angular velocity of $x_0$ is zero then for each non-zero $\mu$ there are at least 6 \re, and if $x_0$ is definite as above, then at least 2 of these are each of Lyapounov stable, elliptic and unstable. Moreover, if the equilibrium in question has symmetry, then the symmetry can be used to give precise estimates for the numbers of \re\ in a neighbourhood. Details can be found in \cite{MR99} and a study of deformations from zero velocity to non-zero velocity can be found in \cite{Mo11}.

%%%%%%%%%%%%%%%%%%%%%%%%%%%%%%%%%%%%%%%%%%%%%%%%%%%%%%%%%%%%%%%%%%%%%%%%%
\section{A ring of identical vortices: $\CC_{nv}(R)$}
\label{sec:Cnv(R)}

The linear stability of $\CC_{nv}(R)$ relative equilibria was
determined by Polvani and Dritschel in \cite{PD93}. A few years ago, Boatto
and Cabral \cite{BC03} studied their Lyapounov stability and found
that the two types of stability coincide: whenever the relative
equilibrium fails to be Lyapounov stable the linearization of $X_H$
has real eigenvalues.  In this section, we give another proof using
the geometric method of this paper. The calculations are also used
in later sections.

\subsection{Symplectic slice}

For $n$ vortices of unit vorticity the Hamiltonian is
\begin{equation} \label{eq:ring Hamiltonian}
H(\theta_j,\phi_j) \ =\ -\sum_{j<k}
\ln\left(1-\sin\theta_j\sin\theta_k\cos(\phi_j-\phi_k) -
  \cos\theta_j\cos\theta_k\right)
\end{equation}
and the augmented Hamiltonian is $H_\xi=H-\xi\sum_j\cos\theta_j$.

Let $x_e$ be a $\CC_{nv}(R)$ relative equilibrium and $\theta_0$ the
co-latitude of the ring.  The angular velocity of $x_e$ is
$$\xi \ =\  \frac{(n-1)\cos\theta_0}{\sin^2\theta_0},$$
since $H_{\xi}$
has a critical point there and $\frac{\partial
  H_\xi}{\partial\theta_j}(x_e) =
\frac{(n-1)\cos\theta_0-\xi\sin^2\theta_0}{\sin\theta_0}$. The
momentum of the relative equilibrium is just $\mu=n\cos\theta_0$.
\bigskip

The second derivatives of $H$ at the relative equilibrium are:
$$\left.\begin{array}{rclcrcl} \frac{\partial^2H}{\partial\theta_j^2}
    &=& -\frac{(n-1)(n-5)}{6\sin^2\theta_0} &\qquad&
    \frac{\partial^2H}{\partial\theta_j\partial\theta_k} &=&
    \frac{1}{2\sin^2\theta_0\sin^2(\pi(j-k)/n)} \\[6pt]
    \frac{\partial^2H}{\partial\theta_j\partial\phi_j} &=& 0 &&
    \frac{\partial^2H}{\partial\theta_j\partial\phi_k} &=& 0 \\[6pt]
    \frac{\partial^2H}{\partial\phi_j^2} &=&
    \frac16(n^2-1) &&
    \frac{\partial^2H}{\partial\phi_j\partial\phi_k} &=&
    -\frac{1}{2\sin^2(\pi(j-k)/n)}.
\end{array}\right.
$$
We note that $\sum_{r=1}^{n-1}1/\sin^2(\pi r/n) =
\frac{1}{3}(n^2-1)$ \cite{H75}. 

\bigskip

\paragraph{Notation}
In order to harmonize the statements of the results between $n$ even
and $n$ odd, we introduce the following notation: let
$\eta^{(\ell)}_1,\eta^{(\ell)}_2,\eta^{(\ell)}_3,\eta^{(\ell)}_4$ be
objects defined for all $2\leq \ell\leq \left[\frac{n-1}{2}\right]$,
where $[m]$ is the integer part of $m\in\N$, and only
$\eta^{(\ell)}_1,\eta^{(\ell)}_2$ for $\ell=n/2$ when $n$ is even.
Then define
$$\left\lbrace
  \eta^{(\ell)}_1,\eta^{(\ell)}_2,\eta^{(\ell)}_3,\eta^{(\ell)}_4\mid
  2\leq\ell\leq [n/2]\right\rbrace^*
$$
to be
\begin{equation}\label{eq:set*}
 \left\{ \begin{array}{lll}
\mbox{for even $n$} &:&
\left\lbrace
  \eta^{(\ell)}_1,\eta^{(\ell)}_2,\eta^{(\ell)}_3,\eta^{(\ell)}_4\mid
  2\leq\ell\leq \frac{n}{2}-1\right\rbrace\cup\left\lbrace
  \eta^{(n/2)}_1,\eta^{(n/2)}_2\right\rbrace   \\
\mbox{for odd $n$} &:& \lbrace
\eta^{(\ell)}_1,\eta^{(\ell)}_2,\eta^{(\ell)}_3,\eta^{(\ell)}_4\mid
2\leq\ell\leq [n/2]\rbrace.
  \end{array}\right.
\end{equation}

In Proposition~\ref{prop:isotypic} a basis is given for each of the
$V_\ell$ except for $\ell=1$.  The following proposition does the
same for $V_1$.

\begin{proposition} \label{diagonering}
For $n=2$ the symplectic slice is 0. For $n\geq 3$, the symplectic
slice decomposes as $\NN = \bigoplus_{\ell=1}^{[n/2]} V_\ell,$ where
for $\mu\neq0$
 \[\dim V_\ell =
\begin{cases}2 & \text{if $\ell=1$ or $n/2$ }\cr
       4 & \text{otherwise}.\end{cases}
\]
If\/ $\mu=0$ then $V_1=0$, while for $\mu\neq0$, $V_1 =
\left<e_1,e_2\right>$, where
$$
\begin{array}{lll}
e_1&=&\sin\theta_0\;\alpha^{(1)}_{\theta}+\cos\theta_0\;\beta^{(1)}_{\phi}\\
e_2&=&\sin\theta_0\;\beta^{(1)}_{\theta}-\cos\theta_0\;\alpha^{(1)}_{\phi}.
\end{array}
$$
With respect to the resulting basis for the symplectic slice, the
Hessian $\dd^2H_\xi|_\NN(x_e)$ is diagonal, and $L_\NN$ block
diagonalizes in $2\times 2$ blocks.
\end{proposition}

\begin{proof}
  It is straightforward to check that the vectors above do form
  a basis for $V_1$ at $x_e$ thanks to Propositions
  \ref{prop:kernel} and \ref{prop:tangentspace}.

  The Hessian $\dd^2H_\xi|_\NN(x_e)$ and the linearization $L_\NN$ are
  both $G_{x_e}^0$-invariant.  Assume $n$ odd.  It follows from
  Section \ref{singlering} and Schur's Lemma (see the introduction) that
  $\dd^2H_\xi|_\NN(x_e)$ and $L_\NN$ both block diagonalize into $4\times 4$
  blocks and one $2\times 2$ block corresponding to the subspaces
  $$V_\ell=\<\alpha^{(\ell)}_{\theta},\beta^{(\ell)}_{\theta},
  \alpha^{(\ell)}_{\phi},\beta^{(\ell)}_{\phi}\m$$ 
   and  $\<e_1,e_2 \m$,
  respectively. See the proof of Theorem~4.5 of \cite{LP02} for a
  detailed proof of a similar assertion.

  Now fix $\ell$ and denote by $s$ an anti-symplectic (time-reversing)
  element of ${G_{x_e}}$.  For example $s$ could be the reflexion
  $y\mapsto -y$ together with an order two permutation of $S_n$.
  The restriction of $H_\xi$ to $V_\ell$ is $\Z_2[s]$-invariant. Moreover
  $\<\alpha^{(\ell)}_{\theta},\beta^{(\ell)}_{\phi}\m$ and
  $\<\beta^{(\ell)}_{\theta},\alpha^{(\ell)}_{\phi}\m$ are
  non-isomorphic irreducible representation of $\Z_2[s]$ on $V_\ell$.
  Hence $\dd^2H_\xi|_\NN(x_e)$ block diagonalizes into $2\times 2$ blocks
  which correspond to subspaces
  $\<\alpha^{(\ell)}_{\theta},\beta^{(\ell)}_{\phi}\m$,
  $\<\beta^{(\ell)}_{\theta},\alpha^{(\ell)}_{\phi}\m$, and $\<
  e_1,e_2 \m$.  This result does not depend on the details of the
  Hamiltonian, only its symmetries.  However taking account of its particular
  form, one can improve the block diagonalization. Indeed one has
  $$\dd^2H_\xi(x_e)\cdot
  (\alpha^{(\ell)}_{\theta},\beta^{(\ell)}_{\phi})\ =\
  \dd^2H_\xi(x_e)\cdot
  (\beta^{(\ell)}_{\theta},\alpha^{(\ell)}_{\phi})\ =\ 0$$
  and
  $\dd^2H_\xi(x_e)\cdot(e_1,e_2)=0$ which gives the desired
  diagonalization of the Hessian.

  The particular form of the symplectic form also enables us to improve
  the diagonalization of $L_\NN$. Among the basis vectors of
  $V_\ell$, only
  $\omega(\alpha^{(\ell)}_{\theta},\alpha^{(\ell)}_{\phi})$ and
  $\omega(\beta^{(\ell)}_{\theta},\beta^{(\ell)}_{\phi})$ do not
  vanish, and so the restriction of $\omega$ to $V_\ell$ block
  diagonalizes into two $2\times 2$ blocks which correspond to the
  subspaces $\<\alpha^{(\ell)}_{\theta},\alpha^{(\ell)}_{\phi}\m$ and
  $\<\beta^{(\ell)}_{\theta},\beta^{(\ell)}_{\phi}\m$.  The block
  diagonalization of $L_\NN$ then follows from
  $L_\NN=\J_{\NN}\dd^2H_{\xi}|_{\NN}(x_e)$.

  The case $n$ even is very similar, except that there is an additional
  $2\times 2$ block in the $G_{x_e}^0$-isotypic decomposition, and
  leads to the same result.
\end{proof}

%%%%%%%%%%%%%%%%%%%%%%
\subsection{Stability}
The block diagonalization of $\dd^2H_\xi|_\NN(x_e)$ and $L_\NN$ enable
us to find formulae for their eigenvalues, and thus to conclude
criteria for both Lyapounov and linear stability.

\begin{theorem}\label{polvani}
  The stability of a ring of $n$ identical vortices depends on $n$ and
  the co-latitude $\theta_0$ as follows:

\begin{tabular}{c|l}
n & condition for stability\\
\hline
2, 3 & all $\theta_0$\\
4 & $\cos^2\theta_0 > 1/3$\\
5 & $\cos^2\theta_0 > 1/2$\\
6 & $\cos^2\theta_0 > 4/5$\\
$\geq$7 & always linearly unstable\\
\hline
\end{tabular}

\medskip

\noindent For $n=4,5,6$, the ring is linearly unstable if the inequality is reversed.

  % \begin{description}
  % \item[n = 2 and 3] is Lyapounov stable at all latitudes;
  % \item[n = 4] is Lyapounov stable if $\cos^2\theta_0 > 1/3$, and
  %   linearly unstable if the inequality is reversed;
  % \item[n = 5] is Lyapounov stable if $\cos^2\theta_0 > 1/2$, and
  %   linearly unstable if the inequality is reversed;
  % \item[n = 6] is Lyapounov stable if $\cos^2\theta_0 > 4/5$, and
  %   linearly unstable if the inequality is reversed;
  % \item[n$\geq$7] is always (linearly) unstable.
  % \end{description}
\end{theorem}

\begin{proof}
  Any arrangement of two vortices is a relative equilibrium
  \cite{KN98}. When perturbing such a relative equilibrium, we obtain
  a new relative equilibrium close to the first. Thus any relative equilibrium of
  two vortices is Lyapounov stable modulo $\SO(2)$ (and modulo $\SO(3)$ if
  $\mu=0$).

  Now assume $n\geq 3$.  We first study Lyapounov stability.
  Suppose further that $\mu\neq 0$ and so the ring is not
  equatorial.  A simple calculation shows that $\dd^2H_\xi(x_e)\cdot
  (\beta^{(\ell)}_{\theta},\beta^{(\ell)}_{\theta})=\dd^2H_\xi(x_e)\cdot
  (\alpha^{(\ell)}_{\theta},\alpha^{(\ell)}_{\theta})$ and
  $\dd^2H_\xi(x_e)\cdot
  (\beta^{(\ell)}_{\phi},\beta^{(\ell)}_{\phi})=\dd^2H_\xi(x_e)\cdot
  (\alpha^{(\ell)}_{\phi},\alpha^{(\ell)}_{\phi})$. Hence it follows
  from Proposition~\ref{diagonering} that
  $$\dd^2H_\xi|_\NN(x_e)\ =\ \diag\left(\lambda_1,\lambda_1,   \lbrace\lambda^{(\ell)}_\theta,\lambda^{(\ell)}_\phi,\lambda^{(\ell)}_\theta,\lambda^{(\ell)}_\phi
    \mid 2\leq \ell\leq [n/2]\rbrace^* \right)
  $$
  (recall notation from (\ref{eq:set*})) where
  $$
\begin{array}{lll}
\lambda_1&=&\sin^2\theta_0\lambda^{(1)}_\theta
+\cos^2\theta_0\lambda^{(1)}_\phi,\\[6pt]
\lambda^{(\ell)}_\theta&=&\dd^2H_\xi(x_e)\cdot
(\alpha^{(\ell)}_{\theta},\alpha^{(\ell)}_{\theta}),\\[6pt]
\lambda^{(\ell)}_\phi&=&\dd^2H_\xi(x_e)\cdot (\alpha^{(\ell)}_{\phi},\alpha^{(\ell)}_{\phi}).
\end{array}
$$
Thanks to the following formula \cite[p.271]{H75}
$$
\sum_{j=1}^{n-1} \frac{\cos(2\pi\ell j/n)}{\sin^2(\pi
  j/n)}\ =\ \frac{1}{3}(n^2-1)-2\ell(n-\ell),
$$
we find after some computations that
$\lambda^{(\ell)}_\phi=n\ell(n-\ell)/2$ and
\begin{equation}\label{eq:theta-evs for ring}
\lambda^{(\ell)}_\theta\ =\ \frac{n}{2\sin^2\theta_0}\left[
  -(\ell-1)(n-\ell-1)+(n-1)\cos^2\theta_0 \right].
\end{equation}
The eigenvalues $\lambda^{(\ell)}_\phi$ are all positive and
$\lambda_1=n(n-1)\cos^2\theta_0>0$. The relative equilibrium is
therefore Lyapounov stable (modulo $\SO(2)$) if $(n-1)\cos^2\theta_0
> (\ell-1)(n-\ell-1)$ for all $\ell=2,\dots,[n/2]$, that is if
$\cos^2\theta_0 > ([n/2]-1)(n-[n/2]-1)/(n-1) = \frac1{n-1}\left[\frac{n^2}4\right]-1$. This gives the
desired values.

We now turn to linear stability.  It follows from Proposition
\ref{diagonering} and the block diagonalization of
$\dd^2H_\xi|_\NN(x_e)$ that
$$
L_\NN=\diag\left( \left( \begin{array}{cc}
      0&-\lambda_1\\
      \lambda_1&0
\end{array}\right),
\left\lbrace \left( \begin{array}{cc}
      0&-\lambda^{(\ell)}_\phi\\
      \lambda^{(\ell)}_\theta&0
\end{array}\right),
\left( \begin{array}{cc}
    0&-\lambda^{(\ell)}_\phi\\
    \lambda^{(\ell)}_\theta&0
\end{array}\right)
\mid 2\leq \ell\leq [n/2]\right\rbrace^* \right)
$$
where the blocks are given up to a strictly positive scalar factor.
The eigenvalues of $L_\NN$ are therefore
$$
\pm \ii \lambda_1, \left\lbrace\pm
\ii\sqrt{\lambda^{(\ell)}_\theta\lambda^{(\ell)}_\phi}\mid
  2\leq \ell\leq [n/2]\right\rbrace,
$$
(up to a positive factor) and so the relative equilibrium is
linearly unstable if $\lambda^{(\ell)}_\theta<0$ for some $\ell$,
that is if
$$
\cos^2\theta_0 <\frac1{n-1}\left[\frac{n^2}4\right]-1.
$$
In particular this inequality is satisfied if $\theta_0=\pi/2$ and
$n>3$.

When the ring is equatorial, one has $\theta_0=\pi/2$ and $\mu=0$.
In particular $\lambda_1=0$.  This is because the symplectic slice
is smaller ($G_\mu=\SO(3)$ for $\mu=0$): it follows from Proposition
\ref{prop:tangentspace} that we have to remove the vectors $e_1$, $e_2$
from the basis for $\mu\neq 0$ (that is to remove $\lambda_1$ from
the previous eigenvalue study).  However, this does not change the
instability results, as the instability is primarily due to the
$\ell=[n/2]$ mode. It follows that the $\CC_{nv}$ equatorial
relative equilibria are linearly unstable for $n>3$, and Lyapounov
stable (modulo $\SO(3)$) for $n=3$.
\end{proof}

\subsection{Bifurcations}
The loss of stability of the ring is as usual accompanied by a bifurcation.  The proof above shows that the `critical mode' for stability is $\ell = [n/2]$. For $n \geq 7$ a ring is always unstable to this mode, while for $4 \leq n \leq 6$ the ring is stable when sufficiently close to the pole, and loses stability to this mode as it moves closer to the equator.  Here we describe the bifurcations that accompany this loss of stability.  

There is another bifurcation that occurs in all rings with $n\geq 6$, namely in the $\ell=2$ mode. Indeed, expression (\ref{eq:theta-evs for ring}) shows that only the modes $\ell=2$ and $\ell=[n/2]$ can satisfy $\lambda_\theta^{(\ell)}=0$ (the former at $\cos^2\theta_0=\frac{n-3}{n-1}$). We are only considering `relative steady-state' bifurcations, so where the eigenvalues pass through zero.  

The dihedral symmetry in each mode dictates the type of bifurcation to expect: if $\ell=n/2$ so $\dim V_\ell=2$, one has a $\Z_2$-pitchfork bifurcation involving a loss of $\Z_2$-symmetry. On the other hand if $\dim V_\ell=4$ then there is a pitchfork bifurcation with dihedral symmetry. See Section \ref{sec:symmetry+bifurcations} for a discussion of these.  Here we describe briefly the bifurcations arising for low values of $n$: the pattern continues for larger $n$.

\begin{figure}[t]
\begin{center}
\psset{unit=1.8,dotsep=2pt,linewidth=0.5pt,dotsize=4pt}
\subfigure[$n=4,\:\ell=2$]{%
\begin{pspicture}(-1,-1)(1,1) %show(4,2,[0.2,0])
 \psline[linestyle=dotted](0., 1.)(-1., 0.)(0., -1.)(1., 0.)(0., 1.)
 \psline[showpoints=true](0., .8)(-1.2, 0.)(0., -.8)(1.2, 0.)(0., .8)
 \psline[linecolor=lightgray](-0.3,0)(0.3,0)
 \psline[linecolor=lightgray](0,-0.3)(0,0.3)
\end{pspicture}\label{fig:perturbations:4,2}}\qquad\qquad
\subfigure[$n=5,\:\ell=2$]{\label{fig:perturbations:5,2}%
\begin{pspicture}(-1,-1)(1,1) %show(5,2,[.2, 0, 0, .1])
 \psline[linestyle=dotted](1., 0.)(.306, .952)(-.809, .588)(-.809, -.588)(.306, -.952)(1., 0.)
 \psline[showpoints=true](.312, .780)(-.914, .546)(-.914, -.546)(.312, -.780)(1.2, 0.)(.312, .780)
 \psline[linecolor=lightgray](-0.3,0)(0.3,0)
\end{pspicture}}
\qquad\quad
\subfigure[$n=6,\:\ell=3$]{%
\begin{pspicture}(-1,-1)(1,1) %show(6,3,0.1,0)
 \psline[linestyle=dotted](.500, .865)(-.500, .865)(-1., 0.)(-.500, -.865)(.500, -.865)(1., 0.)(.500, .865)
 \psline[showpoints=true](.450, .778)(-.550, .952)(-.9, 0.)(-.550, -.952)(.450, -.778)(1.1, 0.)(.450, .778)
 \psline[linecolor=lightgray](-0.3,0)(0.3,0)
 \psline[linecolor=lightgray](-.1500, .2598)(.1500, -.2598)
 \psline[linecolor=lightgray](-.1500, -.2598)(.1500, .2598)
\end{pspicture}\label{fig:perturbations:6,3}}
\caption{Polar view of the $\ell=[n/2]$ bifurcations for $n=4,5,6$. The dotted lines represent the `central' relative equilibrium with $\DD_n=\CC_{nv}$ symmetry, while the dots are the vortices of the (stable) bifurcating relative equilibrium with lower symmetry;  the grey lines in the centre of each represent the lines of reflexion.  The figures would be rotating in time.}
\label{fig:bifurcating re for ell = n/2}
\end{center}
\end{figure}
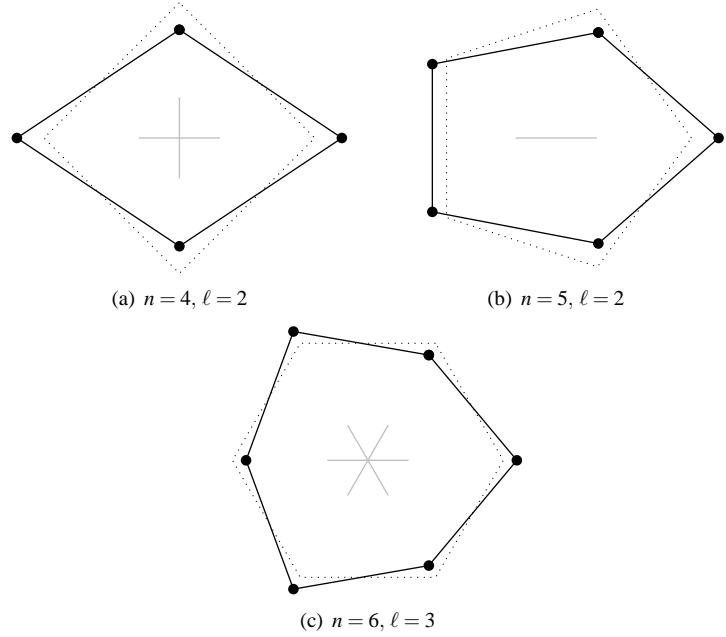

%%%%%%%%%%%%%%%%%%%%%%%%%%%%%%%%%%%%%%%%%%%%%%%%%%%%%%%
\begin{figure}[t]
\psset{unit=1.8,dotsep=2pt,linewidth=0.5pt,dotsize=4pt}
\begin{center}
\subfigure[$n=6,\:\ell=2$ (Fix $\kappa$)]{%
\begin{pspicture}(-1,-1)(1,1) %show(6, 2, [.1, 0, 0, .05])
 \psline[linestyle=dotted](.500, .865)(-.500, .865)(-1., 0.)(-.500, -.865)(.500, -.865)(1., 0.)(.500, .865)
 \psline[showpoints=true](.512, .800)(-.512, .800)(-1.1, 0.)(-.512, -.800)(.512, -.800)(1.1, 0.)(.512, .800)
 \psline[linecolor=lightgray](-0.3,0)(0.3,0)
 \psline[linecolor=lightgray](0,-0.3)(0,0.3)
\end{pspicture}\label{fig:perturbations:6,2}}
\qquad\qquad
\subfigure[$n=7,\:\ell=2$]{%
\begin{pspicture}(-1,-1)(1,1) %show(7, 2, [.1, 0, 0, 0.05])
 \psline[linestyle=dotted](.623, .782)(-.219, .976)(-.901, .434)(-.901, -.434)(-.219, -.976)(.623, -.782)(1., 0.)(.623, .782)
 \psline[showpoints=true](.647, .735)(-.220, .883)(-.972, .425)(-.972, -.425)(-.220, -.883)(.647, -.735)(1.1, 0.)(.647, .735)
 \psline[linecolor=lightgray](-0.3,0)(0.3,0)
\end{pspicture}}

\subfigure[$n=8,\:\ell=2$ (Fix\,$\kappa$)]{%
\begin{pspicture}(-1,-1)(1,1) % show(8, 2, [.1, 0, 0, 0.05])
 \psline[linestyle=dotted](.705, .705)(0., 1.)(-.705, .705)(-1., 0.)(-.705, -.705)(0., -1.)(.705, -.705)(1., 0.)(.705, .705)
 \psline[showpoints=true](.740, .670)(0., .9)(-.740, .670)(-1.1, 0.)(-.740, -.670)(0., -.9)(.740, -.670)(1.1, 0.)(.740, .670)
 \psline[linecolor=lightgray](-0.3,0)(0.3,0)
 \psline[linecolor=lightgray](0,-0.3)(0,0.3)
\end{pspicture}\label{fig:perturbations:8,2a}}
\qquad\qquad
\subfigure[$n=8,\:\ell=2$ (Fix\,$\kappa'$)]{%
\begin{pspicture}(-1,-1)(1,1) 
 \psline[linestyle=dotted,linewidth=1pt](.705, .705)(0., 1.)(-.705, .705)(-1., 0.)(-.705, -.705)(0., -1.)(.705, -.705)(1., 0.)(.705, .705)
 \psline[showpoints=true](.729, .779)(-0.0352, .930)(-.631, .681)(-1.07, 0.0352)(-.729, -.779)(0.0352, -.930)(.631, -.681)(1.07, -0.0352)(.729, .779)
 \psline[linecolor=lightgray](-.277, -.115)(.277, .115)
 \psline[linecolor=lightgray](.115,-.277)(-.115,.277)
\end{pspicture}\label{fig:perturbations:8,2b}}
\caption{Perturbations of the $n$-ring in the mode $\ell=2$. These configurations are invariant under a subgroup isomorphic to $\Z_2$ or $\Z_2\times\Z_2$ according to the parity of $n$, and the grey lines in the centre of each represent the lines of reflexion. The dotted figures are the regular $n$-gons.}
\label{fig:ell=2 bifurcations}
\end{center}
\end{figure}
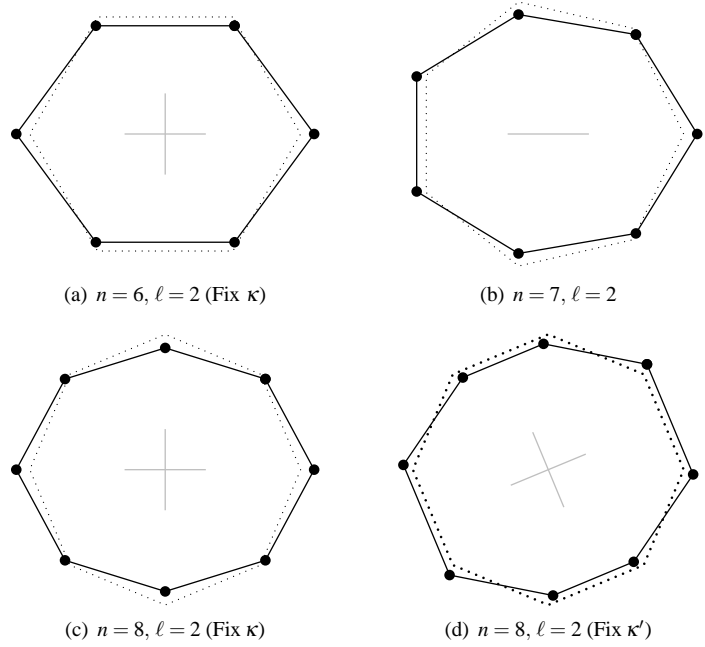
%%%%%%%%%%%%%%%%%%%%%%%%%%%%%%%%%%%%%%%%%%%%%%%%%%%%%%%

\begin{description}
\item[n = 3] There is no bifurcation: the relative equilibrium is always Lyapounov stable (relative to $\SO(2)$ if $\mu\neq0$ and to $\SO(3)$ if $\mu=0$).

\item[n = 4]
The only bifurcation occurring is when $\cos^2\theta_0=1/3$ where the stability of the ring changes. The mode $\ell=2$ which is involved in the bifurcation is of dimension 2, spanned by $\alpha_\theta^{(2)},\,\alpha_\phi^{(2)}$. It is the eigenvalue corresponding to $\alpha_\theta^{(2)}$ that passes through zero, which tells us the symmetry of the bifurcating relative equilibria: they consists of a pair of 2-rings on different latitudes, denoted $\CC_{2v}(R,R')$, as shown in Figure~\ref{fig:bifurcating re for ell = n/2}(a). 

\item[n = 5]
Again the only bifurcation occurs in the $\ell=2$ mode when $\theta_0=\pi/4$ and $3\pi/4$. This mode is 4-dimensional, with an action of $\DD_5$. For $\theta$ close to these values (and either above or below: we have not calculated the higher order terms to determine this) there are 5 saddle points (unstable relative equilibria) and 5 other critical points which will be elliptic (and possibly nonlinearly stable, depending whether the bifurcation is sub- or super-critical).  See Figure~\ref{fig:bifurcating re for ell = n/2}\,(b) for a representation of a typical bifurcating solution.

\item[n = 6]
In this case there are two bifurcations as the co-latitude increases. The first occurs at $\cos^2\theta_0=4/5$ in the $\ell=3$ mode where the regular hexagon loses stability. Here there is a usual $\Z_2$-pitchfork bifurcation, resulting in a pair of staggered 3-rings (that is, $C_{3v}(R,R')$), which will be elliptic; see Figure~\ref{fig:bifurcating re for ell = n/2}\,(c) for a typical bifurcating configuration. The second bifurcation occurs at $\cos^2\theta_0=3/5$ in the $\ell=2$ mode; the effective action will be of $\DD_3$ and so if it is generic one expects a transcritical bifurcation. Figure~\ref{fig:ell=2 bifurcations}\,(a) shows a typical bifurcating configuration in this mode.  Between the two bifurcation values of $\theta_0$, the reduced system has a pair of (double) imaginary eigenvalues and a pair of real eigenvalues. Closer to the equator, all the eigenvalues are real.

\item[n = 7]
Here there is just one bifurcation (the ring is always unstable due to the mode $\ell=3$, where the linear system has real eigenvalues).  It occurs at $\cos^2\theta_0=2/3$ in the 4-dimensional $\ell=2$ mode. This is a $\DD_7$-pitchfork, which is a standard pitchfork, and analogous to that shown (for $n=4$) in Figure~\ref{fig:dihedral bifurcations}\,(b). Whether the bifurcating solutions exist for $\theta>\theta_0$ or for $\theta<\theta_0$ depends on the higher order terms, and we have not analysed these.

\item[n = 8]
Again there is just one bifurcation, occurring at $\cos^2\theta_0=5/7$.  In this case, since again $\ell=2$ the effective action on $V_2$ is one of $\DD_4$. Depending on the higher order terms, the bifurcation is either transcritical, as in Figure~\ref{fig:dihedral bifurcations}\,(a), or standard pitchfork as in Figure~\ref{fig:dihedral bifurcations}\,(b). (They are not sub- or super-critical as the central equilibrium is unstable throughout the bifurcation.)  One type of bifurcating solution will be as depicted in Figure~\ref{fig:dihedral bifurcations}\,(c) while the other will be as in Figure~\ref{fig:dihedral bifurcations}\,(d). 
\end{description}

Similar conclusions can be made for $n\geq9$, where all $\ell=2$ bifurcations will be standard pitchfork, rather than transcritical. 
Further calculations involving the higher order terms in the bifurcations occurring here can be found in \cite{MT11}.

%%%%%%%%%%%%%%%%%%%%%%%%%%%%%%%%%%%%%%%%%%%%%%%%%%%%%%%%%%%%%%%%%%%%%
\section{A ring with a polar vortex: $\CC_{nv}(R,p)$}
\label{sec:Rp-stability}

A configuration consisting of a ring of identical vortices placed regularly around a line of latitude, together with a single pole either at the North or South pole is always a relative equilibrium, rotating steadily about the `vertical' axis. We assume that the polar vortex lies at the North pole and its vorticity is $\kappa$, while the remaining $n$ vortices are all identical with vortex strength 1 and lie in a regular ring on a fixed circle of co-latitude $\theta_0$. 

This relative equilibrium is of symmetry type $\CC_{nv}(R,p)$ and denoted $x_e$. Its momentum is $\Phi(x_e)=(0,0,\mu)$, where $\mu=\kappa+n\cos\theta_0$.

The Hamiltonian is given by
$$H \ =\  H_r + H_p$$
where $H_r$ is the contribution to the total Hamiltonian from the
interactions within the ring (\ref{eq:ring Hamiltonian}), and
$$
H_p(x,y,\theta_i,\phi_i) = -\kappa\sum_{j=1}^n \ln\left(1 -
  x\sin\theta_j\cos\phi_j-y\sin\theta_j\sin\phi_j -
  \sqrt{1-x^2-y^2}\cos\theta_j\right),
$$
is the Hamiltonian responsible for the interaction of the polar
vortex and the ring.

In this case
$$
H_\xi \ =\  H - \xi\left(\sum_j\cos\theta_j +
  \kappa\sqrt{1-x^2-y^2}\right).
$$
The $\CC_{nv}(R,p)$ relative equilibrium at $x=y=0,\
\theta_j=\theta_0,\ \phi_j=2\pi j/n$ has angular velocity
$$
\xi \ =\ \frac{(n-1)\cos\theta_0 +
  \kappa(1+\cos\theta_0)}{\sin^2\theta_0} = \frac{\mu+(\kappa-1)\cos\theta_0}{\sin^2\theta_0}
$$
since $\frac{\partial H_\xi}{\partial\theta_j}(x_e) =
-\frac{(n-1)\cos\theta_0+\kappa(1+\cos\theta_0)-\xi\sin^2\theta_0}{\sin\theta_0}$
must vanish.

\bigskip

The second derivatives at the relative equilibrium of $H$ can be
derived from those for $H_r$ given in the previous section, together
with (for $n>2$):
$$
\left.
  \begin{array}{rclcrcl}
    \frac{\partial^2H_p}{\partial\theta_j^2} &=& \frac{\kappa}{1-\cos\theta_0}
    &\quad&
    \frac{\partial^2H_p}{\partial x^2} &=& \frac{n\kappa}{2} \;=\;
    \frac{\partial^2H_p}{\partial y^2}  \\[6pt]
    \frac{\partial^2H_p}{\partial x\partial \theta_j} &=&
    -\frac{\kappa\cos(2\pi j/n)}{1-\cos\theta_0} &&
    \frac{\partial^2H_p}{\partial y\partial \theta_j} &=&
    -\frac{\kappa\sin(2\pi j/n)}{1-\cos\theta_0}\\[6pt]
    \frac{\partial^2H_p}{\partial x \partial\phi_j} &=&
    -\frac{\kappa\sin\theta_0\sin(2\pi j/n)}{1-\cos\theta_0} &&
    \frac{\partial^2H_p}{\partial y \partial\phi_j} &=&
    \frac{\kappa\sin\theta_0\cos(2\pi j/n)}{1-\cos\theta_0},
  \end{array}\right.
$$
while the other second derivatives all vanish.  Here we have used
that $\sum\cos^2(2\pi j/n) = n/2$ for $n>2$.  For  $n=2$ this sum is
$2$ and one obtains
$$
\textstyle \frac{\partial^2H_p}{\partial x^2} =
\frac{2\kappa}{1-\cos\theta_0} \qquad\qquad
\frac{\partial^2H_p}{\partial y^2} =
-\frac{2\kappa\cos\theta_0}{1-\cos\theta_0}.
$$

\subsection{Symplectic slice}
The following proposition gives the symmetry adapted basis for
$\CC_{nv}(R,p)$ relative equilibria.

\begin{proposition} \label{diag} %
With the symplectic slice decomposition described in Proposition
\ref{prop:isotypic}, that is $\NN = \bigoplus_{\ell=1}^{[n/2]}
V_\ell$, and for $n\geq 3$ and $\mu\neq 0$ a basis for $V_1$ is
$\{e_1,e_2,e_3,e_4\}$, where
$$
\left\{\begin{array}{lll}
 e_1 &=& \cos\theta_0\:\beta^{(1)}_{\phi} + \sin\theta_0\: \alpha^{(1)}_{\theta}\\[6pt]
 e_2 &=& \kappa\: \beta^{(1)}_{\phi} + \frac{n}2 \sin\theta_0\:\delta x\\[6pt]
 e_3 &=&  \cos\theta_0\: \alpha^{(1)}_{\phi} - \sin\theta_0\:\beta^{(1)}_{\theta}\\[6pt]
 e_4 &=&  \kappa\: \alpha^{(1)}_{\phi} - \frac{n}2\sin\theta_0\:\delta y\,.
\end{array}\right.
$$
With respect to the resulting basis for $\NN$, the Hessian
$\dd^2H_\xi|_\NN(x_e)$ block diagonalizes into two $2\times 2$ blocks
(corresponding to $\ell=1$) and $(2n-6)$ $1\times 1$ blocks, and
$L_\NN$ block diagonalizes into one $4\times 4$ block (for $\ell=1$)
and the remaining are $2\times 2$ blocks. %

When $\mu=0$, $V_1$ drops dimension by 2, and in particular the
vectors $e_1+\frac2n e_2$ and $e_3 + \frac2n e_4$ lie in the tangent
space to the group orbit so one can take $V_1=\<e_1,e_3\m$, and in this basis the
Hessian restricted to $V_1$ is diagonal.

The case where $n=2$ is discussed  below.
\end{proposition}

For later use we record the symplectic form on the space $V_1$, the values follow from (\ref{eq:omega on alphas,betas}):
\begin{equation}\label{eq:omega on V1(R,p)}
\begin{array}{lcl}
 \omega(e_1,\,e_2) =0 &\qquad& \omega(e_3,\,e_4)=0\\[6pt]
 \omega(e_1, e_3) =  n\cos\theta_0\sin^2\theta_0 & \qquad&
 \omega(e_1, e_4) =   \half n\,\kappa  \sin^2\theta_0 \\[6pt]
 \omega(e_2, e_3) =  \half n\,\kappa  \sin^2\theta_0 & \qquad&
 \omega(e_2, e_4) =   -\frac14 n^2\kappa\sin^2\theta_0.%\\[6pt]
\end{array}
\end{equation}
It follows that $\left<e_1,\,e_2\right>$ and $\left<e_3,\,e_4\right>$ are invariant Lagrangian subspaces (and the representation is the sum of an irreducible and its dual).

\begin{proof}
The proof is similar to the proof of Proposition~\ref{diagonering}, and we omit the details.
\end{proof}

Using equations (\ref{eq:action on alphas,betas}) one sees that the
subspaces $\left<e_1,\,e_3\right>$ and $\left<e_2,\,e_4\right>$ are
invariant and by (\ref{eq:omega on V1(R,p)} above) they are
symplectic. It follows that $V_1$ is the direct sum of two
representations of \emph{complex dual} type \cite{MRS88} (they are
dual as on one the orientations defined by $\omega$ and by $\CC_n$
coincide, while for the other they are opposite).

\subsection{Stability}

The block diagonalizations of the proposition enable us to prove the
following stability theorem for $n \geq 4$, illustrated by Figure
\ref{fig:CRp}. The cases $n=2$ and $3$ are treated afterwards (and
illustrated in Figure~\ref{fig:CRp2/3}).

\begin{theorem}
\label{ring+pole/stab} A $\CC_{nv}(R,p)$ relative equilibrium with
$n\geq 4$ and $\mu\neq 0$ is \\
(i) Lyapounov stable if
$$
\kappa > \kappa_0\ \ \mbox{and}\ \ \kappa
(\kappa+n\cos\theta_0)(a\kappa-b) < 0,
$$
(ii) spectrally unstable if and only if
$$
\kappa < \kappa_0\ \ \mbox{or}\ \ 8a\kappa >
(n\sin^2\theta_0+4(n-1)\cos\theta_0)^2,
$$
where
\begin{equation}\label{eq:ring+pole/stab}
\begin{array}{lll}
a&=&(1+\cos\theta_0)^2(n\cos\theta_0-n+2)\\[6pt]
b&=& (n-1)\cos\theta_0\; (n\sin^2\theta_0+2(n-1)\cos\theta_0)\\[6pt]
\kappa_0&=&\left(\left[\frac{(n-2)^2}4\right] - (n-1)\cos^2\theta_0\right)/(1+\cos\theta_0)^2.
\end{array}
\end{equation}
\end{theorem}

As will be seen in the proof, the conditions involving $\kappa_0$
arise from the $\ell=[n/2]$ mode, while the others arise from the
$\ell=1$ mode.

\begin{proof}
(i) We first study the Lyapounov stability.  Following the beginning
of the proof of Theorem~\ref{polvani}, we obtain from Proposition
\ref{diag} that
  $$
  \dd^2H_\xi|_\NN(x_e)=\diag(A,A,D)
  $$
where $D=\diag(\lbrace\lambda^{(\ell)}_\theta,\lambda^{(\ell)}_\phi,
\lambda^{(\ell)}_\theta,\lambda^{(\ell)}_\phi \mid 2\leq \ell\leq
[n/2]\rbrace^* )$,
  $$
  A=\left( \begin{array}{cc}
      q_{11}&q_{12}\\
      q_{12}&q_{22}
  \end{array}\right)
  $$
and
$$
\begin{array}{lll}
\lambda^{(\ell)}_\theta&=&\dd^2H_\xi(x_e)\cdot (\alpha^{(\ell)}_{\theta},\alpha^{(\ell)}_{\theta})\\
\lambda^{(\ell)}_\phi&=&\dd^2H_\xi(x_e)\cdot (\alpha^{(\ell)}_{\phi},\alpha^{(\ell)}_{\phi})\\
q_{11}&=&\dd^2H_\xi(x_e)\cdot (e_1,e_1)\\
q_{12}&=&\dd^2H_\xi(x_e)\cdot (e_1,e_2)\\
q_{22}&=&\dd^2H_\xi(x_e)\cdot (e_2,e_2).
\end{array}
$$
Note that $D$ exists only for $n\geq 4$.  From the previous section
one has $\lambda^{(\ell)}_\phi=n\ell(n-\ell)/2$ and some additional
computations give
$$
\lambda^{(\ell)}_\theta=\frac{n}{2\sin^2\theta_0}\left[
  -(\ell-1)(n-\ell-1)+(n-1)\cos^2\theta_0+\kappa(1+\cos\theta_0)^2
\right].
$$
The eigenvalues $\lambda^{(\ell)}_\phi$ are all positive, thus $D$
is definite if and only if, for all $\ell=2,\dots,[n/2]$,
$$ \kappa(1+\cos\theta_0)^2 > (\ell-1)(n-\ell-1) - (n-1)\cos^2\theta_0.
$$
This holds if and only if the inequality is satisfied for
$\ell=[n/2]$, which is equivalent to $\kappa>\kappa_0$.

The relative equilibrium is therefore Lyapounov stable if $A$ is
positive definite, that is if $q_{11}q_{22}-q_{12}^2 >0$ and $q_{11}
> 0$.  Some lengthy computations give
$$
\begin{array}{rcl}
  q_{11}q_{22}-q_{12}^2&=& -\frac18 n^2 \kappa(\kappa+n\cos\theta_0)(a\kappa-b)\\
  q_{11}&=& \half n\kappa(1+\cos\theta_0)^2+ n(n-1)\cos^2\theta_0
\end{array}
$$
where $a,b$ are given in the theorem. Now, if $\kappa
> \kappa_0$, then $q_{11} > 0$: indeed
$$q_{11} = \half n(1+\cos\theta_0)^2(\kappa-\kappa_0) + \half
n(n-1)\cos^2\theta_0+\sfrac{n}{2}\left[\sfrac{(n-2)^2}{4}\right].$$
We proved therefore that this \re\ is Lyapounov stable if
$\kappa > \kappa_0$ and $\kappa (\kappa+n\cos\theta_0)
(a\kappa-b)<0$.

(ii) We now study the spectral stability of the relative
equilibrium. It follows from Proposition~\ref{diag} and the block
diagonalization of $\dd^2H_\xi|_\NN(x_e)$ that
$$
L_\NN=\diag\left(A_L, \left\lbrace \left( \begin{array}{cc}
        0&-\lambda^{(\ell)}_\phi\\
        \lambda^{(\ell)}_\theta&0
\end{array}\right),
\left( \begin{array}{cc}
    0&-\lambda^{(\ell)}_\phi\\
    \lambda^{(\ell)}_\theta&0
\end{array}\right)
\mid 2\leq \ell\leq [n/2]\right\rbrace^* \right)
$$
where the blocks are given up to a positive scalar factor and
$$
A_L=\left(
\begin{array}{cccc}
    0 & 0 & a & b\\
    0 & 0 & c & d\\
    -a & -b & 0 & 0\\
    -c & -d & 0 & 0
\end{array}\right),\ \ \mbox{where} \ \left\lbrace \begin{array}{ccc}
a&=&\beta q_{11}-\gamma q_{12}\\
b&=&\beta q_{12}-\gamma q_{22}\\
c&=&\alpha q_{12}-\gamma q_{11}\\
d&=&\alpha q_{22}-\gamma q_{12}
\end{array}\right.
$$
and
$$
\begin{array}{lll}
 \alpha&=&-\frac{4\cos\theta_0}{n(\kappa+n\cos\theta_0)\sin^2\theta_0} \\[6pt]
 \beta&=& \frac{1}{(\kappa+n\cos\theta_0)\sin^2\theta_0} \\[6pt]
 \gamma &=& -\frac{2}{n(\kappa+n\cos\theta_0)\sin^2\theta_0}
\end{array}
$$
These ($\alpha,\beta$ and $\gamma$) arise from $L=J \dd^2H_\xi$, where
$J$ is the inverse of the matrix of the symplectic form whose
coefficients are given in (\ref{eq:omega on V1(R,p)}).

The eigenvalues (up to a positive factor) of $L_\NN$ are therefore
$$
\pm\frac{1}{\sqrt{2}}\sqrt{\sigma\pm\sqrt{\nu}\;} ,\left\lbrace\pm
  \ii\sqrt{\lambda^{(\ell)}_\theta\lambda^{(\ell)}_\phi}\mid 2\leq
  \ell\leq [n/2]\right\rbrace^*
$$
where $\nu=a^4+4a^2bc-2a^2\dd^2+4bc\dd^2+d^4+8adbc$ and
$\sigma=-a^2-2bc-\dd^2$. The eigenvalues $\pm
i\sqrt{\lambda^{(\ell)}_\theta\lambda^{(\ell)}_\phi}$ are all purely
imaginary if and only if $\kappa > \kappa_0$. After some lengthy but
straightforward computations we obtain that
$$
\begin{array}{lll}
\sigma&=& \frac{1}{8\sin^2\theta_0}\left[4(1+z)^2(nz-n+2)\kappa  \right.\\
   & & \qquad\qquad \left. - n^2z^4 + 4n(n-1)z^3 - 2(3n^2-8n+4)z^2 - 4n(n-1)z -
   n^2\right]\\[12pt]
\nu&=&-\frac9{16\sin^4\theta_0} \left[8(1+z)^2(nz-n+2)\kappa
-\left(n(1-z^2)+4(n-1)z\right)^2\right]
\end{array}
$$
where $z=\cos\theta_0$. One can check that if $\nu\geq0$, then
$\sqrt{\nu}+\sigma\leq0$ and the eigenvalues are purely imaginary.
If $\nu<0$, then the eigenvalues have a non-zero real part. Thus the
eigenvalues $\pm\sqrt{\sigma\pm\sqrt{\nu}\:}$ are purely imaginary
if and only if $\nu\geq0$ which is equivalent to $8a\kappa \leq
(n\sin^2\theta_0+4(n-1)\cos\theta_0)^2$.
\end{proof}

A spectrally stable relative equilibrium for which the Hessian
$\dd^2H_\xi|_\NN(x_e)$ is not definite is said to be \emph{elliptic}.
Note that in principle an elliptic relative equilibrium may be
Lyapounov stable, but if there are more than 4 vortices then it is
expected to be unstable as a result of Arnold diffusion. Moreover an
elliptic relative equilibrium typically becomes linearly unstable
when some dissipation is added to the system \cite{DR02}; however
adding dissipation to the point vortex system would have more
profound effects, such as spreading of vorticity into vortex
patches.

\begin{corollary}
  A $\CC_{nv}(R,p)$ relative equilibrium with $n\geq 4$ and $\mu\neq
  0$ is elliptic if and only if $$
  \kappa \geq \kappa_0,\ \
  \kappa (\kappa+n\cos\theta_0)(a\kappa-b) \geq 0\ \
  \mbox{and}\ \ 8a\kappa \leq
  (n\sin^2\theta_0+4(n-1)\cos\theta_0)^2,
  $$
where $a,b$ and $\kappa_0$ are given in (\ref{eq:ring+pole/stab}).
\end{corollary}

\paragraph{Discussion of results for $n \geq 4$} 
See Figure~\ref{fig:CRp}.

\begin{figure}[htbp] %% ring + pole -- n = 5 & 8
\begin{center}\psset{unit=0.85}
\begin{pspicture}(-1,-1.5)(12,7.3)  %% n = 8
  \psframe(-1,-1.3)(12,7.3)
  \rput(11.3,7){{\large n = 8}}
  \psline(10.6,7.3)(10.6,6.6)(12,6.6)
%% rescale x to degrees and y so that -2 .. 20 fits
  \psset{xunit=0.06, yunit=.3}
%% stability regions:
  \pscustom[fillstyle=solid,fillcolor=LyapounovColour]{ %% upper stability region
    \psplot{25.30756309}{0}{7 x cos mul 8 x sin 2 exp mul 14 x cos mul add mul %
        8 x cos mul 6 sub 1 x cos add 2 exp mul div}
    \psline[linewidth=0.1pt](0,12.25)(0,0.5)
    \psplot{0}{111.3987128}{9 x cos 2 exp 7 mul sub x cos 1 add 2 exp div} %kappa=kappa_0
    }
  \psplot[linewidth=.08,linecolor=EllipticColour]{0}{25.30756309}{7 x cos mul 8 x sin 2 exp mul 14 x cos mul add mul %
        8 x cos mul 6 sub 1 x cos add 2 exp mul div} %% left eliptic sliver
%% eigenvalue diagrams:
  \psline[linewidth=0.2pt](120,10)(140,10)
   \psline[linewidth=0.2pt](130,8)(130,12)
   \psdots[dotstyle=x](125,10)(135,10)
  \psline[linewidth=0.2pt](5,18)(15,18)
   \psline[linewidth=0.2pt](10,16.5)(10,19.5)
   \psdots(7,17)(7,19)(13,17)(13,19)
  \psline[linewidth=0.2pt](35,13)(45,13)
   \psline[linewidth=0.2pt](40,11.5)(40,14.5)
   \psdots(40,12)(40,12.5)(40,13.5)(40,14)
%% axes:
  \psaxes[linewidth=0.5pt,labels=y,ticks=y,Dy=5]{->}(0,0)(0,-2)(180,21)
  \rput(185,0){$\theta_0$}  \rput(-3,22){$\kappa$}
  \rput(90,-0.8){\small $\frac{\pi}{2}$}
  \rput(180,-0.8){\small $\pi$}
  \rput(100,8.7){$\kappa=\kappa_0$}
%% curves xi = 0 and mu = 0:
%  \psplot[linestyle=dotted]{0}{139}{x cos 1 add -1 exp x cos mul -7 mul} \rput(145,22){$\xi=0$}
  \psplot[linestyle=dashed,linewidth=0.6pt]{75.52248783}{182}{x cos -8 mul} \rput(190,8){$\mu=0$}
%% hemispheres:
 \rput(45,-3){\it Northern hemisphere} %
 \rput(90,-3){$\vert$}
 \rput(135,-3){\it Southern hemisphere}
\end{pspicture}
\begin{pspicture}(-1,-4.5)(12,5)  %% n = 5
  \psframe(-1,-4.3)(12,4.3)
  \rput(11.35,4.02){{\large n = 5}}
  \psline(10.6,4.3)(10.6,3.7)(12,3.7)
%% rescale x to degrees and y so that -10 .. 10 fits
  \psset{xunit=0.06, yunit=.35}
%% stability regions:
  \pscustom[fillstyle=solid,fillcolor=LyapounovColour]{ %% upper stability region
    \psplot{40.46315127}{0}{4 x cos mul 5 x sin 2 exp mul 8 x cos mul add mul %
        5 x cos mul 3 sub 1 x cos add 2 exp mul div}
    \psline[linewidth=0.1pt](0,4)(0,0)(45,0)
    \psplot{45}{128.6395024}{2 x cos 2 exp 4 mul sub x cos 1 add 2 exp div} %kappa=kappa_0
    \psplot{128.6395024}{180}{-5 x cos mul} % mu = 0
    \gsave
    \psline(180,5)(180,10)
    \fill[fillstyle=solid]
    \grestore
    }
  \pscustom[fillstyle=solid,fillcolor=LyapounovColour]{ %% lower stability region
    \psline[linewidth=0.1pt](180,0)(135,0)
  \psplot{135}{141.6975218}{2 x cos 2 exp 4 mul sub x cos 1 add 2 exp div} %kappa=kappa_0
    \gsave
    \psline(141.6975218,-10)(180,-10)
    \fill[fillstyle=solid]
    \grestore
    }
  \pscustom[fillstyle=solid,fillcolor=EllipticColour]{ %% right elliptic region
    \psline(180,0)(135,0)
    \psplot{135}{128.6395024}{2 x cos 2 exp 4 mul sub x cos 1 add 2 exp div} %kappa=kappa_0
    \psplot{128.6395024}{180}{-5 x cos mul}
    }
  \pscustom[fillstyle=solid,fillcolor=EllipticColour]{ %% lower left elliptic region
    \psline(0,0)(0,-0.5)
    \psplot{0}{45}{2 x cos 2 exp 4 mul sub x cos 1 add 2 exp div} %kappa=kappa_0
    }
  \psplot[linewidth=.08,linecolor=EllipticColour]{40.46315127}{0}{4 x cos mul 5 x sin 2 exp mul 8 x cos mul add mul %
        5 x cos mul 3 sub 1 x cos add 2 exp mul div 0.1 add} %% upper left sliver
%% eigenvalue diagrams:
  \psline[linewidth=0.2pt](80,-6)(100,-6) \psline[linewidth=0.2pt](90,-4.5)(90,-7.5) %bottom middle
   \psdots[dotstyle=x](85,-6)(95,-6)
  \psline[linewidth=0.2pt](10,8)(30,8)\psline[linewidth=0.2pt](20,6)(20,10) % top left
   \psdots(15,7)(15,9)(25,7)(25,9) %\ell=1
   \psdots[dotstyle=x](20,6.5)(20,9.5) %\ell=2
%% axes:
  \psaxes[linewidth=0.5pt,labels=y,ticks=y,Dy=2]{->}(0,0)(0,-10.2)(180,11)
  \rput(185,0){$\theta_0$}  \rput(-3,11.5){$\kappa$}
  \rput(90,-0.8){\small $\frac{\pi}{2}$}
  \rput(180,-0.8){\small $\pi$}
  \rput(131,-8){$\kappa=\kappa_0$}
%% curves xi = 0 and mu = 0:
  \psplot[linestyle=dotted]{0}{137}{x cos 1 add -1 exp x cos mul -4 mul} \rput(145,11){$\xi=0$}
  \psplot[linestyle=dashed,linewidth=0.6pt]{0}{182}{x cos -5 mul} \rput(190,5){$\mu=0$}
%% hemispheres:
 \rput(45,-11){\it Northern hemisphere} %
 \rput(90,-11){$\vert$}
 \rput(135,-11){\it Southern hemisphere}
%\rput(72,-4){X}
\end{pspicture}
\caption{Bifurcation diagrams for  $\CC_{nv}(R,p)$ relative equilibria.
See Fig.~\ref{fig:CRp2/3} for the key. 
The bifurcation diagrams for $n\geq7$ are similar to that for
$n=8$, while those for $n=4$ and $6$ are similar to that for $n=5$.
The circles represent the eigenvalues of the mode $\ell=1$, while
the crosses represent those of the mode $[n/2]$.  
Notice the sliver of elliptic near the upper left hand corner of both diagrams: these are not drawn to scale as they are too
thin---cf.~Fig.~\ref{fig:CRp2/3}~($n=3$), where it is drawn to scale.
Stability is modulo rotations about the vertical axis 
(i.e., $\SO(2)$), or modulo all rotations if $\mu=0$. See text for
more details.   In the stable regions all eigenvalues are imaginary.}
\label{fig:CRp}
\end{center}
\end{figure}
%%%%%%%%%%%

\begin{itemize}
\item If the sign of the vorticity of the polar vortex is opposite
to that of the ring then there are stable configurations only for
$n\leq 6$. Conversely configurations with $n\leq 6$ and $\theta_0$
close to $\pi$, ie with the ring close to the opposite pole, are
Lyapounov stable for all $\kappa<0$.

\item The region of Lyapounov stability is larger when the
vorticities of the pole and the ring have the same sign
($\kappa>0$). The stability frontiers in the upper-left corners of
Figure~\ref{fig:CRp} go to infinity as $\theta_0$ approaches
$\arccos(1-2/n)$. It follows that for $n\geq 4$ and
$\theta_0>\arccos(1-2/n)$, the relative equilibria are Lyapounov
stable for all sufficiently large $\kappa$. Thus, a ring of vortices
is stabilized by a polar vortex with a sufficiently large vorticity
of the same sign as the vortices in the ring. Note that for $4 \leq
n \leq 6$ and $\kappa$ positive, but sufficiently small, a ring near
the opposite pole is only elliptic and may not be Lyapounov stable.

\item
The limiting stability results for $\theta_0=0$, ie when the ring is
close to the polar vortex, coincide with the stability of a planar
$n$-ring plus a central vortex, see \cite{CS99} and \cite{LP}. This
is also true for $n = 2$ and $n =3$.

\item
The stability boundary where $\kappa=\kappa_0$ corresponds to the
mode $\ell=[n/2]$ and is analogous to the stability boundary for a
single ring. When $n$ is even stability is probably lost through a
pitchfork bifurcation to a relative equilibrium of type
$\CC_{\frac{n}{2}v}(R,R',p)$ consisting of two staggered
$\frac{n}{2}$-rings and a pole as $(\kappa, \theta_0)$ passes
through this boundary. This is illustrated for $n=4$ in Figure~8 of
\cite{LMR01}, where branch (g) meets branch (d). When $n$ is odd
there is an analogous transcritical bifurcation to relative
equilibria with only a single reflexional symmetry which fixes two
vortices and permutes the others. These are denoted by
$\CC_{h}(\frac{n-1}{2}R,2E)$ in \cite{LMR01}.  A nice illustration
in the case $n=3$ can be found in Figure~8 of \cite{CMS03}.  We have
not checked the non-degeneracy conditions for these bifurcations.

\item
The $\ell=1$ mode is responsible for two types of bifurcation. See Section  Section \ref{sec:symmetry+bifurcations} for descriptions of the associated bifurcations.

Firstly, when the eigenvalues become zero, the kernel is a single
irreducible symplectic representation of complex type (it is a plane
with $\CC_n$ acting by rotations) so the eigenvalues pass through
zero and remain on the imaginary axis. This corresponds to a
transition from Lyapounov stable to elliptic (the very thin blue (or light grey) sliver in the top left-hand corners of the figures).

Secondly, at the other side of the blue (or light grey) sliver, the pair of
eigenvalues continues to move away from 0 along the imaginary axis
until it meets the other pair from the $\ell=1$ mode; these then
`collide' and leave the imaginary axis through a Hamiltonian Hopf
bifurcation. 

\item The transition from Lyapounov stable to elliptic when passing
across the curve where $\mu=0$ is a generic phenomenon, and can be
explained by the geometry of the reduced spaces, and occurs when 
the \re\ with zero momentum is in fact an orbit of equilibria. 
See \cite{Mo11} for details.

\item Finally we note that when $\kappa$ crosses zero, eigenvalues
change sign without strictly passing through zero due to the fact that the
symplectic form becomes degenerate for $\kappa=0$.

\item There are other bifurcations which do not involve loss of stability, as they occur in one particular mode, while the relative equilibrium is already unstable because of a different mode. For example, for $n=4$ and 5 (and perhaps others) there is an inverted parabola in the $\kappa<0$ region along which there is a Hamiltonian Hopf bifurcation in the $\ell=1$ mode (similar to that for $n=3$ below). 

\end{itemize}

%%%%%%%%%%%%
\paragraph{Discussion of the case $n=3$}
See Figure~\ref{fig:CRp2/3}.

\medskip

\noindent Suppose first that the relative equilibrium has non-zero
momentum.
For $n=3$ (so 4 vortices altogether), by the proof of Theorem
\ref{ring+pole/stab} we have $\dd^2H_\xi|_\NN(x_e)=\diag(A,A)$ and
$L_N=A_L$. Hence $\CC_{3v}(R,p)$ is Lyapounov stable if $\kappa
(\kappa+3\cos\theta_0)(a\kappa-b) < 0$, and spectrally unstable if
and only if
$$
8a\kappa > (3\sin^2\theta_0+8\cos\theta_0)^2,
$$
where $a=(3\cos\theta_0-1)(1+\cos\theta_0)^2$ as in Theorem
\ref{ring+pole/stab}. These results are illustrated in Figure
\ref{fig:CRp2/3}. 

\begin{itemize}
\item Notice that a polar vortex \emph{destabilizes} a
3-ring if either the polar vortex is in the same hemisphere as the
ring and has a sufficiently strong vorticity of the same sign as the
ring, or the polar vortex has the opposite sign vorticity and the
ring lies in an interval containing $\theta_0 = 2\pi/3$ that grows
as the magnitude of the polar vorticity increases.  These two
regions have a vertical asymptote at $\cos\theta_0 = 1/3$ determined
by the vanishing of $a$.

\item Outside these regions there is a patchwork of regimes in which the
relative equilibrium is either Lyapounov stable or elliptic.

\item The relative equilibria with $\mu=0$ have a 2-dimensional symplectic
slice, and are all Lyapounov stable. The transition point where
$\mu=\xi=0$, $\kappa=1$ and $\cos(\theta_0)=-1/3$ corresponds to the
stable equilibrium consisting of 4 identical vortices placed at the
vertices of a regular tetrahedron \cite{PM98,LMR01,Ku04}.

\item Comparing this case with the diagram for $n=5$ shows large regions of stability with $\kappa<0$ which are unstable for $n=5$. The instability for $n=5$ is due to the $\ell=2$ mode which is absent for $n=3$.
\end{itemize}

%%%%%%%%%%%%%%%%%%%%%%%%%%%%%%%%%%%%%%%%%%
\begin{figure}[htbp] %% ring + pole -- n = 2 & 3
\begin{center}\psset{xunit=0.9,yunit=0.85}
\begin{pspicture}(-1,-5.5)(11,5)  %% n = 3 :
  \psframe(-1,-5)(11,5)
  \rput(10.3,4.7){{\large n = 3}}
  \psline(9.7,5)(9.7,4.3)(11,4.3)
%% rescale x to almost radians:
  \psset{xunit=3}
%% stability regions:
  \pscustom[fillstyle=solid,fillcolor=LyapounovColour]{ %% upper stability region
    \psline(1.57,0)(0,0)(0,0.5)
    \psplot[xunit=0.0174]{0}{60.714}{x cos x cos 2 exp 3 mul x cos -4 mul add -3 add mul %
           x cos 3 mul -1 add 1 x  cos add 2 exp mul div -2 mul}
    \gsave
    \psline(3.14159,4)(3.14159,3)
    \psplot[xunit=0.0174]{180}{109.47}{x cos -3 mul}
    \psplot[xunit=0.0174]{109.47}{90}{x cos x cos 2 exp 3 mul x cos -4 mul add -3 add mul %
           x cos 3 mul -1 add 1 x  cos add 2 exp mul div -2 mul}
    \fill[fillstyle=solid]
    \grestore
    }
  \pscustom[fillstyle=solid,fillcolor=LyapounovColour]{ %% lower left stability region
    \psline[linewidth=0.1pt](0,-4)(0,-3)
    \psplot[xunit=0.0174,linestyle=dashed,linewidth=0.1pt]{0}{90}{x cos -3 mul} % from mu=0
    \psplot[xunit=0.0174]{90}{76.34547974}{x cos x cos 2 exp 3 mul x cos -4 mul add -3 add mul %
           x cos 3 mul -1 add 1 x  cos add 2 exp mul div -2 mul}
    }
  \pscustom[fillstyle=solid,fillcolor=LyapounovColour]{ %% central stability region
    \psplot[xunit=0.0174,linestyle=dashed,linewidth=0.1pt]{90}{109.47}{x cos -3 mul} % from mu=0
    \psplot[xunit=0.0174]{109.47}{122.356}{x cos x cos 2 exp 3 mul x cos -4 mul add -3 add mul %
           x cos 3 mul -1 add 1 x  cos add 2 exp mul div -2 mul}
    }
  \pscustom[fillstyle=solid,fillcolor=LyapounovColour]{ %% lower right stability region
    \psplot[xunit=0.0174]{131.81}{122.356}{x cos x cos 2 exp 3 mul x cos -4 mul add -3 add mul %
           x cos 3 mul -1 add 1 x  cos add 2 exp mul div -2 mul}
    \psline[xunit=0.0174](122.356,0)(180,0)
    \gsave
     \psline[xunit=0.0174](180,0)(180,-4)
     \fill[fillstyle=solid]
    \grestore
    }
  \pscustom[fillstyle=solid,fillcolor=EllipticColour]{ %% lower left elliptic region
    \psline[xunit=0.0174,linewidth=0.1pt](0,-3)(0,0)(90,0)
    \psplot[xunit=0.0174,linestyle=dashed,linewidth=0.1pt]{90}{0}{x cos -3 mul} % from mu=0
    }
  \pscustom[fillstyle=solid,fillcolor=EllipticColour]{ %% central elliptic region
    \psplot[xunit=0.0174,linestyle=dashed,linewidth=0.1pt]{90}{109.47}{x cos -3 mul} % from mu=0
    \psplot[xunit=0.0174]{109.47}{90}{x cos x cos 2 exp 3 mul x cos -4 mul add -3 add mul %
           x cos 3 mul -1 add 1 x  cos add 2 exp mul div -2 mul}
    }
  \pscustom[fillstyle=solid,fillcolor=EllipticColour]{ %% upper right elliptic region
    \psplot[xunit=0.0174]{122.356}{109.47}{x cos x cos 2 exp 3 mul x cos -4 mul add -3 add mul %
           x cos 3 mul -1 add 1 x  cos add 2 exp mul div -2 mul}
    \psplot[xunit=0.0174,linestyle=dashed,linewidth=0.1pt]{109.47}{180}{x cos -3 mul} % from mu=0
    \gsave
    \psline[xunit=0.0174](180,3)(180,0)(122.356,0)
    \fill[fillstyle=solid]
    \grestore
    }
  \pscustom[fillstyle=solid,fillcolor=EllipticColour]{ %% lower central elliptic region
    \psplot[xunit=0.0174]{76.34541523}{90}{x cos x cos 2 exp 3 mul x cos -4 mul add -3 add mul %
           x cos 3 mul -1 add 1 x  cos add 2 exp mul div -2 mul}
    \psline[xunit=0.0174](90,0)(122.356,0)
    \psplot[xunit=0.0174]{122.356}{131.81}{x cos x cos 2 exp 3 mul x cos -4 mul add -3 add mul %
           x cos 3 mul -1 add 1 x  cos add 2 exp mul div -2 mul}
    \psplot[xunit=0.0174]{130.2880417}{79.03921899}{x sin 2 exp 3 mul x cos 8 mul add 2 exp %
           x cos 3 mul -1 add 1 x  cos add 2 exp mul 8 mul div}
    }
  \pscustom[fillstyle=solid,fillcolor=EllipticColour]{ %% upper left elliptic sliver
    \psplot[xunit=0.0174]{60.714}{0}{x cos x cos 2 exp 3 mul x cos -4 mul add -3 add mul %
           x cos 3 mul -1 add 1 x  cos add 2 exp mul div -2 mul}
    \psplot[xunit=0.0174]{0}{59.01256515}{x sin 2 exp 3 mul x cos 8 mul add 2 exp %
           x cos 3 mul -1 add 1 x  cos add 2 exp mul 8 mul div}
    }
%% eigenvalue diagrams:
  \psline[linewidth=0.2pt](1.5,1.8)(1.5,3.2)
   \psline[linewidth=0.2pt](1.3,2.5)(1.7,2.5)
   \psdots(1.5,2)(1.5,3)
   \psdots(1.5,2.2)(1.5,2.8)
  \psline[linewidth=0.2pt](0.5,2)(0.5,3)
   \psline[linewidth=0.2pt](0.3,2.5)(0.7,2.5)
   \psdots[dotstyle=*](0.4,2.2)(0.4,2.8)(0.6,2.2)(0.6,2.8)
  \psline[linewidth=0.2pt](1.8,-2.5)(1.8,-3.5)
   \psline[linewidth=0.2pt](1.6,-3)(2,-3)
   \psdots[dotstyle=*](1.7,-3.4)(1.7,-2.6)(1.9,-3.4)(1.9,-2.6)
%% axes:
  \psaxes[linewidth=0.5pt,labels=y,ticks=y,arrowsize=0.15]{->}(0,0)(0,-4.1)(3.25,4.3)
  \rput(3.35,0){$\theta_0$}  \rput(-0.05,4.5){$\kappa$}
  \rput(1.5708,-0.3){\small $\frac{\pi}{2}$}
%  \rput(2.094,-0.3){\small $\frac{2\pi}{3}$}
  \rput(3.14159,-0.3){\small $\pi$}
%% curves xi = 0 and mu = 0:
  \psplot[xunit=0.0174,linestyle=dotted]{0}{133}{x cos 1 add -1 exp x cos mul -2 mul} \rput(2.5,4.4){$\xi=0$}
  \psplot[xunit=0.0174,linestyle=dashed,linewidth=0.6pt]{0}{182}{x cos -3 mul} \rput(3.35,3){$\mu=0$} \rput(0.5,-3){$\mu=0$}
  \psplot[xunit=0.0174]{76.34547974}{131.81}{x cos x cos 2 exp 3 mul x cos -4 mul add -3 add mul %
           x cos 3 mul -1 add 1 x  cos add 2 exp mul div -2 mul}
\rput(0.7853,-4.5){\it Northern hemisphere} %
 \rput(1.571,-4.5){$\vert$}
\rput(2.356,-4.5){\it Southern hemisphere}
\end{pspicture}
%%%%%%%%%

\begin{pspicture}(-1,-4.5)(11,4.3)  %% n = 2
  \psframe(-1,-4.3)(11,4.3)
  \rput(10.3,4){{\large n = 2}}
  \psline(9.7,4.3)(9.7,3.6)(11,3.6)
%% rescale x to almost radians:
  \psset{xunit=3}
%% stability regions:
  \pscustom[fillstyle=solid,fillcolor=LyapounovColour]{ %% lower left stability region
    \psplot[xunit=0.0174]{0}{120}{x cos 3 mul 2 add x cos mul -1 mul x cos 1 add 2 exp div}
    \psline(2.094,-3.5)
    \gsave
    \psline(0,-3.5)(0,-2)
    \fill[fillstyle=solid]
    \grestore
    }
  \pscustom[fillstyle=solid,fillcolor=LyapounovColour]{ %% right hand stability region - top half
    \psline(2.094,3.4)(2.094,1)
    \psplot[xunit=0.0174]{120}{131.8103148}{x cos 3 mul 2 add x cos mul -1 mul x cos 1 add 2 exp div}
    \gsave
    \psline[xunit=0.0174](131.8103148,0)(180,0)(180,3.4)%(2.094,3.4)
    \fill[fillstyle=solid]
    \grestore
    }
  \pscustom[fillstyle=solid,fillcolor=LyapounovColour]{ %% right hand stability region - lower half
    \psline[xunit=0.0174](180,0)(131.8103148,0)
    \psplot[xunit=0.0174]{131.8103148}{139.232}{x cos 3 mul 2 add x cos mul -1 mul x cos 1 add 2 exp div}
    \gsave
    \psline[xunit=0.0174](139.232,-3.5)(180,-3.5)(180,0)%(2.094,3.4)
    \fill[fillstyle=solid]
    \grestore
    }
%% eigenvalue diagrams:
  \psline[linewidth=0.2pt](1,1)(1,2)
   \psline[linewidth=0.2pt](0.8,1.5)(1.2,1.5)
   \psdots(0.9,1.5)(1.1,1.5)
  \psline[linewidth=0.2pt](2.25,-2.2)(2.25,-2.8)
   \psline[linewidth=0.2pt](2.15,-2.5)(2.35,-2.5)
   \psdots(2.2,-2.5)(2.3,-2.5)
  \psline[linewidth=0.2pt](1.5,-2.5)(1.5,-1.5)
   \psline[linewidth=0.2pt](1.333,-2)(1.666,-2)
   \psdots(1.5,-2.25)(1.5,-1.75)
  \psline[linewidth=0.2pt](2.8,0.5)(2.8,1.5)
   \psline[linewidth=0.2pt](2.65,1)(2.95,1)
   \psdots(2.8,0.7)(2.8,1.3)
  \psline[linewidth=0.2pt](1.8,1.5)(1.8,2.5)
   \psline[linewidth=0.2pt](1.6,2)(2,2)
   \psdot[dotstyle=*,dotsize=.2](1.8,2)
   \pscurve[linewidth=0.5pt]{->}(1.9,2.2)(2,2.3)(2.08,2.2)
   \pscurve[linewidth=0.5pt]{->}(1.7,1.6)(1.7,1)(1.85,0.7)
%% axes:
  \psaxes[linewidth=0.5pt,xunit=1.047,labels=y,ticks=y]{->}(0,0)(0,-3.5)(3.1,3.5)
  \rput(3.35,0){$\theta_0$}  \rput(-0.05,3.7){$\kappa$}
  \rput(1.5708,-0.3){\small $\frac{\pi}{2}$}
  \rput(2.094,-0.3){\small $\frac{2\pi}{3}$}
  \rput(3.14159,-0.3){\small $\pi$}
%% curves xi = 0 and mu = 0:
  \psplot[xunit=0.0174,linestyle=dotted]{0}{142}{x cos 1 add -1 exp x cos mul -1 mul} \rput(2.65,3.8){$\xi=0$}
  \psplot[xunit=0.0174,linestyle=dashed,linewidth=0.6pt]{0}{182}{x cos -2 mul} \rput(3.35,2){$\mu=0$}
%% hemispheres:
 \rput(0.7853,-3.8){\it Northern hemisphere} %
 \rput(1.571,-3.8){$\vert$}
 \rput(2.356,-3.8){\it Southern hemisphere}
\end{pspicture}
{\begin{pspicture}(-5,-0.2)(5,1)
\put(-4,0.55){Key:}
\psframe[fillstyle=solid,fillcolor=LyapounovColour](-3,0.5)(-2.5,1) 
\put(-2,0.52){Lyapounov stable}
\psframe[fillstyle=solid,fillcolor=EllipticColour](1.2,0.5)(1.7,1) \put(1.8,0.52){Elliptic}
\psframe(-1,-0.3)(-0.5,0.2)  \put(-0.2,-0.15){Linearly unstable}
\end{pspicture}}
%%%%%%%%%% 
%
\caption{Bifurcation diagrams for  $\CC_{3v}(R,p)$ and $\CC_{2v}(R,p)$
      relative equilibria (so a total of 4 and 3 vortices respectively);
      the polar vortex of strength $\kappa$ is at the North pole.
      Stability is modulo $\SO(2)$ about the polar axis, or modulo $\SO(3)$ when $\mu=0$ (see text). The circles represent the eigenvalues of
      the mode $\ell=1$. }
    \label{fig:CRp2/3}
  \end{center}
\end{figure}
%%%%%%%%%%%

%%%%%%%%%%%%
\paragraph{Discussion of the case $n=2$}
See Figure~\ref{fig:CRp2/3}.

\medskip

\noindent The $\CC_{2v}(R,p)$ relative equilibria are isosceles
triangles lying on a great circle, and for $\theta_0 = 2\pi/3$ the
triangle becomes equilateral. We again discuss the stability of
those with non-zero momentum. Indeed, any 3-vortex configuration with
zero momentum is a relative equilibrium since the reduced space is
just a point, and is consequently also Lyapounov stable relative to
$\SO(3)$  \cite{Pa92}.

For $n=2$ the symmetry adapted basis for $V_1$ is
$$\begin{array}{lll}
e_1 &=& \kappa\:\alpha^{(1)}_\theta + 2\cos\theta_0\,\delta x,\\
e_2 &=&\kappa\:\alpha^{(1)}_\phi + 2\sin\theta_0\,\delta y.
\end{array}
$$
with $\omega(e_1,e_2) = 2\kappa\sin\theta_0\left(2\cos\theta_0 +
\kappa \right)$ (which vanishes only when $\mu=0$). Following the
proof of Theorem~\ref{ring+pole/stab} we obtain
$$\dd^2H_\xi^{(1)} = 2\kappa\mu\begin{pmatrix}\frac{\kappa(1+\cos\theta_0)^2 + 3\cos^2\theta_0+
2\cos\theta_0}{\sin^2\theta_0}& 0\cr 0& -(1+2\cos\theta_0)\end{pmatrix}.
$$
Consequently, for $\mu\neq0$, the $\CC_{2v}(R,p)$ relative
equilibrium is Lyapounov stable if this matrix is definite, so if
$$
(1+2\cos\theta_0)[(1+\cos\theta_0)^2\kappa+\cos\theta_0(2+3\cos\theta_0)]<0.
$$
It is spectrally unstable if the inequality is reversed. See Figure
\ref{fig:CRp2/3}.
\begin{itemize}
\item
  There are two stable regions.
  For $\theta_0 < 2\pi/3$  the relative equilibria are stable provided the polar
  vorticity is less than a certain $\theta_0$ dependent critical
  value, while for $\theta_0 > 2\pi/3$ they are stable for all
  polar vorticities greater than a critical value. As
  $\theta_0 \to\pi$ this value goes to $-\infty$.

\item
  For $\theta_0=\pi/2$, where the $2$-ring is equatorial and the
  isosceles triangle is right-angled, they are stable if and only if
  $\kappa<0$. This is in agreement with \cite[Theorem~III.3]{PM98},
  with $\Gamma_1=\Gamma_2=1$, and $\Gamma_3=\kappa$.

\item \emph{The restricted three vortex problem} \quad The range
of stability when $\kappa=0$ does not coincide with the range of
stability for a single ring. Indeed the $\CC_{2v}(R)$ relative
equilibria are Lyapounov stable for all co-latitudes (see Theorem
\ref{polvani}) while $\CC_{2v}(R,p)$ is unstable for $\kappa=0$ and
$\theta_0\in (0,\pi/2)$. This means that if we place a \emph{passive
tracer} or \emph{ghost vortex} at the North pole and a ring of two
vortices in the Northern hemisphere, then the passive tracer will be
unstable.
\end{itemize}

\begin{remark}\label{rmk:CMS}
 The stability of $\CC_{nv}(R,p)$ relative equilibria
has also been studied in \cite{CMS03}. However our method differs
significantly from theirs in that we consider the definiteness of
the Hessian $\dd^2H_\xi|_\NN(x_e)$ on the $2n-2$ dimensional
symplectic slice, while in \cite{CMS03} the authors determine
conditions for the Hessian to be definite on the whole $2n+2$
dimensional tangent space. The result is that we prove the relative
equilibria to be Lyapounov stable in a larger region of the
parameter space. Notice in particular that for $n \leq 6$ our
results say that a positive vorticity $n$-ring near the South pole
is Lyapounov stable if the North pole has either negative or
sufficiently positive vorticity. However in \cite{CMS03} only the
case of negative North polar vorticity is shown to be Lyapounov
stable. In this paper we also give criteria for when the relative
equilibria are {\it unstable} by considering the eigenvalues of the
linearization $L_\NN$.
\end{remark}

%%%%%%%%%%%%%%%%%%%%%%%%%%
\section{Two aligned rings: $\CC_{nv}(2R)$}
\label{sec:2R}

In this section we consider relative equilibria $x_e$ of symmetry
type $\CC_{nv}(2R)$, that is configurations formed of two `aligned'
rings of $n$ vortices each as illustrated in Figure
\ref{fig:2rings}. We can assume without loss of generality that the
vorticities of the vortices in the first and second ring are  $1$
and  $\kappa$, respectively, and we denote their co-latitudes by
$\theta_1$ and $\theta_2$. We can also assume that the ring of
vorticity $1$ and co-latitude $\theta_1$ lies in the Northern
hemisphere, $\theta_1\in (0,\,\pi/2]$. The first question to answer
is, \emph{for which values of the parameters
$(\theta_1,\theta_2,\kappa)$ is the configuration $\CC_{nv}(2R)$ a
relative equilibrium?}  It was shown in \cite{LMR01} (p.~126) that
for given $\kappa>0$ and each $\mu$ with $|\mu| < n|1+\kappa|$ there
is at least one solution for $(\theta_1,\theta_2)$ with
$n\cos\theta_1+n\kappa\cos\theta_2=\mu$ and with $\theta_1<\theta_2$
and at least one with $\theta_2<\theta_1$. We now make this more
precise.

The isotropy subgroup $G_{x_e}$ is the dihedral group
$\CC_{nv}$. The fixed point set $\Fix(G_{x_e})$ consists
of all pairs of aligned rings, with one vortex from each ring on a
given meridian, so can be parametrized by $x:=\cos\theta_1$ and
$y:=\cos\theta_2$. Denote by $\tilde F$ the restriction of a
function $F$ to $\Fix(G_{x_e})$. The Hamiltonian can be split in
such a way that
$$
H\ =\ H_{11}+\kappa H_{12}+\kappa^2 H_{22}
$$
where $H_{11},H_{12},H_{22}$ do not depend on $\kappa$, $\tilde
H_{11}$ does not depend on $y$ and $\tilde H_{22}$ does not depend
on $x$ ($H_{11}$ governs the interactions within the first ring, 
$H_{12}$ the interactions between the rings etc).  The following 
proposition shows that for almost every pair
$(\theta_1,\theta_2)$ there exists a unique $\kappa$ such that the
$\CC_{nv}(2R)$ configuration with parameters
$(\theta_1,\theta_2,\kappa)$ is a relative equilibrium.

%%%%%%%%%%%%%%%%%
\begin{proposition} \label{prop:exist kappa 2R}
Let $x_{e}$ be a $\CC_{nv}(2R)$ configuration with parameters
$(\theta_1,\theta_2,\kappa)$.
\begin{enumerate}
\item There exists a unique $\kappa\in\R^*$ such that $x_e$ is a
relative equilibrium if and only if both the following conditions
hold:
$$
\left(\frac{\partial \tilde H_{12}}{\partial y}-\frac{\partial
    \tilde H_{11}}{\partial x}\right) (\cos\theta_1,\cos\theta_2)\neq
0,\ \ \left(\frac{\partial \tilde H_{22}}{\partial y}-\frac{\partial
    \tilde H_{12}}{\partial x}\right) (\cos\theta_1,\cos\theta_2)\neq
0.
$$

\item The configuration $x_e$ is a relative equilibrium for {\it
all} $\kappa\in\R^*$ in the degenerate case when both the following
conditions hold:
$$
\left(\frac{\partial \tilde H_{12}}{\partial y}-\frac{\partial
    \tilde H_{11}}{\partial x}\right) (\cos\theta_1,\cos\theta_2) = 0,
\ \ \left(\frac{\partial \tilde H_{22}}{\partial
    y}-\frac{\partial \tilde H_{12}}{\partial x}\right)
(\cos\theta_1,\cos\theta_2) = 0.
$$
\end{enumerate}

In both cases the angular velocity $\xi$ of $x_e$ satisfies
$$
\xi\ =\ \frac{1}{n}\left(\frac{\partial \tilde H_{11}}{\partial
    x}(x_e)+\kappa\frac{\partial \tilde H_{12}}{\partial
    x}(x_e)\right).
$$
\end{proposition}
%%%%%%%%%%%%%%%%

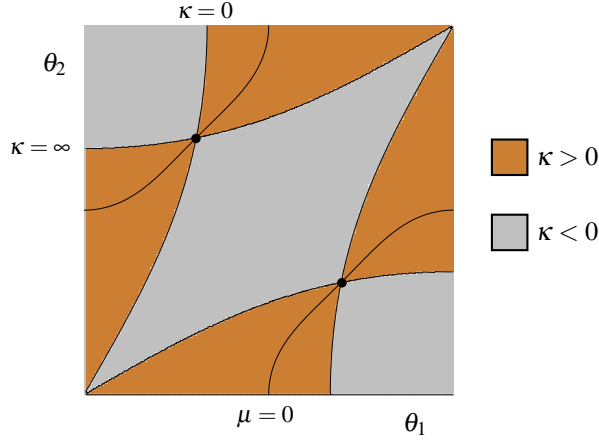
\begin{figure}[t]  % Sign of \kappa for 2 aligned rings
 \begin{center}
\begin{pspicture}(0,-0.5)(6,5) 
\put(0,0){\includegraphics[width=2in,height=2in]{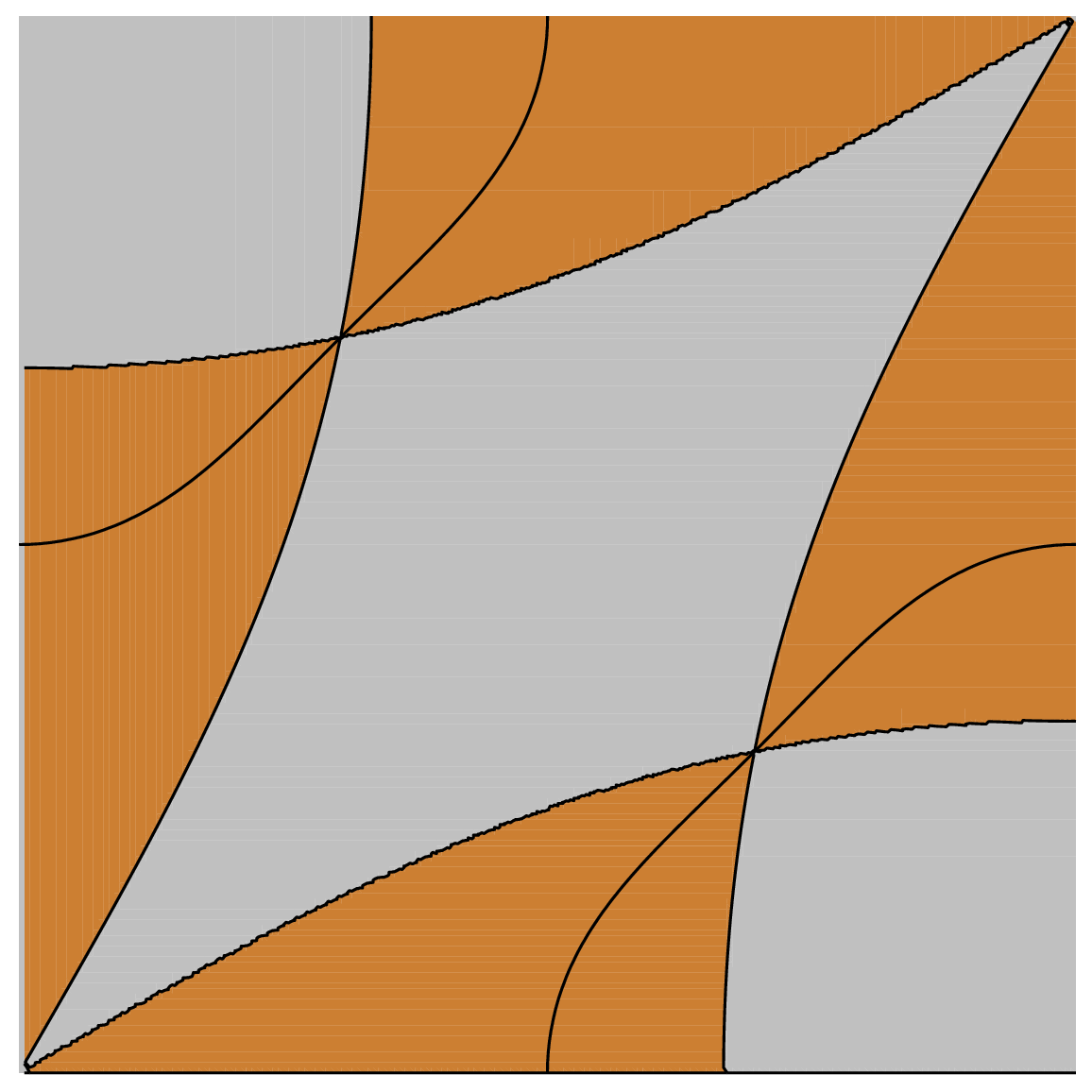}}
\rput(4.5,-0.3){$\theta_1$}\rput(-0.3,4.5){$\theta_2$}
\rput(2.5,-0.2){\small$\mu=0$}
\rput(1.7,5.2){\small$\kappa=0$} 
\put(-0.9,3.3){\small$\kappa=\infty$}
\psdot(1.58,3.5)
\psdot(3.52,1.58)
\psframe[fillstyle=solid,fillcolor=gold](5.5,3)(6,3.5) \put(6.1,3.15){$\kappa>0$}
\psframe[fillstyle=solid,fillcolor=midgray](5.5,2)(6,2.5) \put(6.1,2.15){$\kappa<0$}
\end{pspicture}
    \caption{Sign of $\kappa$ for $\CC_{nv}(2R)$ relative equilibria.  
    The degenerate case occurs where
    the curves $\kappa=0$ and $\kappa=\infty$ intersect.
    The figure is for $n=4$, but is similar for other values of
    $n$. The only region of stability  (see Fig.~\ref{fig:Cnv2R}) lies at the bottom right
    hand and top left-hand corners, corresponding to the rings lying far apart in
    opposite hemispheres, and is contained in the region $\kappa<0$.
    Along both diagonals (where $\theta_1=\theta_2$ and
    $\theta_1+\theta_2=\pi$) one has $\kappa=-1$.
    The two black curves are where $\mu=0$. Note that reflecting along 
    either diagonal exchanges the rings, so corresponds to changing $\kappa$ to $\kappa^{-1}$.
    } \label{fig:signkappa-aligned}
\end{center}
 \end{figure}

The sign of $\kappa$ as a function of $\theta_1,\theta_2$ is shown in Figure~\ref{fig:signkappa-aligned}.

%%%%%%%%%%%%%%%%%
\begin{proof}
  Since $H-\xi\Phi$ is a $G_{x_e}$-invariant function
  (see Section \ref{intro}) the Principle of Symmetric Criticality
  \cite{P79} implies that $x_e$ is a relative equilibrium if and only if it is a critical
  point of $\tilde H-\xi\tilde\Phi$.
  It follows from $\tilde\Phi=n(x+\kappa y)$ that $\dd(\tilde
  H-\xi\tilde\Phi)(x_e)=0$ is equivalent to the pair of equations:
  $$
  \kappa\left(\frac{\partial \tilde H_{22}}{\partial
      y}(x_e)-\frac{\partial \tilde H_{12}}{\partial x}(x_e)\right) +
  \frac{\partial \tilde H_{12}}{\partial y}(x_e)-\frac{\partial \tilde
    H_{11}}{\partial x}(x_e)\ =\ 0
  $$
  $$
  \xi=\frac{1}{n}\left(\frac{\partial \tilde
      H_{11}}{\partial x}(x_e)+\kappa\frac{\partial \tilde
      H_{12}}{\partial x}(x_e)\right).
  $$
  The proposition follows easily from these.
\end{proof}
%%%%%%%%%%%%%%%%

%%%%%%%%%%%
%%%%%%%%%%%

For example, in the case $n=4$ the degenerate case occurs when the
two rings form the vertices of a cube. Hence for any values of the
vorticities of the two rings the `cube configuration' is a
relative equilibrium. However among this family of relative
equilibria only one is an equilibrium, namely the one for which the
two rings have the same vorticities, $\kappa=1$, which corresponds
to the $\mathbb{O}_h(f)$ equilibrium of \cite{LMR01}, a cube formed of
identical vortices, and which is known to be unstable \cite{Ku04}. In Figure~\ref{fig:signkappa-aligned}, the cubic configurations are marked by two dots. For $n=2$, the degenerate case occurs when the vortices form a square lying on a great circle.

For $\theta_2=\pi-\theta_1$, the configuration has an extra symmetry
and its symmetry type is $D_{nh}(2R)$. Such a configuration is a
relative equilibrium if the two rings have opposite vorticities
($\kappa=-1$). The existence and stability of such relative
equilibria were studied in \cite{LP02}.

With the help of the discussion of Section \ref{tworings}, one can
perform a $G_{x_e}$-invariant isotypic decomposition and find that
the symmetry adapted basis for the symplectic slice at a
$\CC_{nv}(2R)$ relative equilibrium with $n\geq3$ and $\mu\neq 0$ is
$$
  \left( e_1,e_2,f_1,f_2,f_3,f_4,f_5,f_6,
  B_2,B_3,\dots,B_{[n/2]}\right)
$$
where
$$
\begin{array}{lll}
e_1&=&\alpha^{(0)}_{0,\phi}-\alpha^{(0)}_{1,\phi}\\
e_2&=&\kappa\sin\theta_2\;\alpha^{(0)}_{0,\theta}-\sin\theta_1\;\alpha^{(0)}_{1,\theta}\\[12pt]
 \end{array}
 $$
% and
 $$
 \begin{array}{lllclll}
 f_1 &=& \sin\theta_1\; \alpha^{(1)}_{0,\theta} + \cos\theta_1\;\beta^{(1)}_{0,\phi}&\qquad&%\\[6pt]
 f_2 &=& \sin\theta_2\; \alpha^{(1)}_{1,\theta} + \cos\theta_2\;\beta^{(1)}_{1,\phi}\\[6pt]
 f_3 &=& \kappa\sin\theta_2 \; \beta^{(1)}_{0,\phi} - \sin\theta_1\;\beta^{(1)}_{1,\phi}&&%\\[6pt]
 f_4 &=&  \cos\theta_1\;\alpha^{(1)}_{0,\phi} - \sin\theta_1\;\beta^{(1)}_{0,\theta},\\[6pt]
 f_5 &=&  \cos\theta_2\;\alpha^{(1)}_{1,\phi} - \sin\theta_2\;\beta^{(1)}_{1,\theta}&&%\\[6pt]
 f_6 &=&  \kappa\sin\theta_2\;\alpha^{(1)}_{0,\phi} - \sin\theta_1\;\alpha^{(1)}_{1,\phi}
\end{array}
$$
and,
$$
\begin{array}{rcl}
B_\ell & = &
\left\lbrace\alpha^{(\ell)}_{0,\theta},\:\alpha^{(\ell)}_{1,\theta},\:
  \alpha^{(\ell)}_{0,\phi},\:\alpha^{(\ell)}_{1,\phi},\:
  \beta^{(\ell)}_{0,\theta},\:\beta^{(\ell)}_{1,\theta},\:
  \beta^{(\ell)}_{0,\phi},\:\beta^{(\ell)}_{1,\phi}
  \right\rbrace\quad \mbox{for
  $2\leq\ell<n/2$}\\[6pt]
B_{n/2} & = &
\left\lbrace\alpha^{(n/2)}_{0,\theta},\:\alpha^{(n/2)}_{1,\theta},\:
  \alpha^{(n/2)}_{0,\phi},\:\alpha^{(n/2)}_{1,\phi}\right\rbrace\quad \mbox{for
  even $n$}.\\
\end{array}
$$

The adapted basis for $n=2$ is simply $(e_1,e_2;f_6,f_7)$, where
$$f_7 = \kappa \cos\theta_2\;\alpha^{(1)}_{0,\theta} - \cos\theta_1\;\alpha^{(1)}_{1,\theta}\,.
$$

\medskip

\noindent\emph{Remark } Almost all $\CC_{nv}(2R)$ relative equilibria have a non-zero
momentum.  Indeed $\mu=0$ if and only if $x+\kappa y=0$, and from
the expression of $\kappa$ one can show that this last equation
defines an algebraic curve in variables $(x,y)\in [0,1)\times
(-1,1)\simeq\Fix \CC_{nv}$ (depicted in Figure~\ref{fig:signkappa-aligned}). Numerics suggest that equilibria ($\xi=0$) occur along curves that are extremely close to these momentum zero curves, and indeed would be indistinguishable in the diagram.

\bigskip

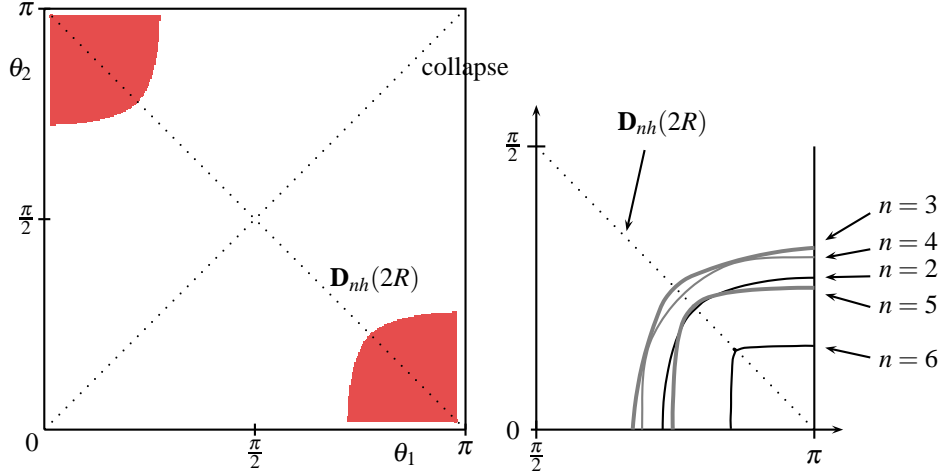
\begin{figure}[t]
\begin{center}
\psset{unit=0.8}
%\fbox
{{\begin{pspicture}(0,-1.5)(7,7.5) % n  = 2 
\put(-0.1,-0.05){\includegraphics[width=2.25in,height=2.25in]{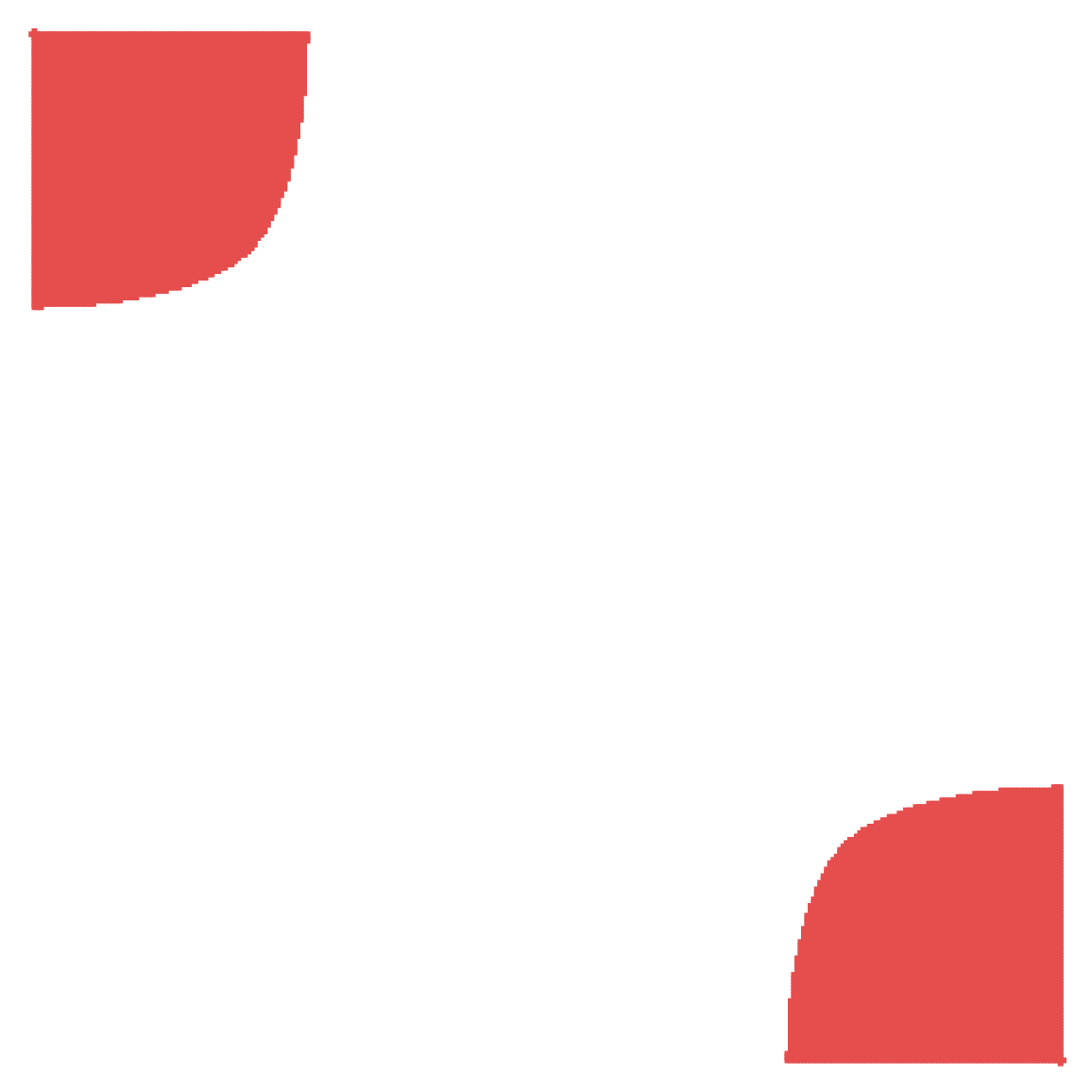}}
\rput(-0.2,-0.2){0}
\psline(0,0)(0,7) \rput(-0.3,7){$\pi$} \rput(-0.3,3.5){$\frac{\pi}2$} 
   \psline(-0.1,3.5)(0.1,3.5)  \psline(-0.1,7)(0.1,7) \rput(-0.4,6){$\theta_2$}
\psline(0,0)(7,0) \rput(7,-0.3){$\pi$} \rput(3.5,-0.4){$\frac{\pi}2$}
   \psline(3.5,-0.1)(3.5,0.1)  \psline(7,-0.1)(7,0.1) \rput(6,-0.4){$\theta_1$}
\psline(7,0)(7,7)(0,7)
\psline[linestyle=dotted](0,0)(7,7) \rput(7,6){collapse}
\psline[linestyle=dotted](0,7)(7,0) \rput(5.5,2.5){$\DD_{nh}(2R)$}
\end{pspicture}}\qquad
\psset{unit=3}
\begin{pspicture}(1.5,-0.5)(3.7,2)
 %axes:
  \psline{->}(1.57,0)(3.3,0) %\rput(3.4,-0.1){$\theta_1$}
  \psline{->}(1.6,0)(1.6,1.8) %\rput(-0.1,1.9){$\theta_2$}
  %% tickmarks
  \rput(1.47,0){0}
  \psline(1.55,1.57)(1.65,1.57) \rput(1.47,1.57){$\frac\pi2$}
  \psline(1.6,-0.05)(1.6,0.05) \rput(1.6,-0.15){$\frac\pi2$}
  \psline(3.14,-0.05)(3.14,1.57) \rput(3.14,-0.15){$\pi$}
  \psline[linestyle=dotted](3.141,0)(1.6,1.57)
  \rput(2.3,1.7){$\mathbf{D}_{nh}(2R)$}
  \psline{->}(2.2,1.55)(2.1,1.1)
  %% n=2 stability
    \pscurve(2.300, 0.00084)(2.302, .1021)(2.310, .2043)%
      (2.326, .3080)(2.349, .4136)(2.384, .5191)(2.436, .6205)%
      (2.515, .7006)(2.521, .7061)(2.622, .7571)(2.728, .7926)%
      (2.834, .8158)(2.937, .8313)(3.039, .8398)(3.141, .8420)
    \psline{<-}(3.2, .8446396945)(3.45,0.85)  \put(3.5,0.85){$n=2$}
  %% n=3 stability
    \pscurve[linecolor=gray,linewidth=1.5pt](2.134, 0.001)%
      (2.145, .1212)(2.163, .2406)(2.187, .3603)(2.220, .4812)%
      (2.261, .6039)(2.315, .7261)(2.415, .8263)(2.538, .8809)%
      (2.660, .9220)(2.781, .9543)(2.901, .9789)(3.020, .9967)(3.141, 1.008)
    \psline{<-}(3.2, 1.05)(3.45,1.2)  \put(3.5,1.2){$n=3$}
  %% n=4 stability
    \pscurve[linecolor=gray](2.186, 0.0009)(2.187, .1161)%
      (2.188, .2344)(2.199, .3560)(2.228, .4765)(2.284, .5881)%
      (2.357, .6892)(2.452, .7842)(2.554, .8578)(2.665, .9131)%
      (2.786, .9430)(2.907, .9536)(3.025, .9550)(3.141, .9560)
    \psline{<-}(3.2, .9538064985)(3.45,1.)  \put(3.5,1.){$n=4$}
  %% n=5 stability
    \pscurve[linecolor=gray,linewidth=1.5pt](2.356, 0.00078)%
       (2.355, 0.095)(2.359, .1924)(2.363, .2939)(2.372, .4016)%
       (2.389, .5157)(2.427, .6284)(2.513, .7151)(2.626, .7522)%
       (2.740, .7695)(2.848, .7784)(2.949, .7827)(3.046, .7862)(3.141, .7860)
    \psline{<-}(3.2, .75)(3.45,0.68)  \put(3.5,0.65){$n=5$}
  %% n=6 stability
    \pscurve(2.678, 0.0046)(2.676, 0.0514)(2.678, 0.0988)%
       (2.679, .1483)(2.679, .2009)(2.681, .2578)(2.684, .3197)%
       (2.688, .3898)(2.697, .4489)(2.693, .4448)(2.752, .4535)%
       (2.822, .4574)(2.884, .4608)(2.941, .4622)(2.993, .4628)%
       (3.043, .4636)(3.090, .4652)(3.137, .4640)
     \psline{<-}(3.2, .4544)(3.45,0.4)  \put(3.5,0.35){$n=6$}
\end{pspicture}}
 \caption{Stability results for $\CC_{nv}(2R)$ relative equilibria; that is for two aligned rings. 
       The coloured regions in the left-hand diagram represents the configurations that are Lyapounov stable. The figure is for $n=2$, but is very similar for $n\leq6$. The right-hand figure illustrates the different sizes of the stable regime for the different values of $n$. If $n$ is odd there is a narrow strip (too narrow to discern on this diagram) between the Lyapounov stable and the linearly unstable configurations where the relative equilibria are elliptic.
       For $n\geq 7$, it seems that all relative equilibria  are unstable.}
      \label{fig:Cnv2R}
\end{center}
\end{figure}

With respect to this basis $\dd^2H_\xi|_\NN(x_e)$ block diagonalizes
into: two $1\times1$ blocks for $\ell=0$,  two $3\times3$ blocks for
$\ell=1$,  two $4 \times 4$ blocks for each of $\ell = 2 \ldots
[(n-1)/2]$), together with two $2 \times 2$ blocks for $\ell = n/2$
when $n$ is even. The linearization $L_\NN$ block diagonalizes into
half as many blocks of twice the size. In order to calculate the
stability of the relative equilibria, we ran a \textsc{Maple}
program to compute numerically the eigenvalues of each of the
blocks of $\dd^2H_\xi|_\NN(x_e)$ and $L_\NN$. The results are
summarized for $n = 2 \ldots 6$ in Figure~\ref{fig:Cnv2R}. Figure
\ref{fig:signkappa-aligned} shows how the sign of $\kappa$ varies for
relative equilibria with different values of $\theta_1$ and
$\theta_2$. A selection of the \textsc{Maple} code is available for download from \cite{Mo-web}; further diagrams are also available from the same site.

\paragraph{Discussion of  results} Here we outline conclusions from a series of numerical calculations, using \textsc{Maple}.  These involved calculating eigenvalues of the Hessian matrix and the linearization using the Fourier bases described above, varying the co-latitudes $\theta_1$ and $\theta_2$ across the range 0 to $\pi$ in steps of $10^{-2}$ (and occasionally smaller steps to investigate specific bifurcations).

\begin{itemize}
\item These numerical calculations suggest strongly that the relative
equilibria $\CC_{nv}(2R)$ are never stable if the two rings lie in
the same hemisphere (Figure~\ref{fig:Cnv2R}) or have the same sign
vorticity (Figure~\ref{fig:signkappa-aligned}).

\item The stable configurations are for one ring close to the North pole
and the other ring close to the South pole, and always with
vorticity of opposite signs. For $n=4$, the furthest from the poles
both rings can be is about $40^\circ$ of latitude. This is in
agreement with Theorem~4.8 of \cite{LP02}. 

\item For $n>2$, as $n$ increases the region of
stability decreases in size. Numerical experiments with $n \geq 7$
suggest that in these cases the relative equilibria are never
stable.

\item The stability boundaries for $2\leq n\leq 6$ approach
configurations where one ring is at a pole and the other lies at a
particular co-latitude $\theta_z$ in the opposite hemisphere. For
$n=3,4$, $\theta_z\approx 0.95 \textrm{rad} \approx 55^\circ$, for
$n=2,5$, $\theta_z\approx 0.8 \textrm{rad} \approx 46^\circ$, and
for $n=6$, $\theta_z\approx 0.45 \textrm{rad} \approx 25^\circ$.

\item For $n = 2$, $4$ and $6$ stability is first lost by a
pair of imaginary eigenvalues of the $\ell = n/2$ block of $L_\NN$
passing through $0$ and becoming real (the `splitting' case
described in Section \ref{sec:symmetry+bifurcations}). The difference between even and odd $n$ (see next point) is that the type of representation on the $\ell=n/2$ mode is `of real type'.

\item For $n = 3$ and $5$ a pair of imaginary eigenvalues of the $\ell
= (n-1)/2$ block passes through $0$ but remains on the imaginary
axis, so the stability changes from Lyapounov to elliptic (the
`passing' case described in Section \ref{sec:symmetry+bifurcations}). This
imaginary pair then collides with another pair, and all move off the
imaginary axis to form a complex quadruplet and create instability (a
Hamiltonian Hopf bifurcation). However, the elliptic regions are
very narrow: for $n=3$ with $\theta_1=0.1$ the elliptic region is
contained in $|\theta_2-2.1438318| < 10^{-7}$. For $\theta_2 =
\pi-\theta_1$ (the $\mathbf{D}_{nh}(2R)$ configurations)
the elliptic range is at its widest, but is still only approximately
$|\theta_2-0.777| < 6\times 10^{-3}$, which is too small to be seen
in Figure~\ref{fig:Cnv2R}. For $n=5$ the elliptic region is even
narrower.

\item As $(\theta_1,\,\theta_2)$ approaches the diagonal (that is, the
rings approach one another) so $\kappa\to -1$, and the configuration approaches one of $n$ dipoles (pairs with equal and opposite vorticities).

\item When $\kappa=0$ (or $\kappa^{-1}=0$) the system is not Hamiltonian, as the symplectic form is degenerate, and since we are using Hamiltonian methods further work would be needed to find the stability at these points.
\end{itemize}

%%%%%%%%%%%%%%%%%%%%%%%%%%%%%%%%%%%%%%%%%%%%%%%%%%%%%%%%%%%%%%%%%%%%%%%%
\section{Two staggered rings: $\CC_{nv}(R,R')$}
\label{sec:RR'}

In this section we consider relative equilibria formed of two rings
of $n$ vortices each of strengths $1$ and $\kappa$ and co-latitude
$\theta_1$ and $\theta_2$ respectively. They differ from those of
the previous section in that the rings here are `staggered', that
is they rotated relative to each other with an offset of $\pi/n$.
Their symmetry type is $\CC_{nv}(R,R')$. As in the previous section
we can assume without loss of generality that the ring of vorticity
$1$ and co-latitude $\theta_1$ lies in the Northern hemisphere. The difference between the stabilities of aligned and staggered rings is striking: compare Figures~\ref{fig:Cnv2R} and \ref{fig:CnvRR'}.

Proposition~\ref{prop:exist kappa 2R} also holds for
$\CC_{nv}(R,R')$ configurations: for almost every pair
$(\theta_1,\theta_2)$ there exists a unique $\kappa$ such that the
corresponding $\CC_{nv}(R,R')$ configuration is a relative
equilibrium, and $\kappa,\xi$ are given by the same expressions in terms of the derivatives of $\tilde H_{ij}$. The only difference from Section \ref{sec:2R} is the expression for $\tilde H_{12}$ in terms of $\theta_1$ and $\theta_2$ (or $x$ and $y$).

When $\theta_2 = \theta_1$ the
configuration forms a single ring with $2n$ vortices with
$\kappa=1$: all the vortices have the same vorticity. These are the
relative equilibria of type $\CC_{2nv}(R)$ studied in Section
\ref{singlering}.  For $\theta_2 = \pi-\theta_1$, the configuration
has an extra symmetry and its symmetry type is $D_{nd}(R,R')$. In
this case $\kappa=-1$, the two rings have opposite vorticities. The
existence and stability of such relative equilibria were studied in
\cite{LP02}.

There exist also degenerate cases in the sense of Proposition~\ref{prop:exist kappa 2R}. For $n=2$ these are, a square on a great circle which is a relative equilibrium whenever the opposite vortices have the same vorticity and the tetrahedral configuration, which is shown in \cite{PM98} to be a relative equilibrium for any values of the four vorticities. These are discussed further below.

With the help of the discussion of Section \ref{tworings}, one finds
that a symmetry adapted basis for the symplectic slice at a
$\CC_{nv}(R,R')$ relative equilibrium with $n\geq 3$ and $\mu\neq 0$
is given by:
$$
  \left( e_1,e_2,e_3,e_4,e_5,e_6,e_7,e_8,
  B_2,B_3,\dots,B_{[n/2]}\right)
$$
for $n$ odd, while for $n$ even one is given by: %
$$
\begin{array}{r}
\left( e_1,e_2,e_3,e_4,e_5,e_6,e_7,e_8,\left\lbrace B_\ell \mid
    2\leq \ell\leq n/2-1 \right\rbrace,\right.\quad \quad \quad \\
  \left.\alpha^{(n/2)}_{0,\theta}-\alpha^{(n/2)}_{1,\theta}, 
  \alpha^{(n/2)}_{0,\phi}-\alpha^{(n/2)}_{1,\phi},
  \beta^{(n/2)}_{1,\theta},\beta^{(n/2)}_{1,\phi} \right),
\end{array}
$$
where the expressions of $e_1,\dots,e_8$ and $B_\ell$ remain as in
the previous section. The corresponding symmetry adapted basis for
$n=2$ is simply $(e_1,e_2,e_3,e_6)$. As in the previous section, it
can readily be seen that almost all $\CC_{nv}(R,R')$ relative
equilibria have non-zero momenta.

As for the aligned rings, we ran a \textsc{Maple} program to
determine the stability of the relative equilibria. The results
are summarized in Figure~\ref{fig:CnvRR'} for $n$ from 2 to 6.

%%%%%%%%%%%
\begin{figure}[t] 
 \begin{center}
\begin{pspicture}(0,-1)(5,5.5) 
\put(0,0){\includegraphics[width=2in,height=2in]{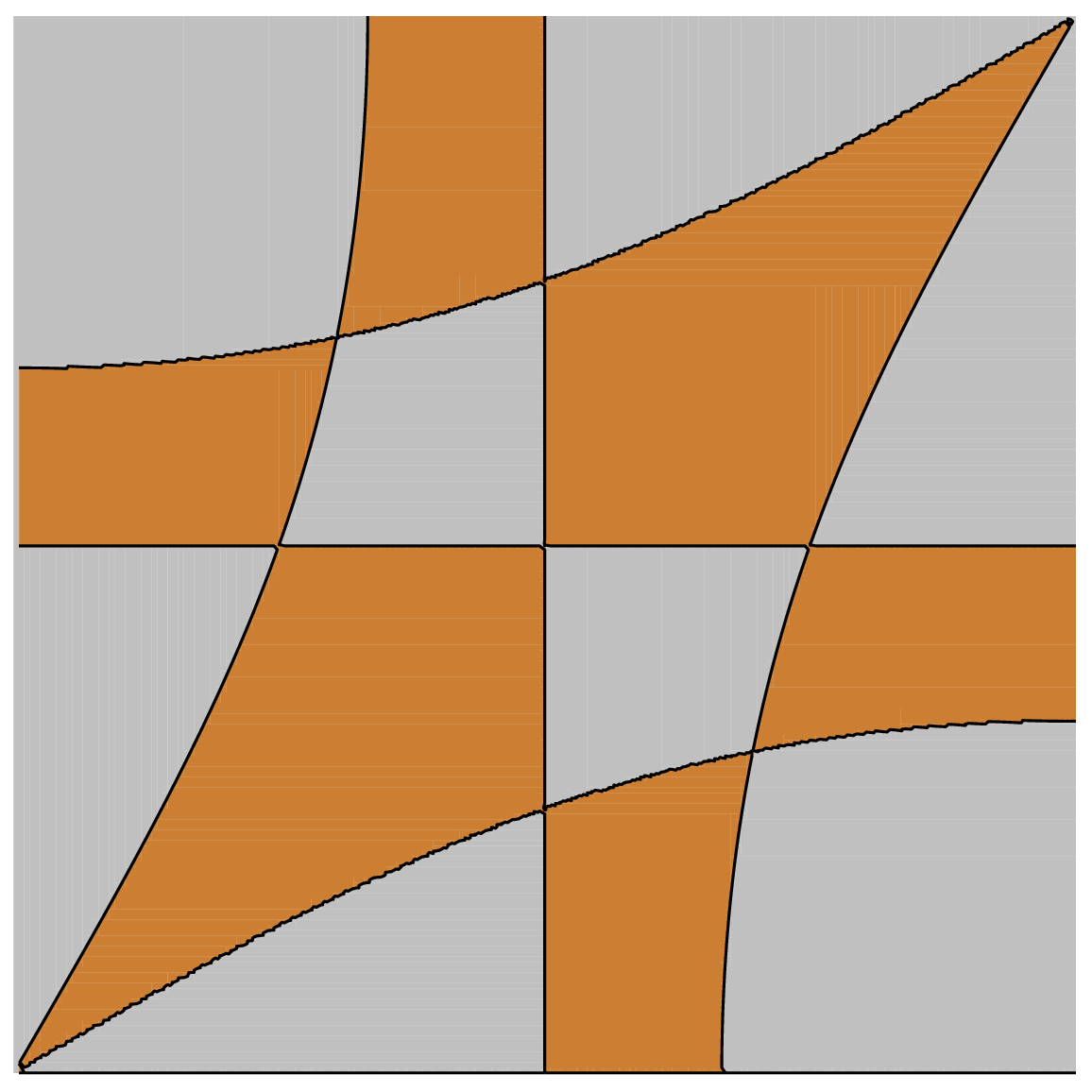}}
\rput(2.5,-0.5){$n=2$}
\rput(1.7,5.3){$0$} \rput(2.5,5.3){0}
\put(-0.3,2.5){$\infty$}\put(-0.3,3.3){$\infty$}
\psframe[fillstyle=solid,fillcolor=gold](4.5,-0.7)(5,-0.2) \rput(4.75,-0.9){$\kappa>0$}
\end{pspicture}
\qquad
\begin{pspicture}(0,-1)(5,5.5) 
\put(0,0){\includegraphics[width=2in,height=2in]{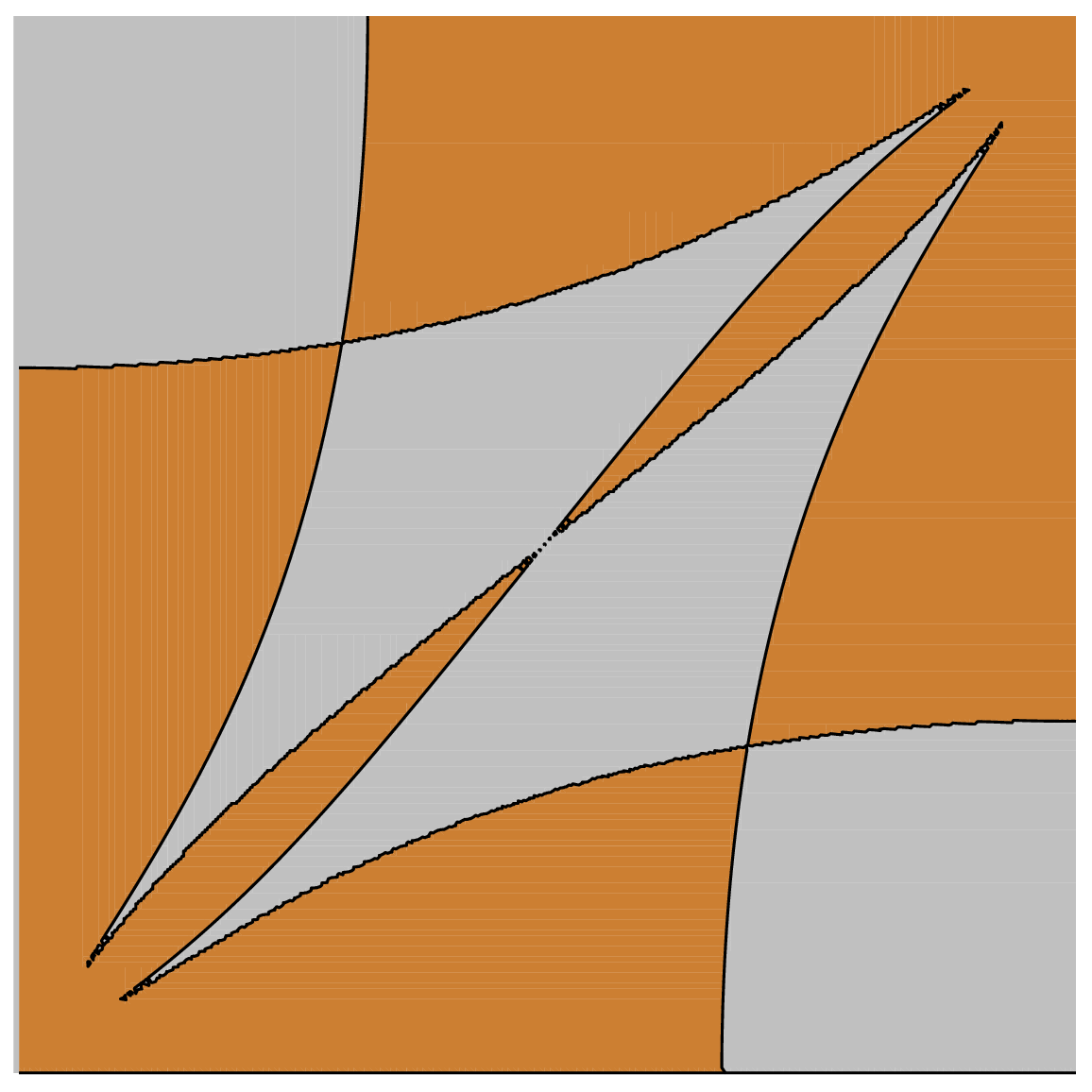}}
\rput(2.5,-0.5){$n=4$}
\rput(1.7,5.3){$0$}
\rput(5.3,1.7){$\infty$}
\psframe[fillstyle=solid,fillcolor=midgray](0,-0.7)(0.5,-0.2) \rput(0.25,-0.9){$\kappa<0$}
\end{pspicture}
    \caption{Sign of $\kappa$ for $\CC_{nv}(R,R')$ relative equilibria for $n=2$ and 4. The diagrams for $n\geq3$ are all similar, with the central strip where $\kappa>0$ getting thinner as $n$ increases. See the discussion of results below for more details.
    } \label{fig:signkappa-staggered}
\end{center}
 \end{figure}
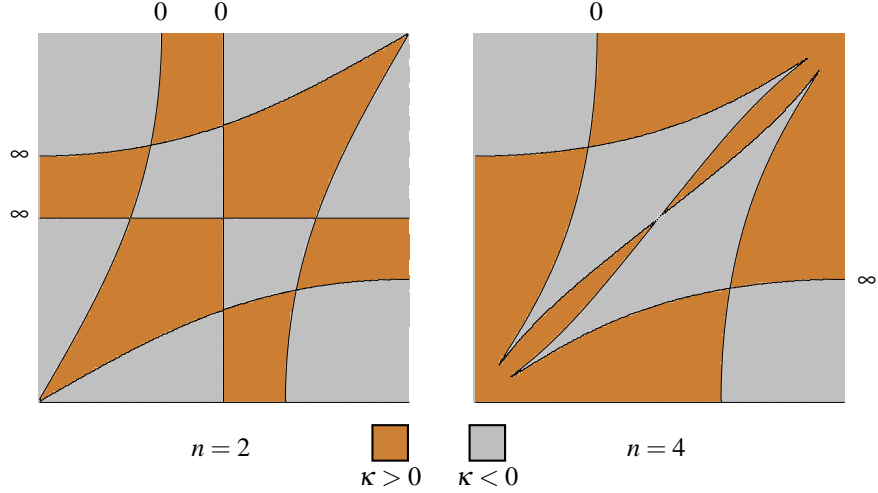
%%%%%%%%%%%

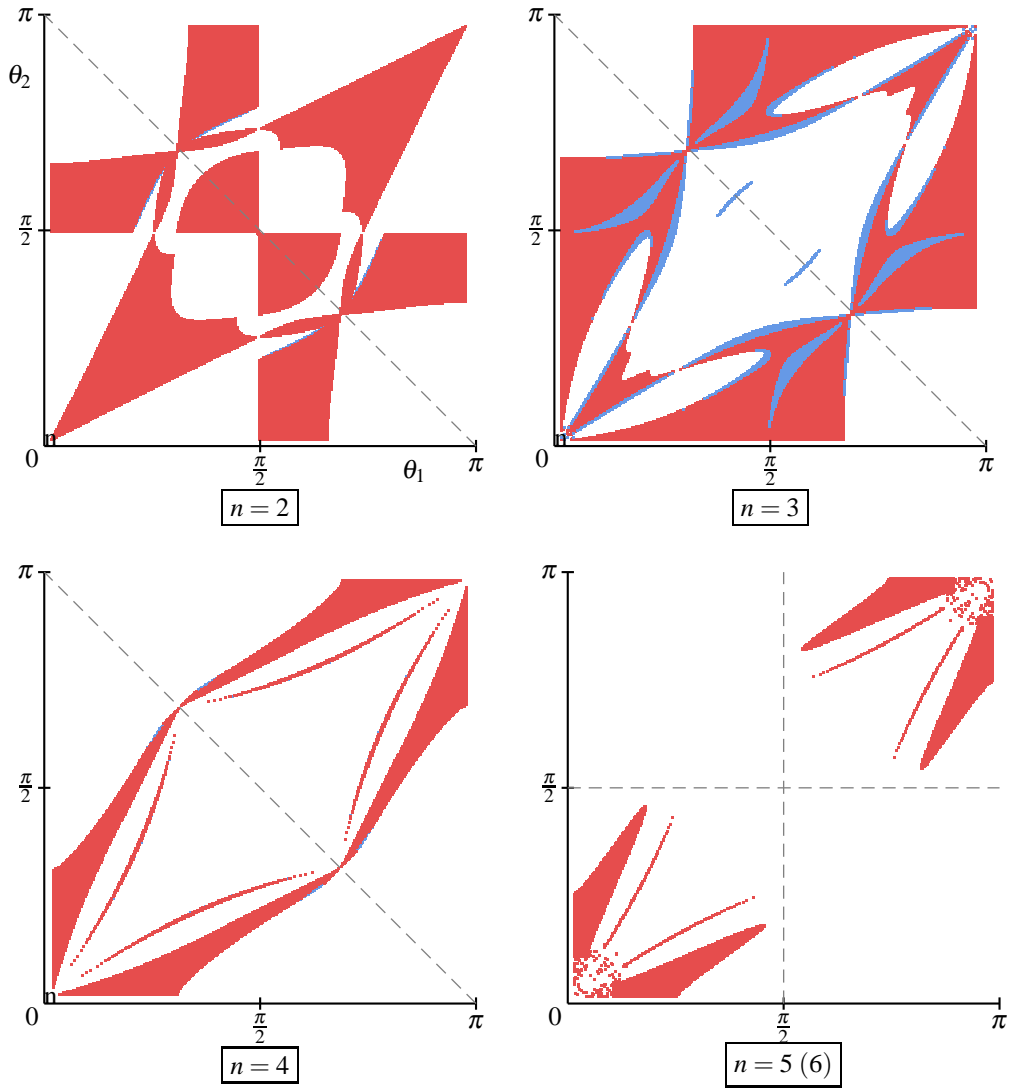
\begin{figure}[tp]
\begin{center}\psset{unit=0.82}
{\begin{pspicture}(0,-1.5)(7,7.5) % n  = 2 
\put(-0.1,-0.1){\includegraphics[width=2.3in]{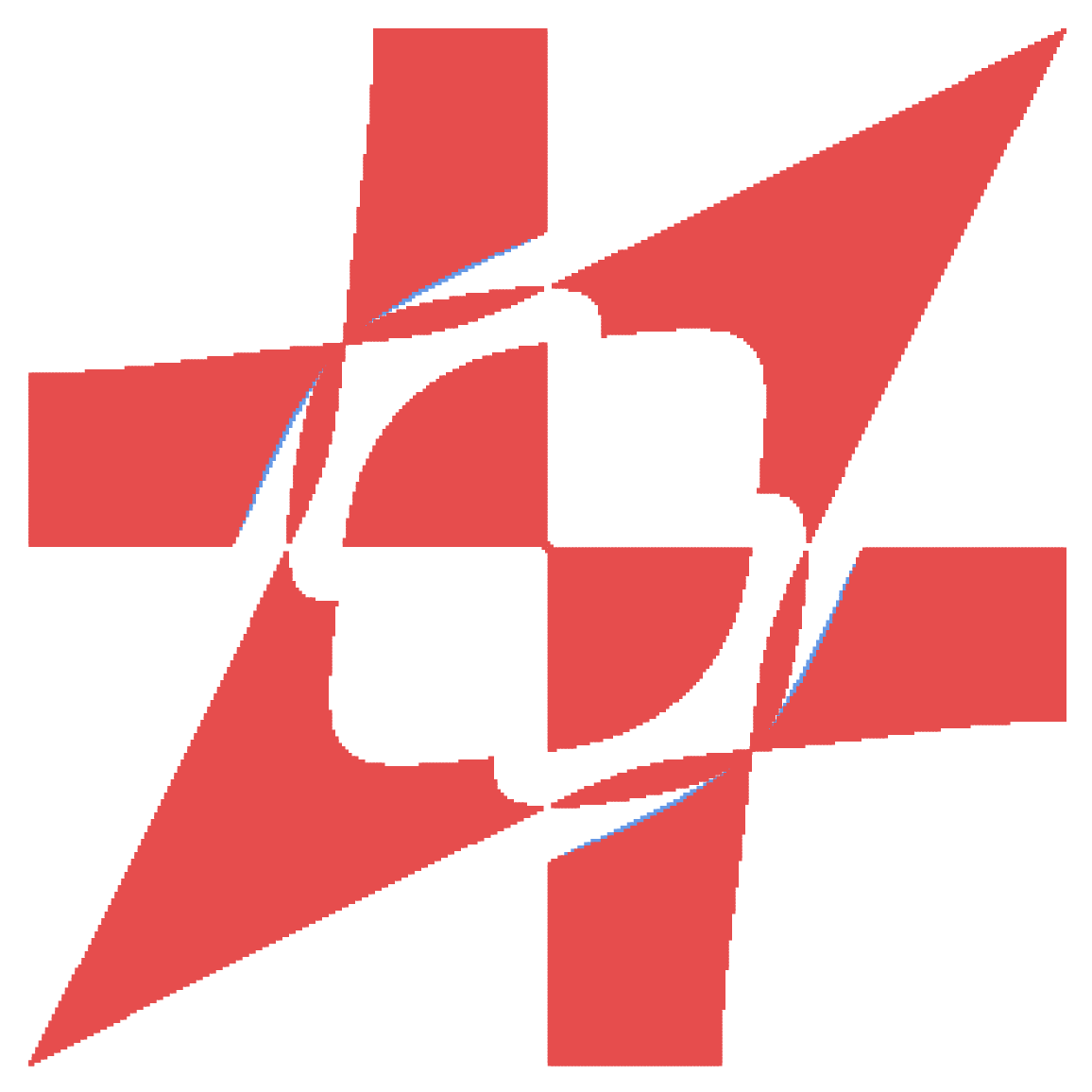}}
\rput(3.5,-1){\fbox{$n=2$}}
\rput(-0.2,-0.2){0}
\psline(0,0)(0,7) \rput(-0.3,7){$\pi$} \rput(-0.3,3.5){$\frac{\pi}2$} 
   \psline(-0.1,3.5)(0.1,3.5)  \psline(-0.1,7)(0.1,7) \rput(-0.4,6){$\theta_2$}
\psline(0,0)(7,0) \rput(7,-0.3){$\pi$} \rput(3.5,-0.4){$\frac{\pi}2$}
   \psline(3.5,-0.1)(3.5,0.1)  \psline(7,-0.1)(7,0.1) \rput(6,-0.4){$\theta_1$}
\psline[linecolor=gray,linestyle=dashed,linewidth=0.5pt](7,0)(0,7)n
\end{pspicture}}\hfill
{\begin{pspicture}(0,-1.5)(7,7.5) % n = 3 
\put(-0.1,-0.1){\includegraphics[width=2.3in]{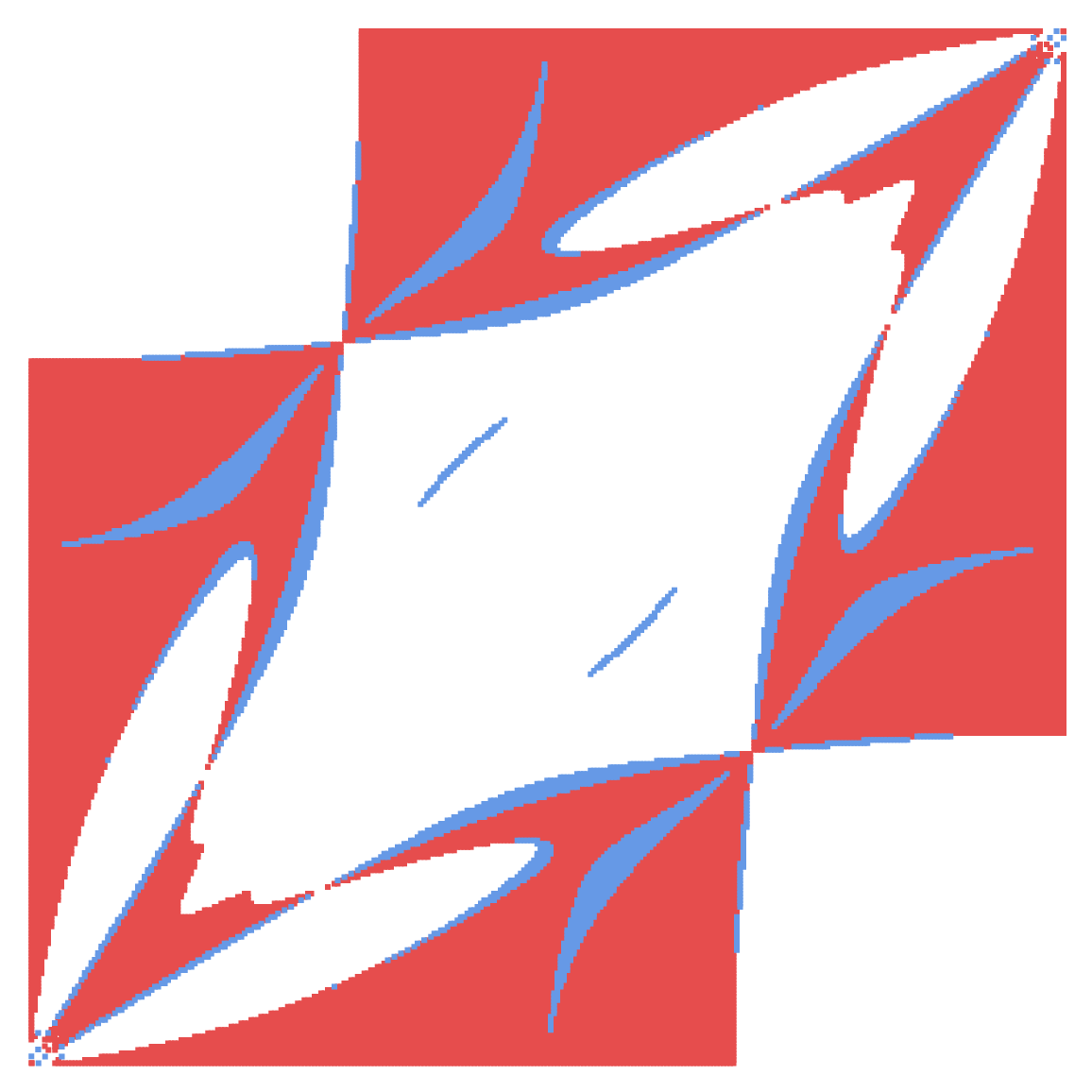}}
\rput(3.5,-1){\fbox{$n=3$}}
\rput(-0.2,-0.2){0}
\psline(0,0)(0,7) \rput(-0.3,7){$\pi$} \rput(-0.3,3.5){$\frac{\pi}2$} 
   \psline(-0.1,3.5)(0.1,3.5)  \psline(-0.1,7)(0.1,7)
\psline(0,0)(7,0) \rput(7,-0.3){$\pi$} \rput(3.5,-0.4){$\frac{\pi}2$}
   \psline(3.5,-0.1)(3.5,0.1)  \psline(7,-0.1)(7,0.1)
\psline[linecolor=gray,linestyle=dashed,linewidth=0.5pt](7,0)(0,7)n
\end{pspicture}}

{\begin{pspicture}(0,-1.5)(7,7.5) % n = 4 
\put(-0.07,-0.05){\includegraphics[width=2.3in]{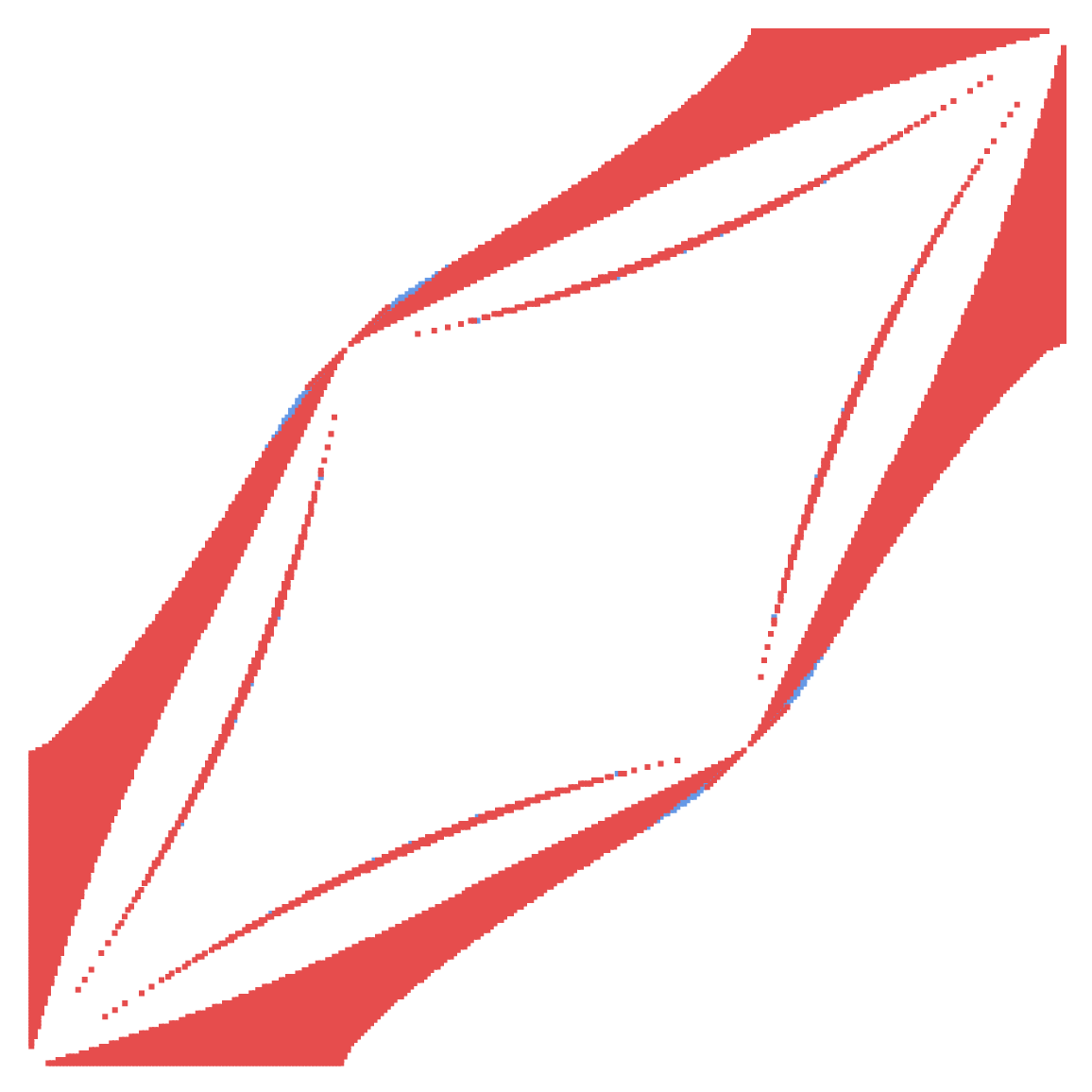}}
\rput(3.5,-1){\fbox{$n=4$}}
\rput(-0.2,-0.2){0}
\psline(0,0)(0,7) \rput(-0.3,7){$\pi$} \rput(-0.3,3.5){$\frac{\pi}2$} 
   \psline(-0.1,3.5)(0.1,3.5)  \psline(-0.1,7)(0.1,7)
\psline(0,0)(7,0) \rput(7,-0.3){$\pi$} \rput(3.5,-0.4){$\frac{\pi}2$}
   \psline(3.5,-0.1)(3.5,0.1)  \psline(7,-0.1)(7,0.1)
\psline[linecolor=gray,linestyle=dashed,linewidth=0.5pt](7,0)(0,7)n
\end{pspicture}}
\hfill
{\begin{pspicture}(0,-1.5)(7,7.5) % n = 5 (and 6)
\put(-0.07,-0.05){\includegraphics[width=2.3in]{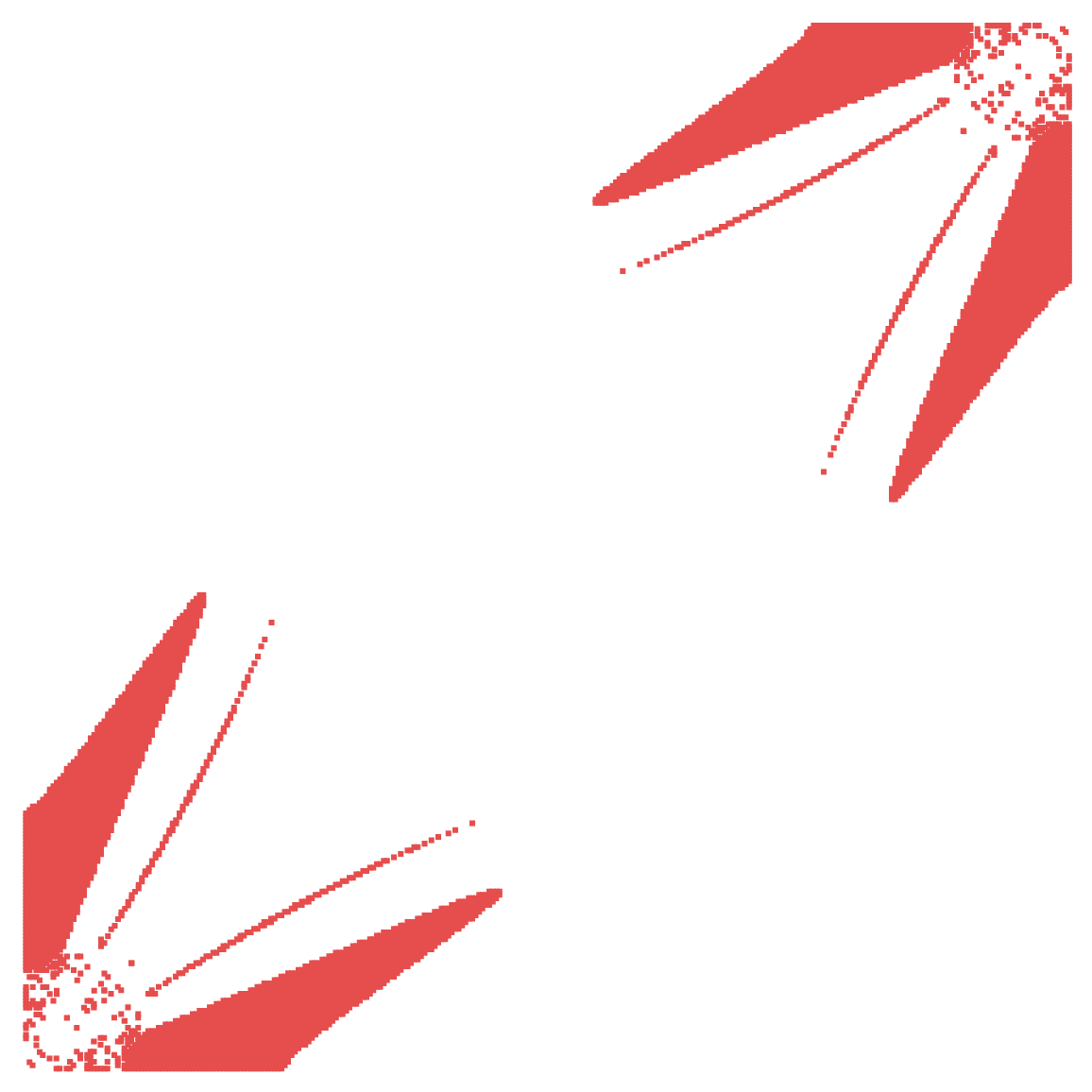}}
\rput(3.5,-1){\fbox{$n=5\;(6)$}}
\rput(-0.2,-0.2){0}
\psline(0,0)(0,7) \rput(-0.3,7){$\pi$} \rput(-0.3,3.5){$\frac{\pi}2$} 
   \psline(-0.1,3.5)(0.1,3.5)  \psline(-0.1,7)(0.1,7)
\psline(0,0)(7,0) \rput(7,-0.3){$\pi$} \rput(3.5,-0.4){$\frac{\pi}2$}
   \psline(3.5,-0.1)(3.5,0.1)  \psline(7,-0.1)(7,0.1)
\psline[linecolor=gray,linestyle=dashed,linewidth=0.5pt](3.5,0)(3.5,7)
\psline[linecolor=gray,linestyle=dashed,linewidth=0.5pt](0,3.5)(7,3.5)
\end{pspicture}}
\caption{Stability results for two staggered rings with $n$ identical vortices in each (the $\CC_{nv}(R,R')$ relative equilibria). See Figure~\ref{fig:CRp2/3} for the meaning of the colours. The diagram for $n=6$ (two rings with 6 vortices in each) is very similar to that for $n=5$, but with a smaller region of stability (reduced by about 70\%).  For $n\geq 7$, it seems that the relative equilibria are all unstable. The diagonal dashed lines represent configurations with $\theta_2=\pi-\theta_1$, where $\kappa=-1$; the  $\DD_{nd}(R,R')$ configurations of \cite{LP02}. The figures are produced using \textsc{Maple}, where the entire square is scanned with a step-size of $10^{-2}$, so nearly $10^5$ data points in each diagram.}
\label{fig:CnvRR'}
\end{center}
\end{figure}
%%%%%%%%%%%

\paragraph{Discussion of results}

\begin{itemize}
\item Referring first to Figure~\ref{fig:signkappa-staggered},  the diagonal $\theta_1=\theta_2$ corresponds to a single ring with $2n$ vortices, so must allow $\kappa=1$ for all $n$; for $n\geq3$ the two interior grey regions appear to extend to the corners of the square, so that for all $n$ there are relative equilibria close to the poles with all possible values of $\kappa$.  Each curve is labeled with a 0 or an $\infty$, corresponding to the `value' of $\kappa$ along that curve.    The degenerate case occurs where the curves $\kappa=0$ and $\kappa=\infty$ intersect. (Note that rotating the diagram by $\pi$ corresponds to turning the sphere upside down, which is a symmetry of the system and so preserves the value of $\kappa$.)

\item Numerical experiments suggest that stable relative
equilibria only exist for $n \leq 6$.
 
\item Refer now to Figure~\ref{fig:CnvRR'}. For $n=5$ and $6$ the relative equilibria $\CC_{nv}(R,R')$
are stable only if the  two rings lie in the same hemisphere but
are sufficiently far apart.

\item For $n\leq4$ these stable regions extend to include relative
equilibria with the rings in different hemispheres. However,
contrary to the case $\CC_{nv}(2R)$, the stable regions are far from
the line $\theta_2=\pi-\theta_1$ corresponding to $\DD_{nd}(R,R')$
relative equilibria.

\item For $n = 2$ and $3$ there is also a stable region with the
two rings in the same hemisphere and close to each other. This
includes the stable $\CC_{4v}(R)$ and $\CC_{6v}(R)$ relative
equilibria discussed in Section \ref{singlering}.

\item Note also that for $n\leq6$, there exist stable relative
equilibria (for some values of $\kappa$) in any neighbourhood of
$(\theta_1,\theta_2) = (0,0)$, that is with the two rings close to
a pole. 

\item A study of the sign of $\kappa$ shows that when $3 \leq n
\leq 6$ the relative equilibria with $\kappa<0$ are almost all unstable. The only ones that are not are elliptic, and occur in very narrow strips around the `hairs' towards the diagonal in the diagrams for $n=4,5,6$, and along the corresponding parts of the diagram for $n=3$, as well as the elliptic crescents near the centre of the $n=3$ diagram.

However for $n=2$ there exist relative equilibria with $\kappa>0$ in the Lyapounov stable region corresponding to the two rings both being relatively close to the equator, but in opposite hemispheres (see also the discussion below on the tetrahedral configuration).

\item It follows from Section \ref{sec:Cnv(R)} that the equatorial square is unstable if all the vorticities are equal.  However, the equatorial square with opposite vortices having the same vorticity is always an equilibrium, regardless of the two values of vorticity, and always with momentum $\mu=0$. This is the point in the centre of the first diagram in Figure~\ref{fig:CnvRR'}. The regions of stability and instability neighbouring this central point correspond to different vorticity ratios. 
Since $\mu=0$ the only relevant mode is $\ell=0$, and the corresponding Hessian is $\mathrm{diag}[16\kappa,\;-4\kappa^2]$ so the equatorial square is (linearly and non-linearly) stable if and only if $\kappa<0$. (This agrees with the conclusion in \cite[Theorem~4.6]{LP02} where the case $\kappa=-1$ is considered.)

\item It is shown by Pekarsky and Marsden \cite{PM98} that the tetrahedral configuration (where $\cos\theta_1= 1/\sqrt3= -\cos\theta_2$) is a relative equilibrium for all values of the four vorticities, though they do not discuss the stabilities. Kurakin shows that the tetrahedral equilibrium with all 4 vortices identical is Lyapounov stable \cite{Ku04}.  

In Figure~\ref{fig:CnvRR'} with $n=2$, the tetrahedral configurations lie at the two points where the two stable regions meet two unstable regions on the dashed line (but not in the centre, which corresponds to the square).  In particular here we have two pairs of identical vortices (with vorticities 1 and $\kappa$) and these are relative equilibria for all values of $\kappa$. For some values these are stable, and will have the stable regions nearby, while for others they are unstable, and will have the corresponding unstable regions nearby. Simple calculations show that they are Lyapounov stable when $|\kappa+5|>2\sqrt6$ (but $\kappa\neq0$) and linearly unstable in the remaining interval $|\kappa+5|<2\sqrt6$ (the bifurcation points are  $\kappa\approx -9.9$ and $-0.1$). 

\item The analogous point with $n=3$ corresponds to the configuration with 6  vortices lying at the vertices of an octahedron, with one ring (of unit vorticity) with $\cos\theta_1 = 1/\sqrt3$ forming one face, and the other ring (of vorticity $\kappa$) forming the opposite face of the octahedron, where $\cos\theta_2 = -1/\sqrt3$. This is a relative equilibrium for all values of $\kappa$. Kurakin \cite{Ku04} has shown this is stable when all vorticities are equal (the $\mathbb{O}(v)$ relative equilibrium of \cite{LMR01}). Calculations using \textsc{Maple} show that this is in fact Lyapounov stable if $\kappa>0$. It is linearly unstable if $|\kappa+7|<4\sqrt{3}$, otherwise it is elliptic. 

\item Some of the changes in stability as $\theta_1,\theta_2$ are varied coincide with a change of sign of $\kappa$ (for example, when $n=2$ and one ring lies on the equator). As already pointed out, when $\kappa=0$ the system is not Hamiltonian, and the methods used do not immediately apply.  To our knowledge, bifurcations involving a degeneracy of the symplectic form have not been investigated. 
\end{itemize}

%%%%%%%%%%%%%%%%%%%%%%%%%%%%%%%%%%%%%%%%%%%%%%%%%%%%%%%%%%%%%%%%%%%%%%%%%%%%%
\section{A ring with two polar vortices: $\CC_{nv}(R,2p)$}
\label{sec:R2p}

In this final section we consider relative equilibria $x_e$ of symmetry
type $\CC_{nv}(R,2p)$; that is configurations formed of a ring of
$n$ vortices of unit strength, together with two polar vortices
$p_N$ and $p_S$ of strengths $\kappa_N$ and $\kappa_S$, respectively
at the North and South poles. We may assume without loss of generality
that the ring lies in the Northern hemisphere.  There is thus a 3-parameter family of relative equilibria to consider: the parameters being $\kappa_N,\kappa_S$ and either the momentum $\mu$ or the co-latitude  $\theta_0$ of the ring.

The complexity of the study in this case is increased on two counts,
compared with the earlier setting of a ring with a single polar vortex:
firstly there are now three parameters $\kappa_N,\kappa_S$ and the
co-latitude $\theta_0$ of the relative equilibrium, and secondly (for
$n>2$) the $\ell=1$ mode of the symplectic slice is of dimension 6
rather than 4. So we need to study a 3-parameter family of 3-degree
of freedom systems. We proceed to obtain analytic (in)stability
criteria for the relative equilibria with respect to the $\ell\geq2$
modes, which of course give sufficient conditions for genuine
instability. We then treat the remaining $\ell=1$ mode numerically for a few low values of $n$. If $n=2$ or 3 then there is only the $\ell=1$ mode; the results for $n=3$ are very similar to larger rings (apart from the `cut-off' due to the higher modes when $n>3$). On the other hand, for $n=2$ the results are quite different, so we describe this last case separately at the end of this section.

\medskip

The Hamiltonian is given by
$$H = H_r + H_{p_N} + H_{p_S} + H_{NS}$$
where $H_r$ is given in Section \ref{sec:Cnv(R)} and
$$
\begin{array}{lll}
  H_{p_N}&=&-\kappa_N\sum_{i=1}^{n}\ln (1-\sin\theta_i \cos\phi_i\; x_N -\sin\theta_i \sin\phi_i\; y_N-\cos\theta_i\;z_N)\\
  H_{p_S}&=&-\kappa_S\sum_{i=1}^{n} \ln (1-\sin\theta_i \cos\phi_i\; x_S -\sin\theta_i \sin\phi_i\; y_S-\cos\theta_i\;z_S)\\
  H_{NS}&=&-\kappa_N\kappa_S\ln (1-x_N x_S - y_N y_S - z_N z_S),
\end{array}
$$
where $z_N=\sqrt{1-x_N^2-y_N^2}$ and $z_S=-\sqrt{1-x_S^2-y_S^2}$, and we have assumed the $n$ vortices in the ring are of unit vorticity. The augmented Hamiltonian is:
$$
H_{\xi}=H-\xi\left(\sum_{j=1}^n\cos\theta_j + \kappa_N\,z_N + \kappa_S\,z_S\right).
$$
The angular velocity of the relative equilibrium $x_e$ at
$x_N=y_N=x_S=y_S=0,\; \theta_j=\theta_0,\; \phi_j=2\pi j/n$ is found to be
$$
\xi = \frac{(n-1)\cos\theta_0 +
  \kappa_N(1+\cos\theta_0)-\kappa_S(1-\cos\theta_0)}{\sin^2\theta_0}.
$$
The momentum for this configuration is $(0,0,\mu)$ with $\mu= \kappa_N-\kappa_S+n\cos\theta_0$, which gives $$\xi=\frac1{\sin^2\theta_0}\left(\mu+(\kappa_N+\kappa_S-1)\cos\theta_0\right).$$

The second derivatives of $H$ at the relative equilibrium can be
derived from those for $H_r$ (Section \ref{sec:Cnv(R)}), those
for $H_p$ (Section \ref{sec:Rp-stability}), together with
\begin{equation}\label{eq:2poles hessian}
\frac{\partial^2 H_{NS}}{\partial x_N^2}=\frac{\partial^2
  H_{NS}}{\partial y_N^2}= \frac{\partial^2 H_{NS}}{\partial
  x_S^2}=\frac{\partial^2 H_{NS}}{\partial y_S^2}= \frac{\partial^2
  H_{NS}}{\partial x_N\partial x_S}=\frac{\partial^2 H_{NS}}{\partial
  y_N\partial y_S}= \frac12\kappa_N\kappa_S\,,
\end{equation}
while the other second derivatives of $H_{NS}$ vanish.

As in the previous sections, we can choose a symmetry adapted basis
of the symplectic slice such that the matrices $\dd^2H_\xi|_\NN(x_e)$
and $L_\NN$ block diagonalize. Bases of $V_\ell$ for $\ell=2,\dots
\frac{n}2$ are given in Proposition~\ref{prop:isotypic}.

\begin{proposition}
\label{prop:basis-2p}
 The symplectic slice decomposes as
$\NN = \bigoplus_{\ell=1}^{[n/2]} V_\ell,$ where for $\mu\neq0$
$$\dim V_1 =
\begin{cases}4 & \text{if\, $n=2$} \cr
       6 & \text{if\, $n\geq3$}.\end{cases}
$$
If $\mu=0$ the dimension of\/ $V_1$ is reduced by 2.

For $n\geq 3$ and $\mu\neq 0$ a basis for $V_1$ is
$\{e_1,e_2,\dots,e_6\}$, where
$$
\begin{array}{lll}
e_1 &=& \sin\theta_0\: \alpha_\theta^{(1)} + \cos\theta_0\: \beta_\phi^{(1)} \\
e_2 &=& 2 \kappa_N\, \beta_\phi^{(1)} + N \sin\theta_0\: \delta x_1\\
e_3 &=& \kappa_S\, \delta x_1 - \kappa_N\, \delta x_2 \\[4pt]
e_4 &=& \sin\theta_0\: \beta_\theta^{(1)} - \cos\theta_0\: \alpha_\phi^{(1)} \\
e_5 &=& 2 \kappa_S\, \alpha_\phi^{(1)} - N \sin\theta_0\: \delta y_2 \\
e_6 &=& \kappa_S\, \delta y_1 - \kappa_N\, \delta y_2
\end{array}
$$
With respect to the resulting basis for $\NN$, the Hessian
$\dd^2H_\xi|_\NN(x_e)$ block diagonalizes into two $3\times 3$ blocks
(for the $\ell=1$ mode) and the remainder is diagonal, while $L_\NN$
block diagonalizes into one $6\times 6$ block (for $\ell=2$), and
the remainder into $2\times 2$ blocks.

If $\mu=0$ and $n\geq3$ one can take for example $\{e_2,e_3,e_5,e_6\}$ as a basis for $V_1$.
\end{proposition}

\begin{proof}
  The proof is similar to that for a single ring (see Section
  \ref{singlering} and Proposition~\ref{diagonering}).
\end{proof}

%%%%%%%%%%%%
\subsection{The higher modes $\ell\geq2$}
\begin{figure}[t]
    \begin{center}
      \includegraphics[width=3in,height=3in]{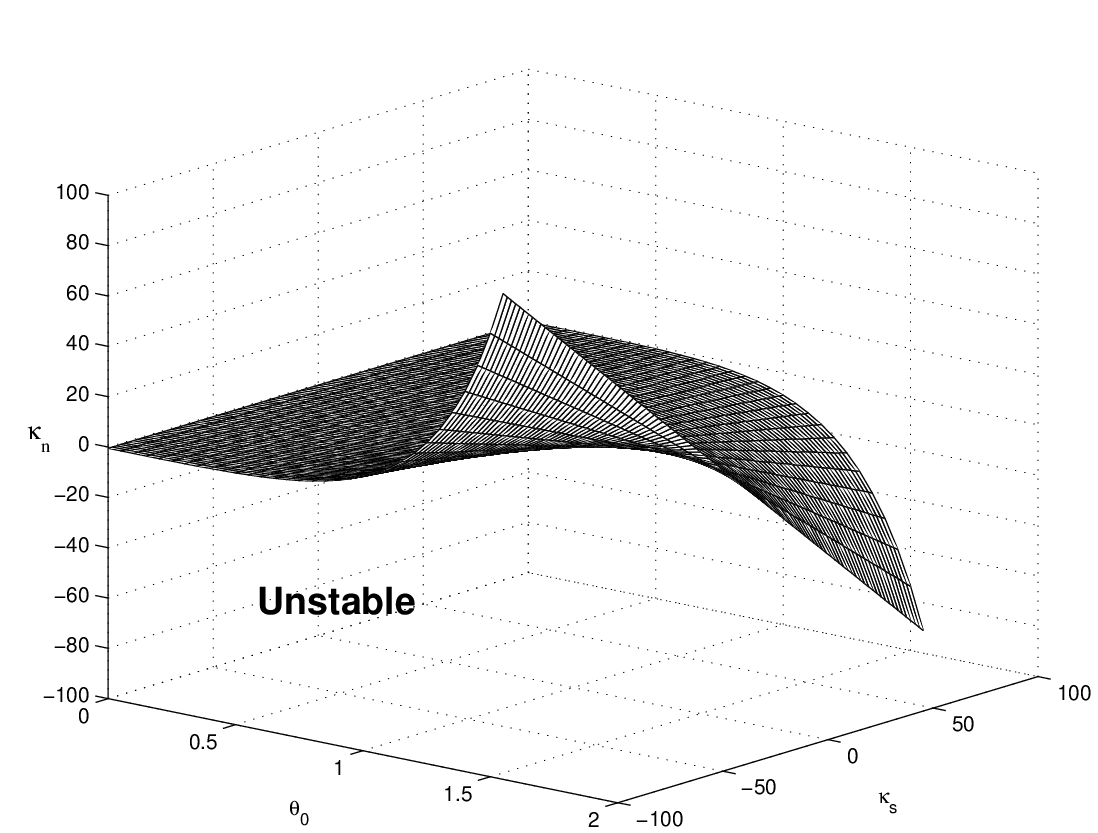}
      \caption{The relative equilibria $\CC_{nv}(R,2p)$ in the Northern
       hemisphere are unstable `below' this ruled surface in 
       $(\theta_{0},\kappa_S,\kappa_N)$-space,
       shown in the figure for $n=4$. Above the surface the relative
       equilibrium is stable with respect to all the $\ell \geq 2$ modes.}
     \label{fig:ring2poles}
     \end{center}
 \end{figure}

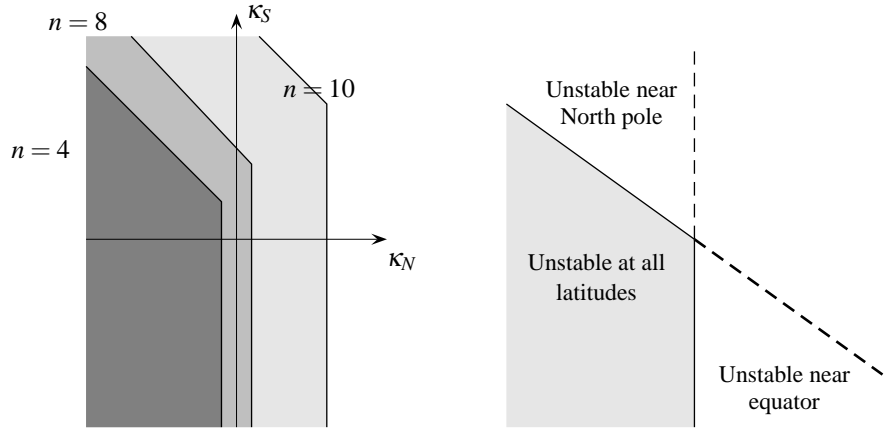
\begin{figure}[ht]
\begin{center}
\newgray{palergray}{0.9}
\begin{pspicture}(-3,-3)(3,3)
 \pspolygon[linestyle=none,fillstyle=solid,fillcolor=palergray](-2,-2.5)(-2,2.7)(0.3,2.7)(1.2,1.8)(1.2,-2.5)%n=10
 \pspolygon[linestyle=none,fillstyle=solid,fillcolor=lightgray](-2,-2.5)(-2,2.7)(-1.4,2.7)(0.2,1)(0.2,-2.5)%n=8
 \pspolygon[linestyle=none,fillstyle=solid,fillcolor=gray](-2,-2.5)(-2,2.3)(-0.2,0.5)(-0.2,-2.5)%n=4
 \psline[linewidth=0.5pt](0.3,2.7)(1.2,1.8)(1.2,-2.5)
 \psline[linewidth=0.5pt](-1.4,2.7)(0.2,1)(0.2,-2.5)
 \psline[linewidth=0.5pt](-2,2.3)(-0.2,0.5)(-0.2,-2.5)
 \psline[linewidth=0.3pt,arrowsize=.15]{->}(0,-2.5)(0,3)  \rput(0.3,3){$\kappa_S$}
 \psline[linewidth=0.3pt,arrowsize=.15]{->}(-2,0)(2,0)  \rput(2.2,-0.3){$\kappa_N$}
 \rput(-2.6,1.2){$n=4$}
 \rput(-2.1,2.9){$n=8$}
 \rput(1.1,2.){$n=10$}
\end{pspicture}
\begin{pspicture}(-3,-3)(2,3)
 \pspolygon[linestyle=none,fillstyle=solid,fillcolor=palergray](-2.5,-2.5)(-2.5,1.8)(0,0)(0,-2.5)
 \psline[linewidth=0.5pt](0,-2.5)(0,0)%(0,2.5)
 \psline[linewidth=0.5pt,linestyle=dashed](0,0)(0,2.5)
 \psline[linewidth=0.5pt](-2.5,1.8)(0,0)
 \psline[linewidth=1pt,linestyle=dashed](0,0)(2.5,-1.8)
 \rput(-1.1,2){\small Unstable near}
 \rput(-1.1,1.6){\small North pole}
 \rput(-1.3,-0.3){\small Unstable at all}
 \rput(-1.3,-0.7){\small latitudes}
 \rput(1.2,-1.8){\small Unstable near}
 \rput(1.2,-2.2){\small equator}
\end{pspicture}
\caption{Schematic diagram showing the instabilities of the $\CC_{nv}(R,2p)$ configurations due to the $\ell\geq2$ modes: (a) The shaded regions depict the values of the polar vorticities for which \emph{all} the relative equilibria in the Northern hemisphere are unstable: the darkest region represents $n=4$, the next $n=8$ and the lightest $n=10$.
  (b) demonstrates that above each shaded region of (a) the corresponding relative equilibria near the North pole are unstable, while to the right it is the relative equilibria near the equator which are unstable.}
  \label{fig:r2p-instability}
\end{center}
\end{figure}

The mode $\ell=1$ gives a $3\times 3$ block (for $\mu\neq 0$ and $n>2$) from which,
unfortunately, we can not derive a useful formula for stability analogous
to that for a single polar vortex. However, we can derive formulae
for the stability of the other modes, and thereby obtain the
following sufficient condition for instability, illustrated by
Figures \ref{fig:ring2poles} and \ref{fig:r2p-instability}. These modes occur (with the same bases) for $\mu=0$ and $\mu\neq0$ and the results here are valid in both cases, but only for $n\geq4$: for $n=2,3$ the only mode is $\ell=1$.

\begin{theorem}
\label{thm:R2p-stability} A $\CC_{nv}(R,2p)$ relative equilibrium
with $n\geq 4$ and $\mu\neq 0$ is linearly unstable if
$$
\kappa_N (1+\cos\theta_0)^2+\kappa_S (1-\cos\theta_0)^2 <
{\textstyle\left[\frac{n^2}4\right]}-(n-1)(1+\cos^2\theta_0),
$$
and is stable with respect to the $\ell \geq 2$ modes if
this inequality is reversed.
\end{theorem}

%%%%%%%%%%%%%%%%%%%%%%%%%%%%%%%%%%%%%%%%%
%%%%%%%%%%%

Thus one finds increasing either polar vorticity tends to stabilize the relative equilibrium, and if both polar vorticities are negative then the system is linearly unstable (this is only for $n\geq4$).

%%%%%%%%%%%
\begin{proof}
  The proof is similar to that for Theorem~\ref{ring+pole/stab}.
  Following the notation of the proof of that theorem, we have that the eigenvalues of the reduced Hessian are given by
  $\lambda^{(\ell)}_\phi=n\ell(n-\ell)/2$ and
  $$
  \lambda^{(\ell)}_\theta=\frac{n}{2\sin^2\theta_0}\left[
    -(\ell-1)(n-\ell-1)+(n-1)\cos^2\theta_0+\kappa_N(1+\cos\theta_0)^2+\kappa_S
    (1-\cos\theta_0)^2 \right].
  $$
  The relative equilibrium is linearly unstable if there exists
  $\ell\geq 2$, such that $\lambda^{(\ell)}_\theta<0$.  Since the
  least $\lambda^{(\ell)}_\theta$ is for $\ell=[n/2]$, the relative
  equilibrium is linearly unstable if
  $$
  -([n/2]-1)(n-[n/2]-1)+(n-1)\cos^2\theta_0+\kappa_N(1+\cos\theta_0)^2+\kappa_S
  (1-\cos\theta_0)^2 < 0,
  $$
  and is stable with respect to the $\ell \geq 2$ modes if
 this inequality is reversed. This gives the desired criterion.
\end{proof}

\paragraph{Stability of the $\ell \geq 2$ modes}
From the theorem we can deduce the following results about the
(in)stability of the $\CC_{nv}(R,2p)$ relative equilibria with
respect to the $\ell \geq 2$ modes. These modes only occur for $n
\geq 4$.  We continue to assume the ring lies in the Northern
hemisphere.

\begin{itemize}

\item In the limiting case as the ring converges to the North pole
($\theta_0=0$), for all values of $\kappa_S$ the relative equilibria
are linearly unstable if $\kappa_N <
\frac14(\left[{n^2}/4\right]-2n+2)$. This agrees with the
instability of a ring and single pole when `$\kappa < \kappa_0$' in
Proposition~\ref{ring+pole/stab}.

\item At the opposite extreme, when the ring is at the equator
($\theta_0=\pi/2$) they are linearly unstable if $\kappa_N+\kappa_S
< \left[{n^2}/4\right]-n+1$.  The right hand side of this inequality
is non-negative for all positive integers $n$, and so the
`equatorial' $\CC_{nv}(R,2p)$ relative equilibria are unstable if
the total polar vorticity has opposite sign to that of the ring. If
$\kappa_N+\kappa_S > 0$ then the critical ratio of the total polar
vorticity to the total ring vorticity needed to stabilize the $\ell
\geq 2$ modes grows linearly with $n$.

\item For all $n\geq 4$ the relative equilibria are unstable for
all latitudes in the Northern hemisphere if $\kappa_N <
\frac14(\left[{n^2}/4\right]-2n+2)$ \emph{and} $\kappa_N+\kappa_S <
\left[{n^2}/4\right]-n+1$.  In particular, for $n \geq 7$ the
relative equilibria are unstable for all $\theta_0$ if $\kappa_N <
0$ and $\kappa_S < 0$.
\end{itemize}

To determine which of these relative equilibria that are stable to the higher modes, are in  fact \emph{genuinely stable} \re\ it is necessary to evaluate the eigenvalues arising from the $\ell=1$ mode.  This we do numerically except in the special case of an equatorial configuration with polar vortices of equal strength.

%%%%%%%%%%%%%%%%%%%%
\subsection{Numerical study of the mode $\ell=1$}
\label{sec:ell=1}
\begin{figure}[tp] 
\begin{center}\psset{unit=1.5}
{\begin{pspicture}(-1.7,-2.6)(1.7,2) 
\rput(0,0){\fbox{\includegraphics[width=2in]{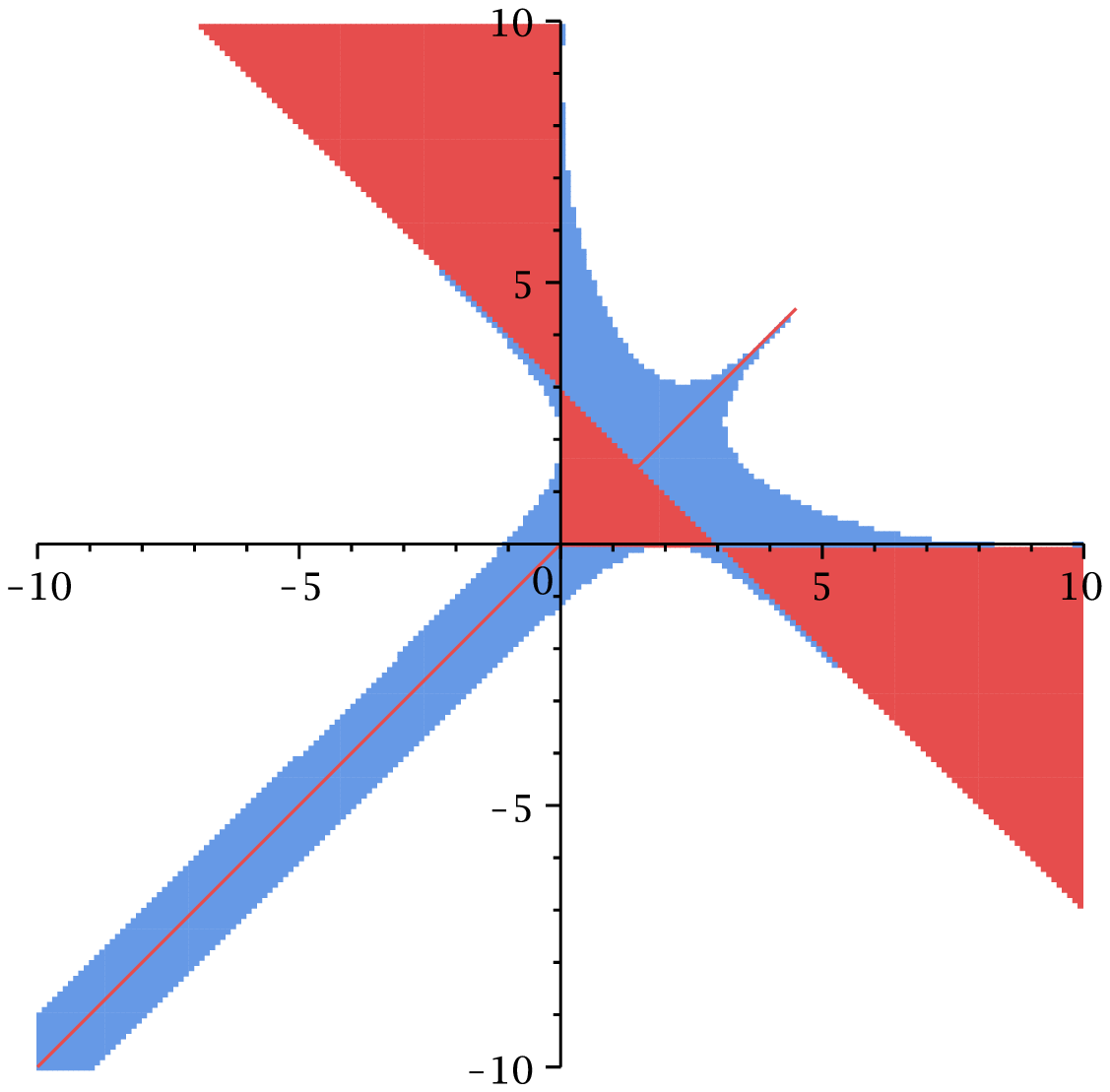}}}
\rput(0,-2.2){{$\theta_0=\pi/2$}}
\end{pspicture}}
 \qquad
{\begin{pspicture}(-1.7,-2.6)(1.7,2) 
\rput(0,0){\fbox{\includegraphics[width=2in]{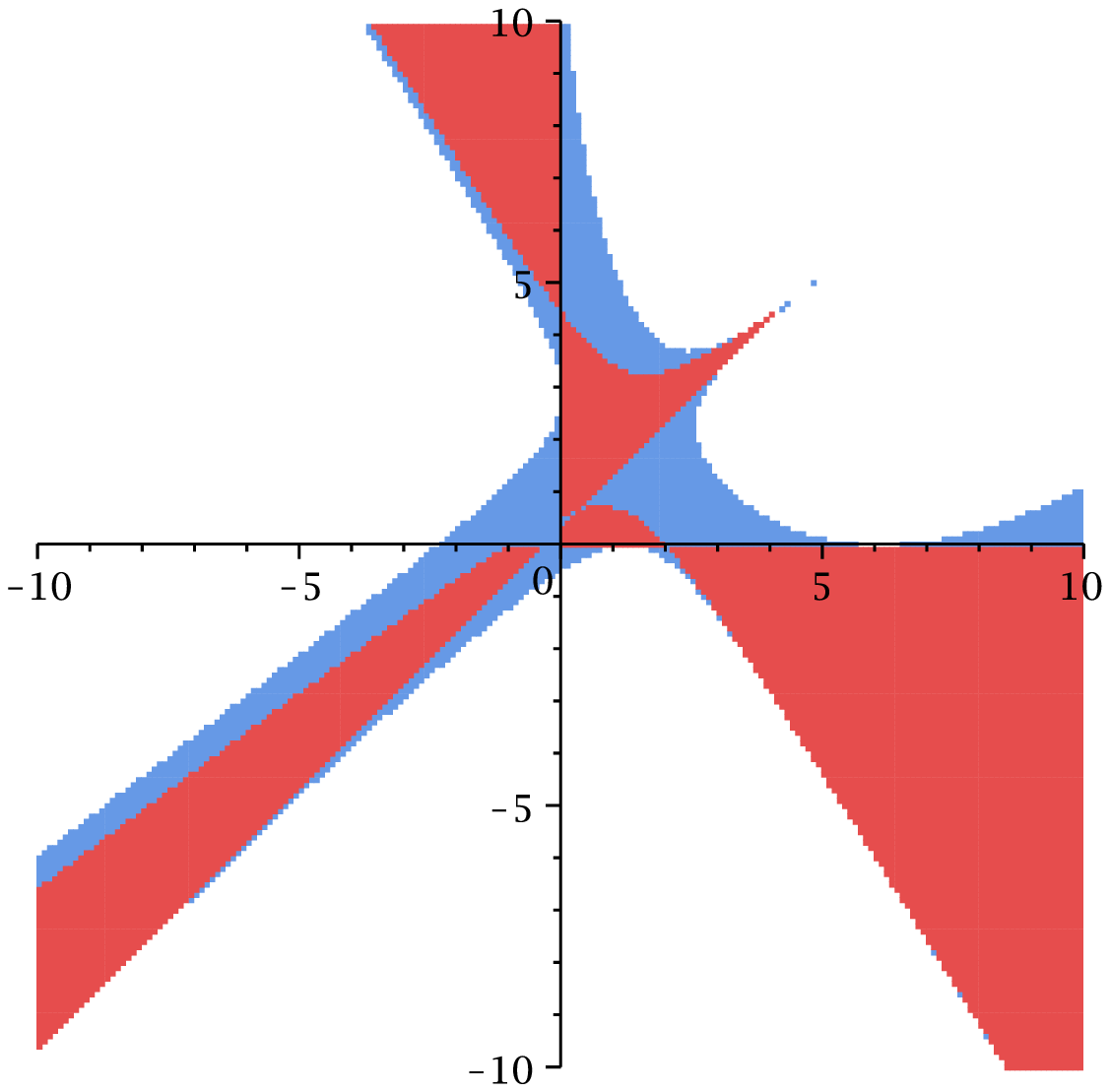}}}
\rput(0,-2.2){{$\theta_0=1.45$}}
\end{pspicture}}
{\begin{pspicture}(-1.7,-2.6)(1.7,2) 
\rput(0,0){\fbox{\includegraphics[width=2in]{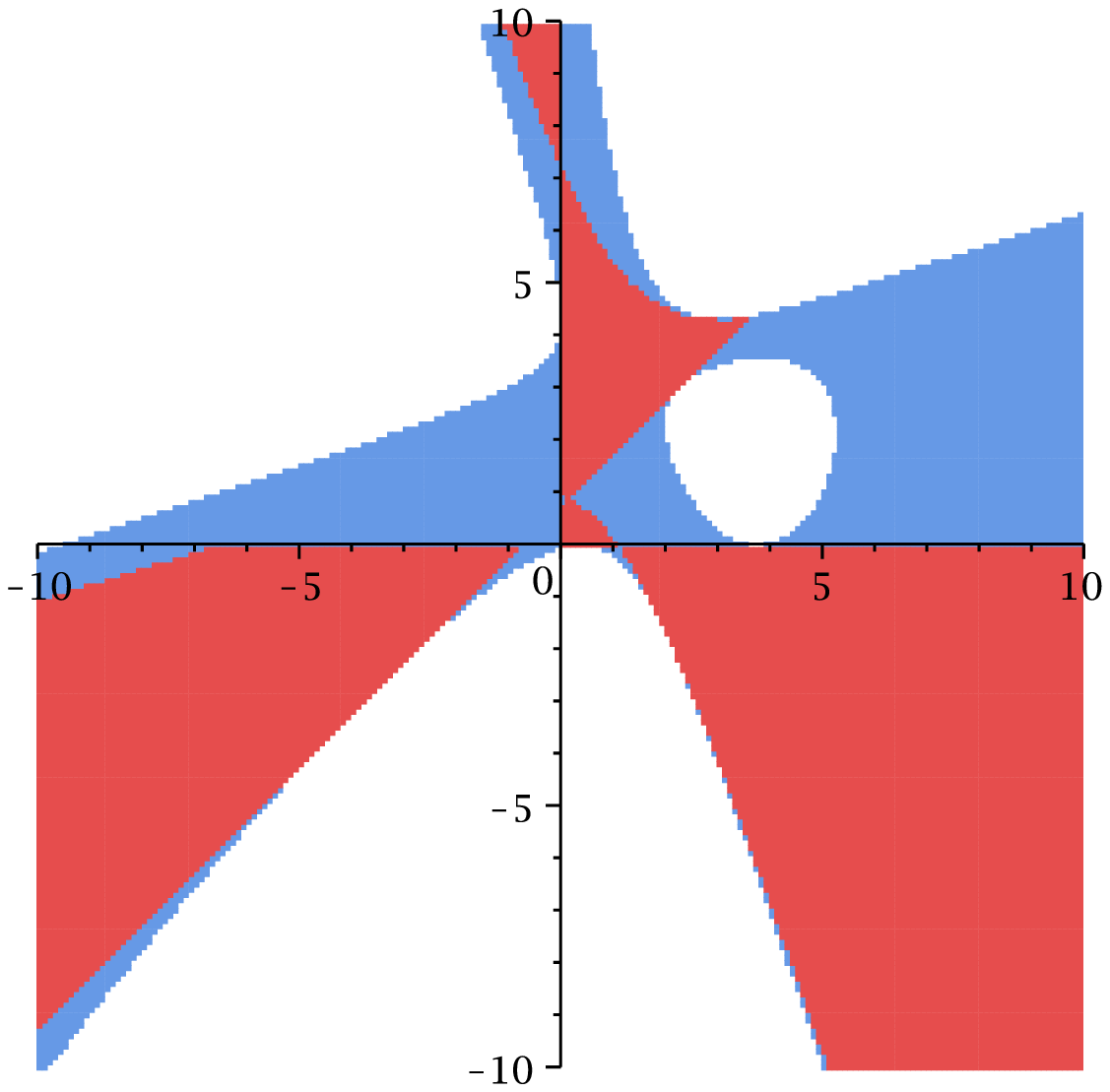}}}
\rput(0,-2.2){{$\theta_0=1.3$}}
\end{pspicture}}
 \qquad
{\begin{pspicture}(-1.7,-2.6)(1.7,2) 
\rput(0,0){\fbox{\includegraphics[width=2in]{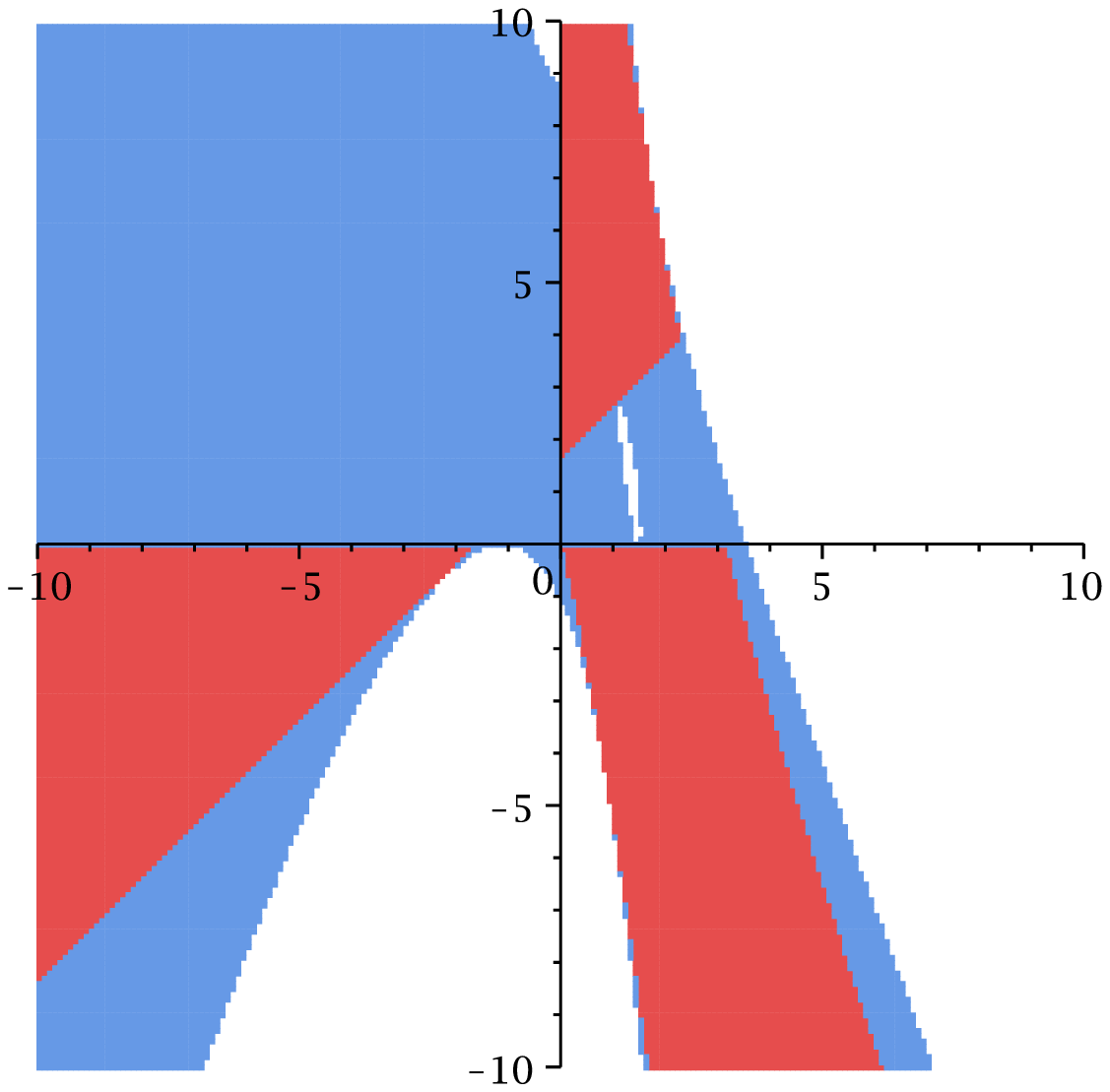}}}\rput(0,-2.2){{$\theta_0=1.0$}}
\end{pspicture}}
\caption{Stability diagrams for a 3-ring with two polar vortices --- $\CC_{3v}(R,2p)$, showing vorticities in the range  $|\kappa_N|,\,|\kappa_S|\leq10$, at four specific values of the co-latitude.  The horizontal axis is $\kappa_N$, the vertical $\kappa_S$. See Fig.~\ref{fig:CRp2/3} for the meaning of the colours.}
\label{fig:C3vR2p}
\end{center}
\end{figure}

\begin{figure}[tp] 
\begin{center}\psset{unit=1.5}
{\begin{pspicture}(-1.7,-2.6)(1.7,2) 
\rput(0,0){\fbox{\includegraphics[width=2in]{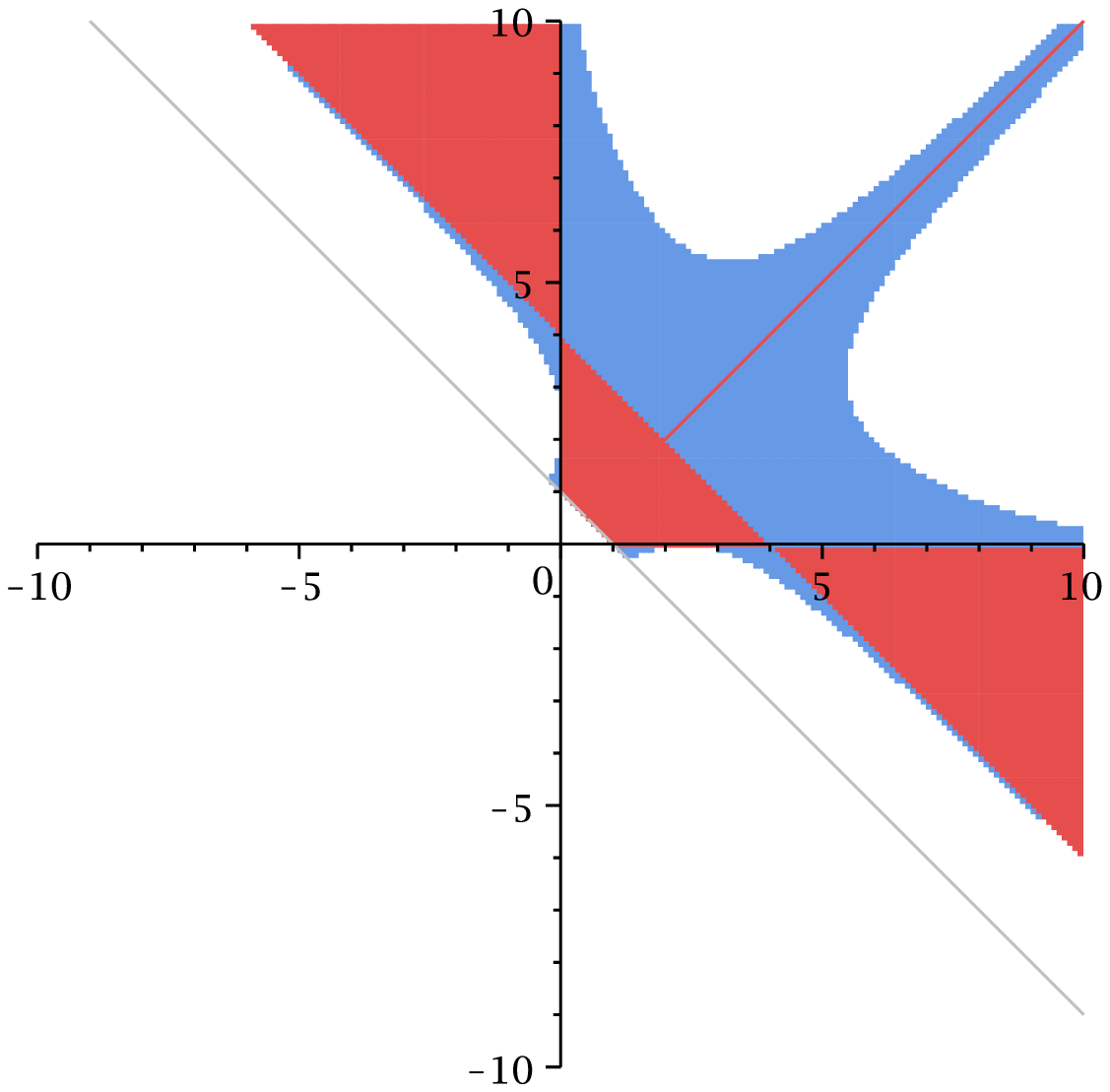}}}
\rput(0,-2.2){{$\theta_0=\pi/2$}}
\end{pspicture}}
 \qquad
{\begin{pspicture}(-1.7,-2.6)(1.7,2) 
\rput(0,0){\fbox{\includegraphics[width=2in]{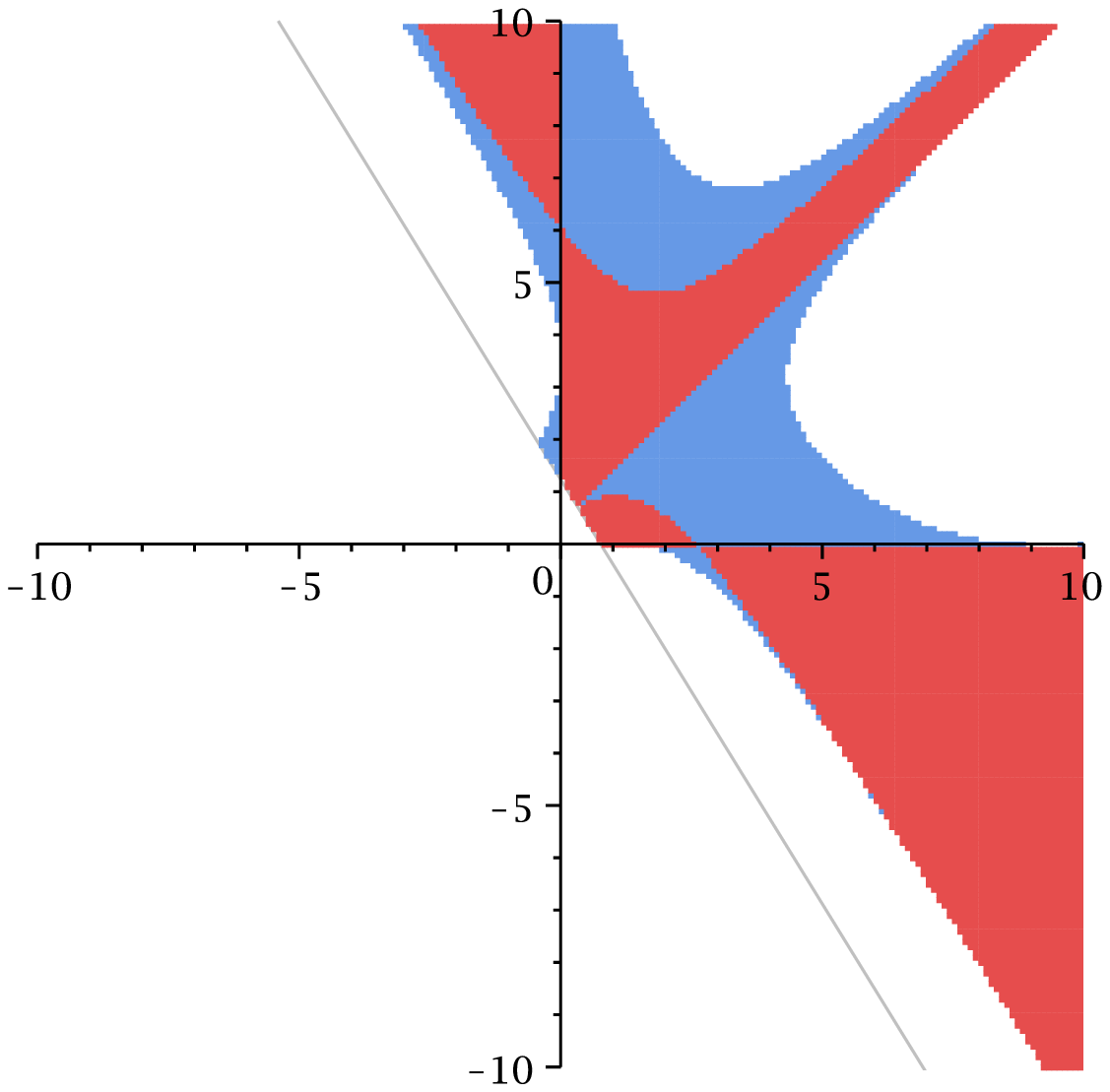}}}
\rput(0,-2.2){{$\theta_0=1.45$}}
\end{pspicture}}
{\begin{pspicture}(-1.7,-2.6)(1.7,2) 
\rput(0,0){\fbox{\includegraphics[width=2in]{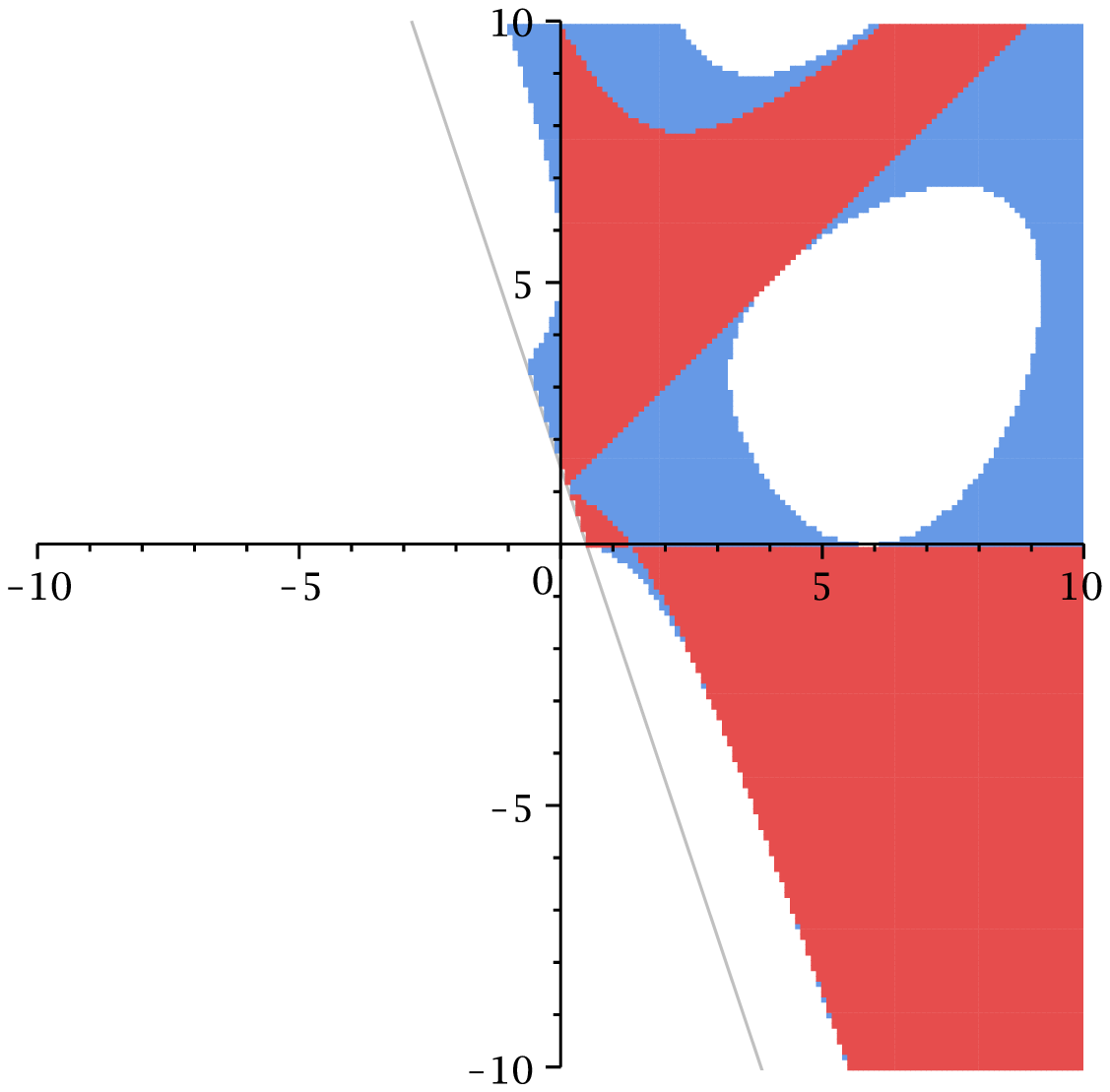}}}
\rput(0,-2.2){{$\theta_0=1.3$}}
\end{pspicture}}
 \qquad
{\begin{pspicture}(-1.7,-2.6)(1.7,2) 
\rput(0,0){\fbox{\includegraphics[width=2in]{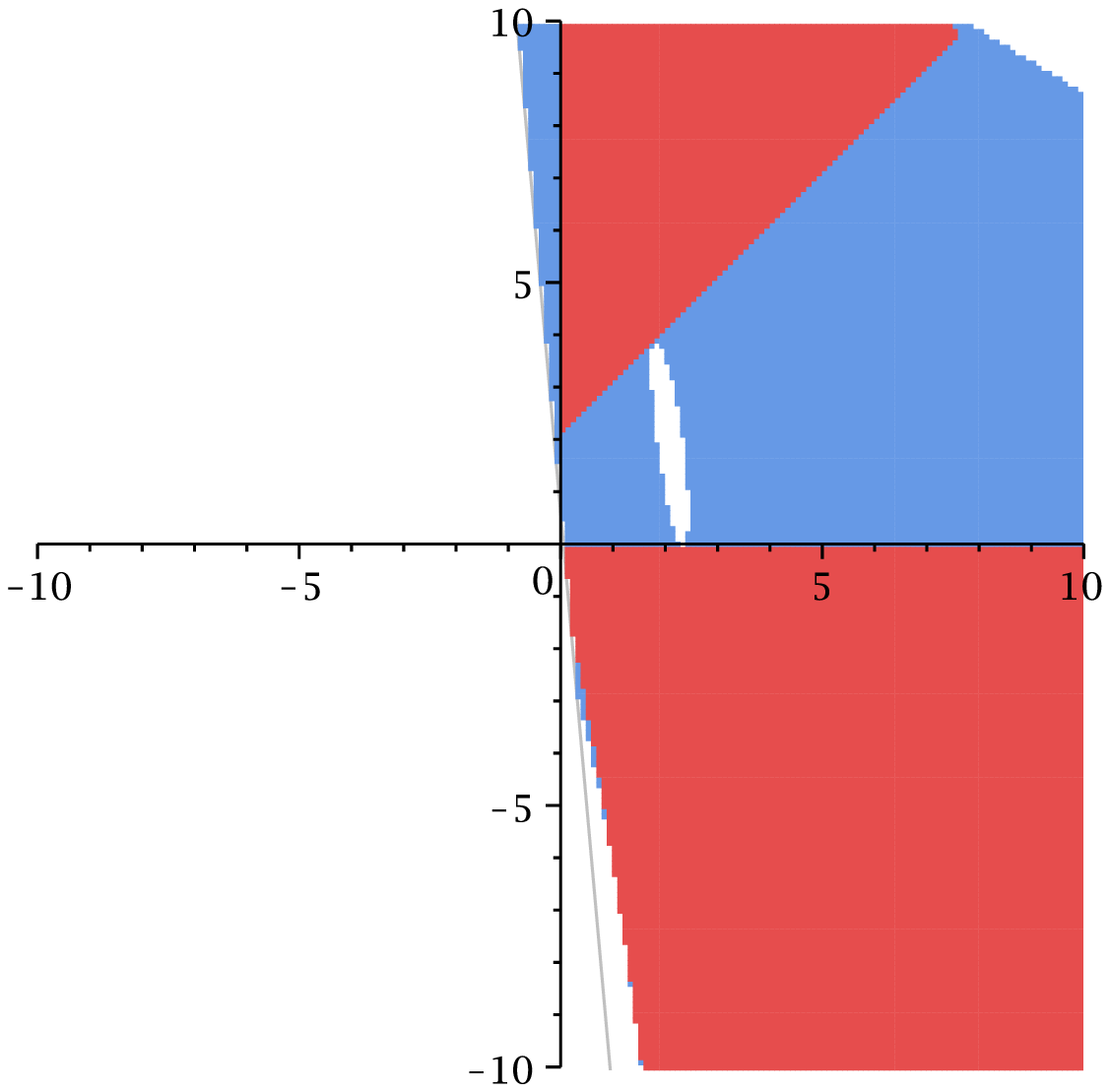}}}
\rput(0,-2.2){{$\theta_0=1.0$}}
\end{pspicture}}
\caption{Stability diagrams for a 4-ring with 2 poles (that is, $\CC_{4v}(R,2p)$), showing vorticities in the range  $|\kappa_N|,\,|\kappa_S|\leq10$, at four specific values of the co-latitude. Below and to the left of the grey line the configuration is unstable to the higher modes. The horizontal axis is $\kappa_N$, the vertical $\kappa_S$. The stability diagrams for $n>4$ are all similar (after a change in the scale of the vorticity). See Fig.~\ref{fig:CRp2/3} for the meaning of the colours.}
\label{fig:C4vR2p}
\end{center}
\end{figure}

For $n>2$ and $\mu\neq0$ the subspace $V_1$ is 6-dimensional, with basis given in Proposition~\ref{prop:basis-2p}. We are only able to obtain analytical results in the special case that the configuration is an equatorial \re\ with zero momentum. In other cases we have performed a  numerical study using \textsc{Maple}.  In this numerical study, we only investigate the `possibly stable' regime given by Theorem~\ref{thm:R2p-stability}\,; that is, we assume, for $n>3$, that
$$\kappa_N (1+\cos\theta_0)^2+\kappa_S (1-\cos\theta_0)^2 >
{\textstyle\left[\frac{n^2}4\right]}-(n-1)(1+\cos^2\theta_0).
$$
For $n=2$ and 3 there is no higher mode. There are many details here that invite  further investigation.  The case $n=2$ is rather different from the others, and we treat it in a separate section.  The principal difference between $n=3$ and $n>3$ is the presence or not of higher modes. We show the figures for $n=3$ in some detail (see Figure~\ref{fig:C3vR2p}).  The figures for $n>3$ are similar (after a rescaling of the polar vorticities), but have a cut-off given by the higher modes. Figure~\ref{fig:C4vR2p} shows the case $n=4$.

\paragraph{Equatorial \re\ with  zero momentum}
The configuration with an equatorial ring ($\theta_0=\pi/2$) has momentum zero when $\kappa_S=\kappa_N\,(=\kappa)$. Since in this case the relative equilibrium is in fact an equilibrium, a method for the full analysis of the relative equilibria in a neighbourhood can be found in \cite{MR99}, particularly Theorems 2.1 and 2.7.  The fact that the symmetry group of the relative equilibrium is $\DD_n\times\Z_2$ (where $\Z_2$ is reflexion in the equator, which acts antisymplectically, as in \cite{MR99}), the results of that paper show that, assuming generic higher order terms, for each $\mu$ close to zero there are (generically) $2n+2$ relative equilibria, two of which have $\DD_n$ symmetry and these are the points which lie on the $\CC_{nv}(R,2p)$-stratum we are considering. If the equatorial \re\ is Lyapounov stable (see below), then the two nearby ones on this stratum can be either elliptic or Lyapounov stable depending on higher order terms (and we see from the numerics they are elliptic).  

For this zero-momentum equatorial equilibrium, the $\ell=1$-mode is just 4-dimensional (for $n>2$), and can be spanned by the basis $\{e_2,e_3,e_5,e_6\}$ from the proposition above. The Hessian on this mode block diagonalizes as two identical 2x2 matrices.  

\begin{proposition}
For $n\geq3$, the eigenvalues of the linear system at the equilibrium with zero momentum on the equator are 
$$\lambda=\pm\half\sqrt{-n^2-2\kappa(n-4)\,}.$$
In particular, combining with the higher modes when $n>3$, we deduce that this configuration is Lyapounov stable under the following circumstances:
$$\begin{array}{c|c|c|c|c|c}
n=3&n=4&n=5&n=6&\dots&n=12\\
\kappa<\frac92&\kappa>\half&\kappa>1&\kappa>4&&\kappa>\frac{25}2
\end{array}$$
If the inequality is reversed, the configuration is linearly unstable.
\end{proposition}

\paragraph{Other configurations} For these we have not succeeded in obtaining analytic conditions.  The numerical calculations suggest the following:

\begin{itemize}   
\item For all $n$ and for all sufficiently large and positive polar vorticities there are ranges of $\theta_0$ with elliptic relative equilibria;
\item For all $n$ and for $\kappa_N$ sufficiently positive and $\kappa_S<0$ there are Lyapounov stable relative equilibria with the ring in the Northern hemisphere.
\item  When the equatorial ring configuration has non-zero momentum, one finds several regions of stability, shown in the first plot in each of Figures~\ref{fig:C3vR2p} and \ref{fig:C4vR2p}. In particular the reduced Hessian is degenerate if $\kappa_N+\kappa_S=n$.  
\item The straight line between the Lyapounov and elliptic regions in the two figures is where $\mu=0$. 
\item The transition between elliptic and unstable in the upper right hand part of the diagram is through a Hamiltonian Hopf bifurcation.
\item The point where the unstable region is tangent to the line $\mu=0$ (most visible for $\theta_0=1.3$) arises where there is a ro-vibrational resonance, as described in Section~\ref{sec:zero momentum}. The transitions between elliptic and unstable \re\ arbitrarily close to this point arise as Hamiltonian Hopf bifurcations, a fact which resonates with the description in \cite{Mo97}.
\end{itemize}

%%%%%%%%%%%%%%%%%%%%
\subsection{Kite configurations: $\CC_{2v}(R,2p)$}
Here we consider the remaining case $n=2$, which consists of a 2-ring with two poles: this is a kite-shaped configuration of vortices lying on a meridian, where the meridian rotates in time about the axis through the poles. The only mode in the symplectic slice is $\ell=1$, which is 4-dimensional, or 2-dimensional when $\mu=0$. It has a basis
$$\begin{array}{lclclcl}
e_1 &=&  \kappa_N\,\alpha_\theta^{(1)} -2\cos\theta_0\,\delta x_1\,, &&
e_2 &=&  \kappa_S\,\delta x_1 - \kappa_N\,\delta x_2\\[4pt]
e_3 &=&  \kappa_N\,\alpha_\phi^{(1)} - 2 \sin\theta_0\,\delta y_1,\, &&
e_4 &=&  \kappa_S\,\delta y_1 - \kappa_N\,\delta y_2\,.
\end{array}
$$
The subspaces $\left<e_1,\,e_2\right>$ and $\left<e_3,\,e_4\right>$ are both Lagrangian and invariant under the group. The reflexion $\kappa:(x,y,z)\to(x,-y,z)$ acts trivially on the first of these subspaces and by $-I$ on the second, which implies that the Hessian on this slice block diagonalizes. The entries are too long to usefully reproduce here.

When $\mu=0$, the symplectic slice is 2-dimensional, and we can use $\{e_1,e_3\}$ as a basis. Because of the invariance under the reflexion $\kappa$, the reduced Hessian on this slice is a diagonal matrix, but one obtains expressions too long to usefully reproduce here. In the special case of a momentum zero equatorial configuration, which has momentum zero if and only if $\kappa_N=\kappa_S\;(=\kappa)$, the diagonal Hessian simplifies to $\diag[4\kappa^3, -4\kappa^2]$, so the configuration is Lyapounov stable if $\kappa<0$ and linearly unstable if $\kappa>0$.

The angular velocity of the relative equilibrium is given by
$$ \xi=\frac1{ \sin^2\theta_0} \left((\kappa_N + \kappa_S +1) \cos\theta_0 +\kappa_N-\kappa_S\right).$$

\begin{figure}[tp] 
\begin{center}\psset{unit=1.5}
{\begin{pspicture}(-1.7,-2.1)(1.7,2) 
\rput(0,0){\fbox{\includegraphics[width=2in]{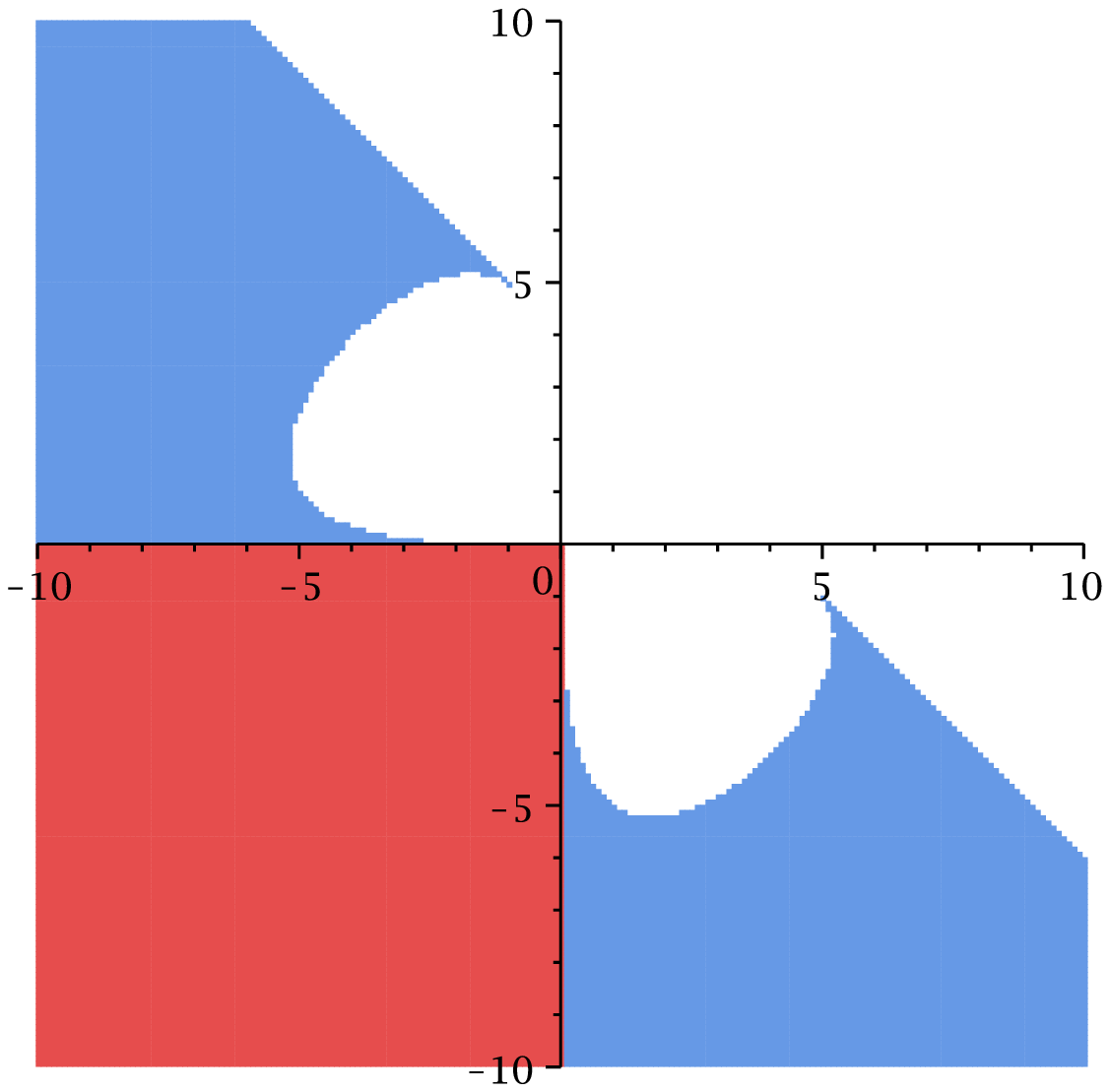}}}
\rput(0,-2){$\theta_0=\pi/2$}
\end{pspicture}}
 \qquad
{\begin{pspicture}(-1.7,-2.1)(1.7,2) 
\rput(0,0){\fbox{\includegraphics[width=2in]{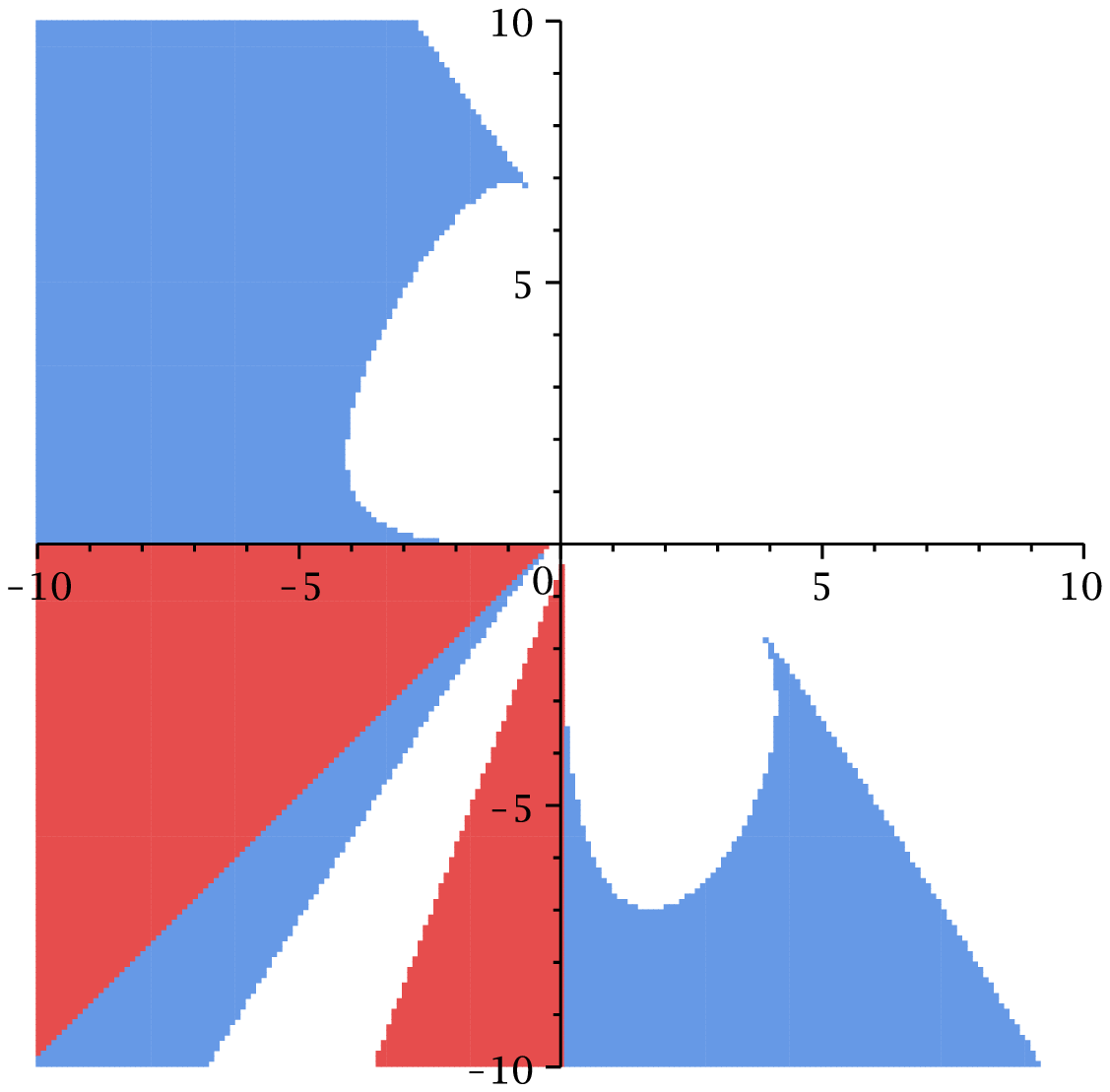}}}
\rput(0,-2){{$\theta_0=1.45$}}
\end{pspicture}}
{\begin{pspicture}(-1.7,-2.1)(1.7,2) 
\rput(0,0){\fbox{\includegraphics[width=2in]{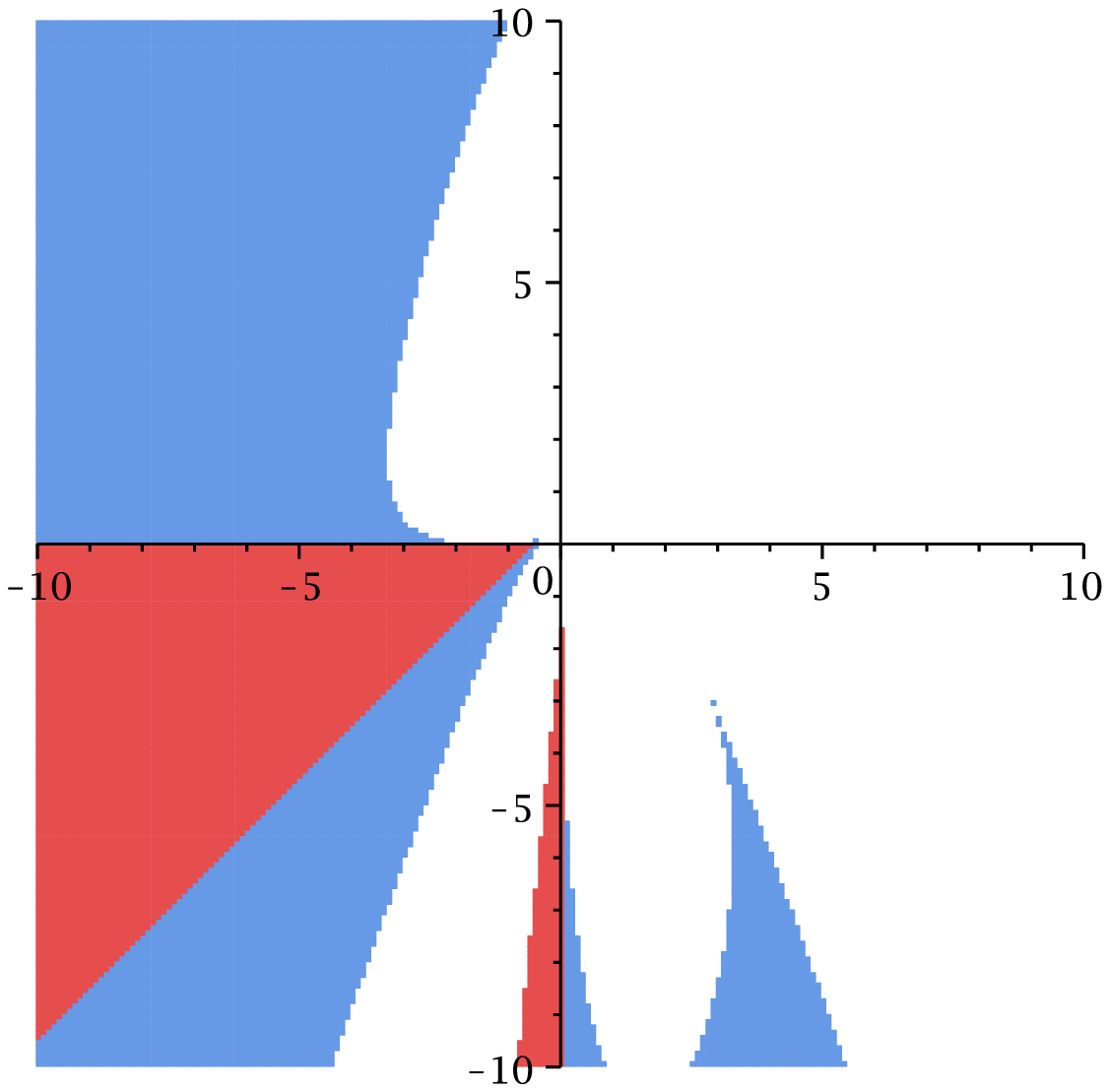}}}
\rput(0,-2){{$\theta_0=1.3$}}
\end{pspicture}}
 \qquad
{\begin{pspicture}(-1.7,-2.1)(1.7,2) 
\rput(0,0){\fbox{\includegraphics[width=2in]{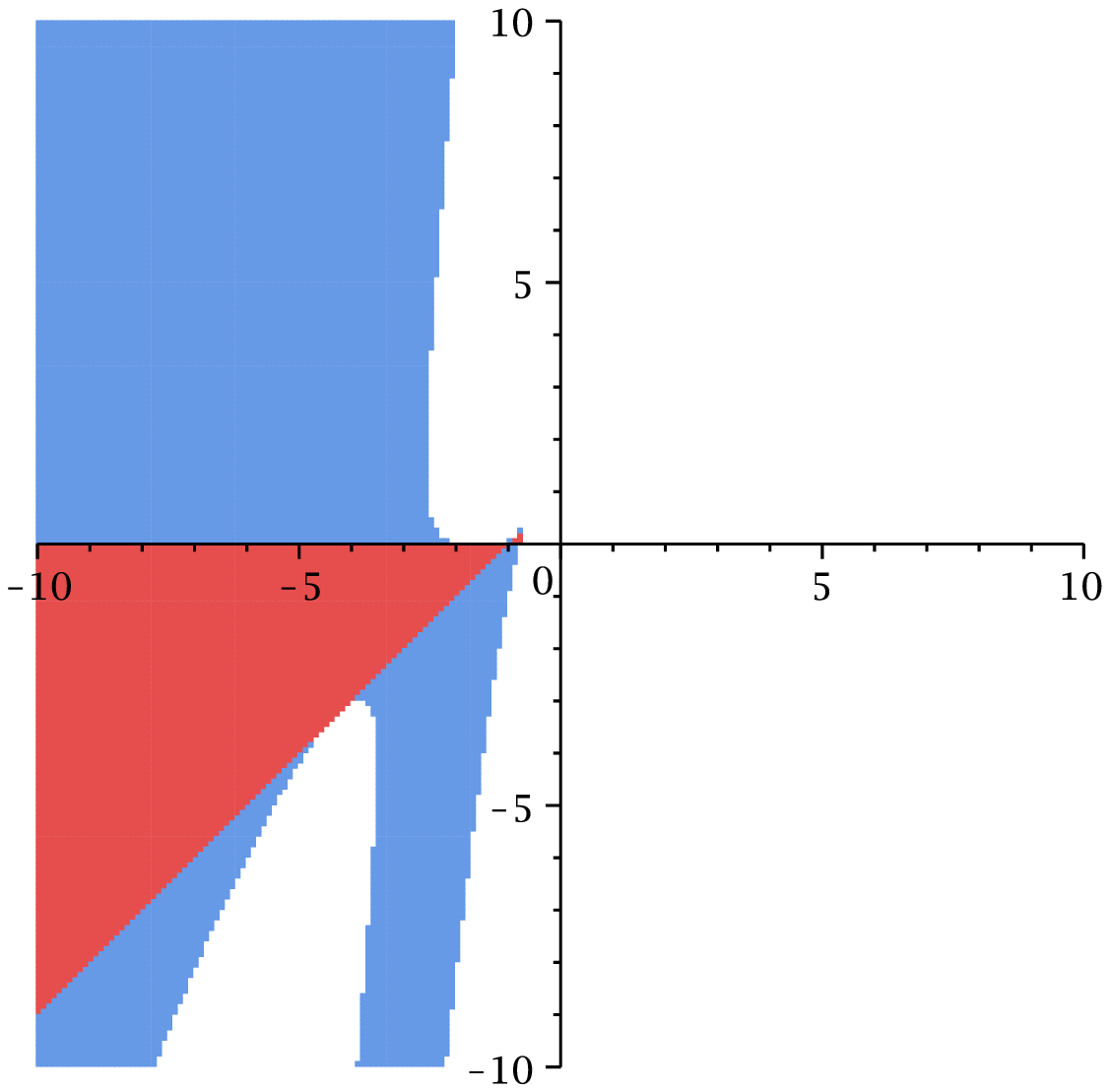}}}
\rput(0,-2){{$\theta_0=1.0$}}
\end{pspicture}}

{\begin{pspicture}(-1.7,-2.1)(1.7,2) 
\rput(0,0){\fbox{\includegraphics[width=2in]{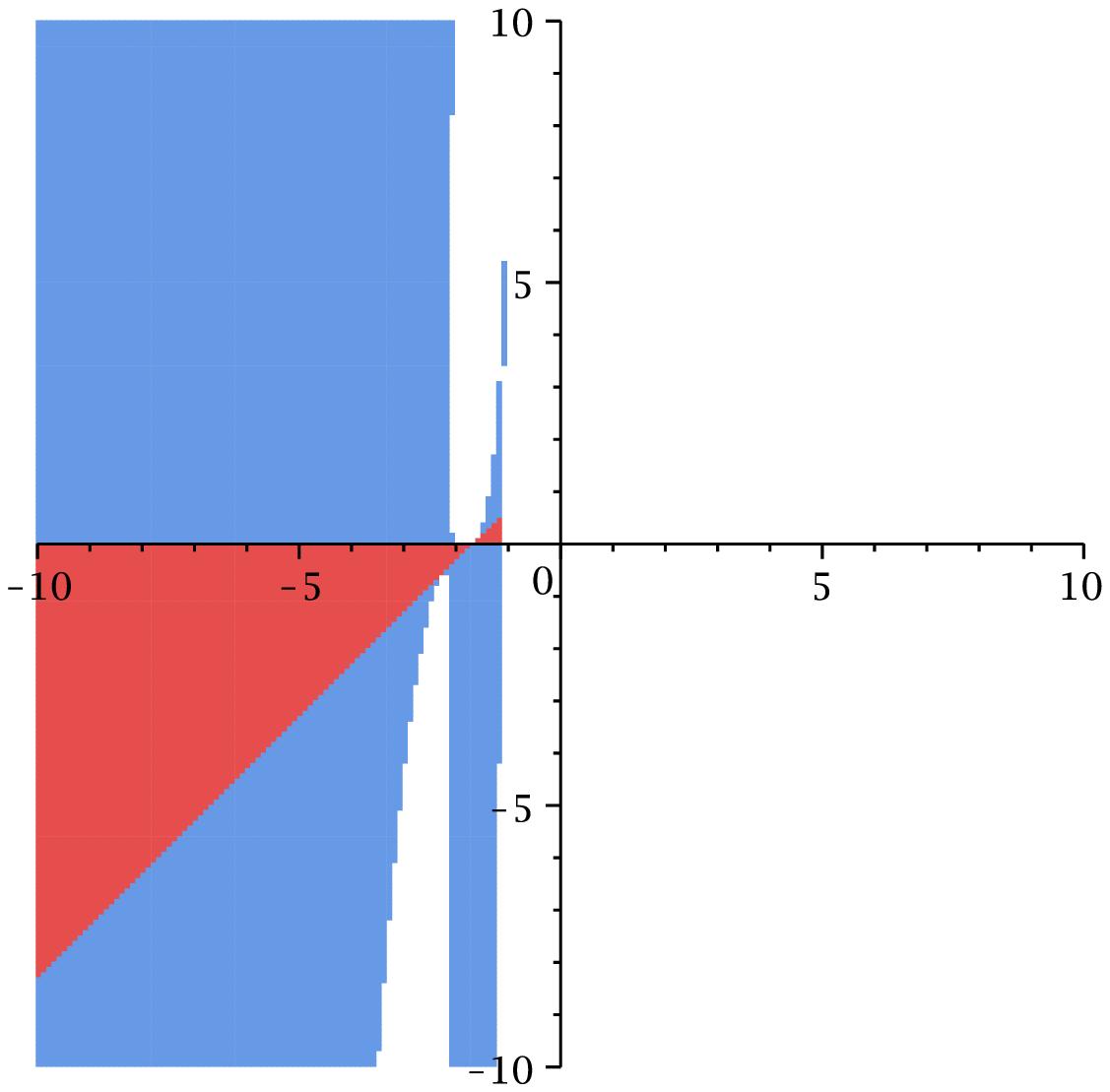}}}
\rput(0,-2){{$\theta_0=0.5$}}
\end{pspicture}}
\caption{Stability diagrams for the kite-shaped configuration $\CC_{2v}(R,2p)$. Showing $|\kappa_N|,\,|\kappa_S|\leq10$, at five specific values of the co-latitude. The horizontal axis is $\kappa_N$, the vertical $\kappa_S$. See Fig.~\ref{fig:CRp2/3} for the meaning of the colours.}
\label{fig:C2vR2p}
\end{center}
\end{figure}

\paragraph{Conclusions from numerics} See Figure~\ref{fig:C2vR2p}.
Note that with just 4 vortices, the reduced spaces have dimension 4, so any elliptic (relative) equilibria has a good chance of being Lyapounov stable, by KAM confinement. We have not checked any of the non-degeneracy conditions for this to hold, so we continue to distinguish between \re\ that we know to be Lyapounov stable by energy-momentum confinement (definiteness of the Hessian on the symplectic slice), and those we we only know to be elliptic.

\begin{itemize}
\item If both polar vorticities are positive then all configurations are linearly unstable.

\item Apart from a very small region, the only Lyapounov stable configurations occur with both polar vorticities negative. The very small region has $\kappa_S>0$ but small, and $\kappa_N<0$ (for the ring in the Northern hemisphere); it is close to the line where $\mu=0$.

\item If the ring is equatorial ($\theta_0=\pi/2$), so the configuration is a square lying on a great circle, with one pair of opposite vortices having unit vorticity, and the other two having vorticity $\kappa_N$ and $\kappa_S$, one sees that the relative equilibrium is Lyapounov stable if and only if  both $\kappa_N,\kappa_S<0$.
There is a transition across the line $\kappa_N+\kappa_S=4$, to the right of which the system has a pair of real eigenvalues and a pair of imaginary ones.  In the unstable `bulbs' near the origin, the linear system has a quadruplet of eigenvalues, and the transition from the nearby elliptic regions is via a Hamiltonian Hopf bifurcation..

\item As the ring is placed nearer the North pole, there opens up a region of instability near the negative part of the diagonal.  The straight line boundary between the Lyapounov stable and elliptic regions corresponds to $\mu=0$ (here we have $\mu=\kappa_N-\kappa_S+2\cos\theta_0$), and so the transition is that described in Section \ref{sec:zero momentum}. The transition from elliptic to unstable is through a pitchfork of splitting type, giving a pair of imaginary and a pair of real eigenvalues. Continuing the path to the right, the two real eigenvalues meet again at 0 to become imaginary, giving a Lyapounov stable configuration via another pitchfork of splitting type.

\item There is a point of interest on the $\mu=0$ locus (particularly visible in the figure for $\theta_0=1.0$) where the region of linear stability (in white) is tangent to the region of Lyapounov stability (in red). This occurs at points where there is a ro-vibrational resonance, as described in Section~\ref{sec:zero momentum} (just before the beginning of this section).  

\item In the region of linear instability shown in the bottom left quadrant of the $\theta_0=1.0$ diagram (in the $\theta_0=1.3$ diagram the corresponding region would be further out, and so is not seen), the linear system has a quadruplet of eigenvalues, and the transition from the elliptic region into this unstable one is by a Hamiltonian Hopf bifurcation, corresponding to the fact that the ro-vibrational resonance of the previous point is essentially a form of Hamiltonian-Hopf bifurcation.

\end{itemize}

%%%%%%%%%%%%%%%%%%%%%%%%%%%%%%%%%%%%%%%%%%%%
%%%%%%%%%
%%          N = 4  :
%%%%%%%%%
%%%%%%%%%%%%%%%%%%%%%%%%%%%%%%%%%%%%%%%%%%%%

\begin{figure}[h] % N=4
\begin{center}\psset{unit=0.75}
\fbox{\begin{minipage}{0.8\textwidth}
 \psset{unit=0.8}
\begin{center}
$n=4$\\
\begin{pspicture}(-1.5,-2)(1.5,2)
\rput(-4,0){$\ell=0$:}
  \pscircle(0,0){1}
  \psdots(1,0)(0,1)(-1,0)(0,-1)
  \psline{->}(1,0)(1.5,0)
  \psline{->}(0,1)(0,1.5)
  \psline{->}(-1,0)(-1.5,0)
  \psline{->}(0,-1)(0,-1.5)
  \rput(0,-1.8){$\alpha_\theta^{(0)}$}
\end{pspicture}
\begin{pspicture}(-1.5,-2)(1.5,2)
  \pscircle(0,0){1}
  \psdots(1,0)(0,1)(-1,0)(0,-1)
  \psline{->}(1,0)(1,0.5)
  \psline{->}(0,1)(-0.5,1)
  \psline{->}(-1,0)(-1,-0.5)
  \psline{->}(0,-1)(0.5,-1)
  \rput(0,-1.8){$\alpha_\phi^{(0)}$}
\end{pspicture}

\begin{pspicture}(-1.5,-2)(1.5,2)
\rput(-2,0){$\ell=1$:}
  \pscircle(0,0){1}
  \psdots(1,0)(0,1)(-1,0)(0,-1)
  \psline{->}(1,0)(1.5,0)
  \psline{->}(-1,0)(-0.5,0)
%  \psline{->}(0,-1)(0,-1.5)
  \rput(0,-1.8){$\alpha_\theta^{(1)}$}
\end{pspicture}
\begin{pspicture}(-1.5,-2)(1.5,2)
  \pscircle(0,0){1}
  \psdots(1,0)(0,1)(-1,0)(0,-1)
%  \psline{->}(1,0)(1.5,0)
  \psline{->}(0,1)(0,1.5)
%  \psline{->}(-1,0)(-1.5,0)
  \psline{->}(0,-1)(0,-0.5)
  \rput(0,-1.8){$\beta_\theta^{(1)}$}
\end{pspicture}
\begin{pspicture}(-1.5,-2)(1.5,2)
  \pscircle(0,0){1}
  \psdots(1,0)(0,1)(-1,0)(0,-1)
  \psline{->}(1,0)(1,0.5)
%  \psline{->}(0,1)(-0.5,1)
  \psline{->}(-1,0)(-1,0.5)
%  \psline{->}(0,-1)(0.5,-1)
  \rput(0,-1.8){$\alpha_\phi^{(1)}$}
\end{pspicture}
\begin{pspicture}(-1.5,-2)(1.5,2)
  \pscircle(0,0){1}
  \psdots(1,0)(0,1)(-1,0)(0,-1)
%  \psline{->}(1,0)(1,0.5)
  \psline{->}(0,1)(-0.5,1)
%  \psline{->}(-1,0)(-1,-0.5)
  \psline{->}(0,-1)(-0.5,-1)
  \rput(0,-1.8){$\beta_\phi^{(1)}$}
\end{pspicture}

\begin{pspicture}(-1.5,-2)(1.5,2)
\rput(-4,0){$\ell=2$:}
  \pscircle(0,0){1}
  \psdots(1,0)(0,1)(-1,0)(0,-1)
  \psline{->}(1,0)(1.5,0)
  \psline{->}(0,1)(0,0.5)
  \psline{->}(-1,0)(-1.5,0)
  \psline{->}(0,-1)(0,-0.5)
  \rput(0,-1.8){$\alpha_\theta^{(2)}$}
\end{pspicture}
\begin{pspicture}(-1.5,-2)(1.5,2)
  \pscircle(0,0){1}
  \psdots(1,0)(0,1)(-1,0)(0,-1)
  \psline{->}(1,0)(1,0.5)
  \psline{->}(0,1)(0.5,1)
  \psline{->}(-1,0)(-1,-0.5)
  \psline{->}(0,-1)(-0.5,-1)
  \rput(0,-1.8){$\alpha_\phi^{(2)}$}
\end{pspicture}
\end{center}
\end{minipage}}
%%%%%%%%%
%%       N = 5  :
%%%%%%%%%
\vskip5mm

\fbox{\begin{minipage}{0.8\textwidth}
\begin{center}\psset{unit=0.8}
$n=5$\\
\begin{pspicture}(-1.5,-2)(1.5,2)
\rput(-4,0){$\ell=0$:}
  \pscircle(0,0){1}
  \psdots(1,0)(0.3090, 0.9511)(-0.8090, 0.5878)(-.8090, -.5878)(.3090, -.9511)%
  \psline{->}(1,0)(1.5,0)
  \psline{->}(0.3090, 0.9511)(0.4635, 1.427)
  \psline{->}(-0.8090, 0.5878)(-1.214, 0.8817)
  \psline{->}(-.8090, -.5878)(-1.214, -0.8817)
  \psline{->}(.3090, -.9511)(0.4635, -1.427)
  \rput(0,-1.6){$\alpha_\theta^{(0)}$}
\end{pspicture}
\begin{pspicture}(-1.5,-2)(1.5,2)
  \pscircle(0,0){1}
  \psdots(1,0)(0.3090, 0.9511)(-0.8090, 0.5878)(-.8090, -.5878)(.3090, -.9511)%
   \psline{->}(1,0)(1., 0.5)
   \psline{->}(0.3090, 0.9511)(-.1669, 1.106)
   \psline{->}(-0.8090, 0.5878)(-1.103, .1834)
   \psline{->}(-.8090, -.5878)(-.5150, -.9924)
   \psline{->}(.3090, -.9511)(.7843, -.7968)
  \rput(0.1,-1.6){$\alpha_\phi^{(0)}$}
\end{pspicture}

\begin{pspicture}(-1.5,-2)(1.5,2)
\rput(-2,0){$\ell=1$:}
  \pscircle(0,0){1}
  \psdots(1,0)(0.3090, 0.9511)(-0.8090, 0.5878)(-.8090, -.5878)(.3090, -.9511)%
   \psline{->}(1,0)(1.5, 0)
   \psline{->}(0.3090, 0.9511)(.35681, 1.0980)
   \psline{->}(-0.8090, 0.5878)(-.48176, .35003)
   \psline{->}(-.8090, -.5878)(-.48176, -.35003)
   \psline{->}(.3090, -.9511)(.35681, -1.0980)
  \rput(0.1,-1.6){$\alpha_\theta^{(1)}$}
\end{pspicture}
\begin{pspicture}(-1.5,-2)(1.5,2)
  \pscircle(0,0){1}
  \psdots(1,0)(0.3090, 0.9511)(-0.8090, 0.5878)(-.8090, -.5878)(.3090, -.9511)%
   \psline{->}(1,0)(1., .5)
   \psline{->}(0.3090, 0.9511)(.16209, .99881)
   \psline{->}(-0.8090, 0.5878)(-.57126, .91505)
   \psline{->}(-.8090, -.5878)(-1.0468, -.26053)
   \psline{->}(.3090, -.9511)(.45601, -.90329)
  \rput(0.1,-1.6){$\alpha_\phi^{(1)}$}
\end{pspicture}
\begin{pspicture}(-1.5,-2)(1.5,2)
  \pscircle(0,0){1}
  \psdots(1,0)(0.3090, 0.9511)(-0.8090, 0.5878)(-.8090, -.5878)(.3090, -.9511)%
%   \psline{->}(1,0)(0,0)
   \psline{->}(0.3090, 0.9511)(.45601, 1.4033)
   \psline{->}(-0.8090, 0.5878)(-1.0468, .76054)
   \psline{->}(-.8090, -.5878)(-.57126, -.41504)
   \psline{->}(.3090, -.9511)(.16209, -.49880)
  \rput(0.1,-1.6){$\beta_\theta^{(1)}$}
\end{pspicture}
\begin{pspicture}(-1.5,-2)(1.5,2)
  \pscircle(0,0){1}
  \psdots(1,0)(0.3090, 0.9511)(-0.8090, 0.5878)(-.8090, -.5878)(.3090, -.9511)%
%   \psline{->}(1,0)(0,0)
  \psline{->}(0.3090, 0.9511)(-.14320, 1.0980)
  \psline{->}(-0.8090, 0.5878)(-.98177, .35003)
  \psline{->}(-.8090, -.5878)(-.98177, -.35003)
  \psline{->}(.3090, -.9511)(-.14320, -1.0980)
  \rput(0.1,-1.6){$\beta_\phi^{(1)}$}
\end{pspicture}

\begin{pspicture}(-1.5,-2)(1.5,2)
\rput(-2,0){$\ell=2$:}
  \pscircle(0,0){1}
  \psdots(1,0)(0.3090, 0.9511)(-0.8090, 0.5878)(-.8090, -.5878)(.3090, -.9511)%
  \psline{->}(1,0)(1.5, 0)
  \psline{->}(0.3090, 0.9511)(.18403, .56634)
  \psline{->}(-0.8090, 0.5878)(-.93404, .67862)
  \psline{->}(-.8090, -.5878)(-.93404, -.67862)
  \psline{->}(.3090, -.9511)(.18403, -.56634)
  \rput(0.1,-1.6){$\alpha_\theta^{(2)}$}
\end{pspicture}
\begin{pspicture}(-1.5,-2)(1.5,2)
  \pscircle(0,0){1}
  \psdots(1,0)(0.3090, 0.9511)(-0.8090, 0.5878)(-.8090, -.5878)(.3090, -.9511)%
  \psline{->}(1,0)(1., .5)
  \psline{->}(0.3090, 0.9511)(.69376, .82603)
  \psline{->}(-0.8090, 0.5878)(-.89985, .46277)
  \psline{->}(-.8090, -.5878)(-.71819, -.71281)
  \psline{->}(.3090, -.9511)(-0.07566, -1.0761)
  \rput(0.1,-1.6){$\alpha_\phi^{(2)}$}
\end{pspicture}
\begin{pspicture}(-1.5,-2)(1.5,2)
  \pscircle(0,0){1}
  \psdots(1,0)(0.3090, 0.9511)(-0.8090, 0.5878)(-.8090, -.5878)(.3090, -.9511)%
%   \psline{->}(1,0)(0,0)
  \psline{->}(0.3090, 0.9511)(0.39988,1.2306)
  \psline{->}(-0.8090, 0.5878)(-.42431, .30828)
  \psline{->}(-.8090, -.5878)(-1.1937, -.86730)
  \psline{->}(.3090, -.9511)(.21822, -.67154)
  \rput(0.1,-1.6){$\beta_\theta^{(2)}$}
\end{pspicture}
\begin{pspicture}(-1.5,-2)(1.5,2)
  \pscircle(0,0){1}
  \psdots(1,0)(0.3090, 0.9511)(-0.8090, 0.5878)(-.8090, -.5878)(.3090, -.9511)%
%   \psline{->}(1,0)(0,0)
   \psline{->}(0.3090, 0.9511)(0.02954, 1.0419)
  \psline{->}(-0.8090, 0.5878)(-.52951, .97250)
  \psline{->}(-.8090, -.5878)(-.52951, -0.97250)
  \psline{->}(.3090, -.9511)(0.02954, -1.0419)
  \rput(0.1,-1.6){$\beta_\phi^{(2)}$}
\end{pspicture}
\end{center}
\end{minipage}}
\end{center}
\caption{The \emph{Fourier bases} for rings of type $R$ with $n=4$ and $5$
points.}\label{fig:Fourier}
\end{figure}
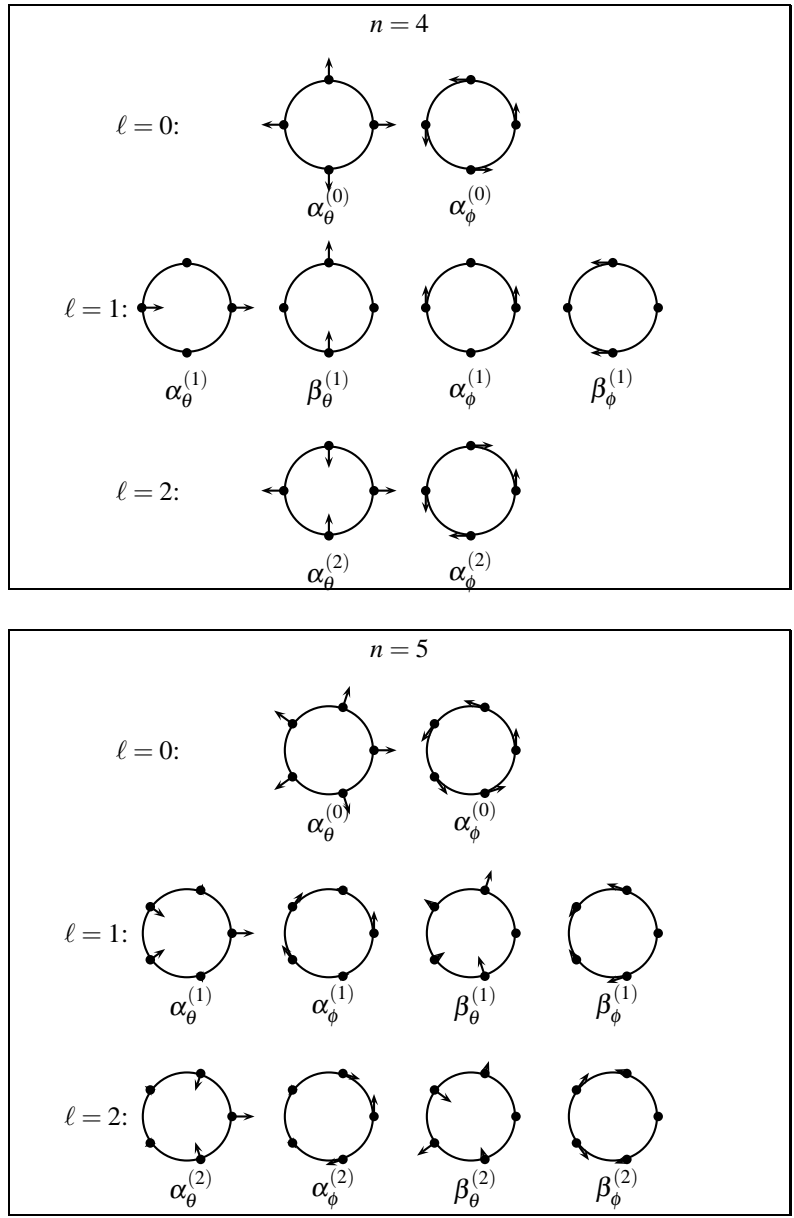

\end{document}